\pdfoutput=1
\documentclass[leqno,12pt]{article}
\usepackage{amsmath,amssymb,rotating,a4wide,color}
\usepackage{wasysym}
\usepackage{enumerate}
\usepackage{mathtools}
\usepackage[font=footnotesize,labelfont=bf,margin=2cm]{caption}

\usepackage{tocloft}
\cftsetindents{section}{0em}{2em}

\usepackage{array}
\newcolumntype{M}{>{\centering\arraybackslash}m{\dimexpr.25\linewidth-2\tabcolsep}}

\usepackage{tikz}
\usetikzlibrary{circuits.logic.US,circuits.logic.IEC,fit}
\usetikzlibrary{cd}
\usetikzlibrary{backgrounds}
\usepackage{subcaption}
\usepackage{lscape}

\usepackage[all,cmtip]{xy}

\usepackage{verbatim}

\usepackage{bbding}             
\usepackage{marginnote}         

\usepackage{hyperref}           

\newdir^{ (}{{}*!/-3pt/\dir^{(}}
\newdir_{ (}{{}*!/-3pt/\dir_{(}}

\entrymodifiers={+!!<0pt,\fontdimen22\textfont2>}

\def\bigtimes{\mathop{\raise-2pt\hbox{\huge$\times$}}}

\newbox\circbulletbox
\setbox\circbulletbox=\hbox{\,{\large$\circledcirc$}\kern-8.7pt\raise .6pt\hbox{$\bullet$}\,}

\let\le\leqslant
\let\ge\geqslant

\def\circVbig{\hbox{\text{\it\r{V}}}}
\def\circVscript{\hbox{\scriptsize\text{\it\r{V}}}}
\def\circVscriptscript{\mbox{\tiny\text{\it\r{V}}}}

\def\circVlimits_#1^#2{{\mathchoice%
{\circVbig{}^{\kern2pt #2}_{\kern-2pt #1}}%
{\circVbig{}^{\kern2pt #2}_{\kern-2pt #1}}%
{\scriptstyle\circVscript{}^{\kern1.7pt #2}_{\kern-1pt #1}}%
{\scriptscriptstyle\circVscriptscript{}^{\kern1.5pt #2}_{\kern-1pt #1}}%
}}
\def\circVr_#1{\circVlimits_#1^r}
\def\circVs_#1{\circVlimits_#1^s}

\def\circWbig{\hbox{\text{\it\r{W}}}}
\def\circWscript{\hbox{\scriptsize\text{\it\r{W}}}}
\def\circWscriptscript{\mbox{\tiny\text{\it\r{W}}}}

\def\circWlimits_#1^#2{{\mathchoice%
{\circWbig{}^{\kern2pt #2}_{\kern-2pt #1}}%
{\circWbig{}^{\kern2pt #2}_{\kern-2pt #1}}%
{\scriptstyle\circWscript{}^{\kern1.7pt #2}_{\kern-1pt #1}}%
{\scriptscriptstyle\circWscriptscript{}^{\kern1.5pt #2}_{\kern-1pt #1}}%
}}

\def\OM{\mathchoice
{\rlap{\kern3.2pt$\overline{\phantom{L}}$}M}
{\rlap{\kern3.2pt$\overline{\phantom{L}}$}M}
{\rlap{\kern2.4pt$\scriptstyle\overline{\phantom{L}}$}M}
{\rlap{\kern1.8pt$\scriptscriptstyle\overline{\phantom{L}}$}M}}

\def\mycirc{{\kern1pt\circ\kern2pt}}

\def\Mdirtext{\kern6pt\widetilde{\phantom{t}}\kern-10pt M}%
\def\Mdirscript{\kern0pt\widetilde{\phantom{N}}\kern-9pt M}%
\def\Mdirscriptscript{\kern.5pt\widetilde{\phantom{N}}\kern-7.5pt M}%
\newcommand{\Mdir}{{\mathchoice{\Mdirtext}{\Mdirtext}{\Mdirscript}{\Mdirscriptscript}}}%
\newcommand{\Mdiri}{\Mdir_{\!i}}%

\def\res{{\rm res}}

\def\cont{{\rm cont}}

\def\Div{\mathop{\rm Div}\nolimits}

\def\kalg{{k^{\rm alg}}}
\def\Kalg{{K^{\rm alg}}}
\def\Ksep{{K^{\rm sep}}}
\def\Knr{{K^{\rm nr}}}
\def\Kperf{{K^{\rm perf}}}
\def\Knrperf{{K^{\rm nr,perf}}}

\def\tame{{\rm tame}}

\def\pptauone{(\kern-1pt(\tau^{-1})\kern-1pt)}

\def\gr{{\rm gr}}

\def\Aut{\mathop{\rm Aut}\nolimits}

\def\Hom{\mathop{\rm Hom}\nolimits}

\def\Gal{\mathop{\rm Gal}\nolimits}
\def\End{\mathop{\rm End}\nolimits}

\def\coker{\mathop{\rm coker}\nolimits}

\def\Mat{\mathop{\rm Mat}\nolimits}

\def\Quot{\mathop{\rm Quot}\nolimits}

\def\rank{\mathop{\rm rank}\nolimits}

\def\ab{{\rm ab}}

\def\alg{{\rm alg}}

\def\ad{{\rm ad}}
\def\nr{{\rm nr}}
\def\perf{{\rm perf}}

\def\tor{{\rm tor}}
\def\sep{{\rm sep}}

\let\phi\varphi
\let\epsilon\varepsilon
\let\setminus\smallsetminus
\let\oldnmid\nmid
\def\nmid{\kern-1pt\oldnmid\kern-1pt}
\def\ndiv\nmid
\let\emptyset\varnothing
\def\barbar#1{\bar{\bar#1}}

\newcommand{\BF}{{\mathbb{F}}}
\newcommand{\BG}{{\mathbb{G}}}

\newcommand{\BQ}{{\mathbb{Q}}}

\newcommand{\BZ}{{\mathbb{Z}}}

\newcommand{\Ff}{{\mathfrak{f}}}

\newcommand{\Fm}{{\mathfrak{m}}}
\newcommand{\Fn}{{\mathfrak{n}}}

\newcommand{\Fp}{{\mathfrak{p}}}
\newcommand{\Fq}{{\mathfrak{q}}}
\newcommand{\Fr}{{\mathfrak{r}}}

\newcommand{\CO}{{\cal O}}
\newcommand{\CP}{{\cal P}}

\newcommand{\CV}{{\cal V}}

\newbox\mybox
\def\arrover#1{\mathrel{
\setbox\mybox=\hbox spread 1.4em
{\hfil$\scriptstyle#1$\hfil}
\vbox{\offinterlineskip\copy\mybox
\hbox to\wd\mybox{\rightarrowfill}}}}

\def\larrover#1{\mathrel{
\setbox\mybox=\hbox spread 1.4em
{\hfil$\scriptstyle#1\vphantom{g}$\hfil}
\vbox{\offinterlineskip\copy\mybox
\hbox to\wd\mybox{\leftarrowfill}}}}

\def\ontoover#1{\mathrel{
\setbox\mybox=\hbox spread 1.4em
{\hfil$\scriptstyle#1\vphantom{g}$\hfil}
\vbox{\offinterlineskip\copy\mybox
\hbox to\wd\mybox{\rightarrowfill\hskip-2.8mm
$\rightarrow$}}}}
\def\leftontoover#1{\mathrel{
\setbox\mybox=\hbox spread 1.4em
{\hfil$\scriptstyle#1\vphantom{g}$\hfil}
\vbox{\offinterlineskip\copy\mybox
\hbox to\wd\mybox{$\leftarrow$\hskip-2.8mm
\leftarrowfill}}}}
\let\longto\longrightarrow
\let\into\hookrightarrow
\let\onto\twoheadrightarrow
\def\longonto{\ontoover{\ }}
\def\longinto{\lhook\joinrel\longrightarrow}

\def\isoto{\mathrel{
\setbox\mybox=\hbox spread 0.9em
{\hfil$\scriptstyle\sim$\hfil}
\vbox{\offinterlineskip\copy\mybox
\hbox to\wd\mybox{\rightarrowfill}}}}

\def\invlim{\mathop{\vtop{\hbox{\rm lim}\vskip-8pt
\hbox{\hskip1pt$\scriptstyle\longleftarrow$}\vskip-1pt}}}

\def\Bigskip{\bigskip\bigskip}

\newtheorem{Thm}{Theorem}[section]
\newtheorem{Prop}[Thm]{Proposition}

\newtheorem{Lem}[Thm]{Lemma}

\newtheorem{Cor}[Thm]{Corollary}

\newtheorem{Def}[Thm]{Definition}

\newtheorem{Rem}[Thm]{Remark}
\newtheorem{Ex}[Thm]{Example}

\newtheorem{Cons}[Thm]{Construction}

\newtheorem{Caut}[Thm]{Caution}

\numberwithin{Thm}{section}
\def\UseTheoremCounterForNextEquation{\setcounter{equation}{\value{Thm}}\addtocounter{Thm}{1}}

\def\qed{{\hskip0pt\unskip\unskip\nobreak\hfil\penalty50
\hskip1em\hbox{}\nobreak\hfil
{$\square$}
\parfillskip=0pt\finalhyphendemerits=0
\par}\medskip}

\newenvironment{Proof}
{\noindent{\bf Proof.}}
{\qed}

{\noindent{\bf Proof of #1.}\ }
{\qed}

\newcommand{\StatementTagList}{{\upshape(\it\alph{enumi}\kern1pt\upshape)}}
\newcommand{\StatementTagText}{{(\alph{enumi})}}

\newcommand{\StatementLabels}{%
  \renewcommand{\labelenumi}{\StatementTagList}%
  \renewcommand{\theenumi}{\StatementTagText}%
}

\newcommand{\arxiv}[2]{\href{https://arxiv.org/abs/#1.#2}{\texttt{arXiv:}\hspace{0pt}\texttt{#1.}\hspace{0pt}\texttt{#2}}}

\newcommand{\isochar}{\ensuremath{\sim}}
\newcommand{\genisosign}[1]{\smash{\raisebox{-0.65ex}{#1}}}
\newcommand{\isosign}{\genisosign{\isochar}} 
\let\isoarrow\isoto
\newcommand{\longisospace}{\hspace{.375ex}}
\newcommand{\longisoarrow}{\xrightarrow{\longisospace\isosign\longisospace}}
\newcommand{\longisoto}{\longisoarrow}

\newcommand{\tholine}[1]{
  \vbox{%
    \hrule height .06em
    \kern.25ex%
    \hbox{%
      \kern-.15ex%
      \ensuremath{#1}%
      \kern-.15ex%
    }
  }
}

\newcommand{\vartholine}[1]{
  \vbox{%
    \hrule height .06em
    \kern.25ex%
    \hbox{%
      \kern-.1ex%
      \ensuremath{#1}%
    }
  }
}

\newcommand{\olsi}[2]{\hspace{#1}\overline{\hspace{-#1}{#2}}}

\newcommand{\GK}{\varGamma_{\hspace{-2pt}K}}    
\newcommand{\GKn}{\varGamma_{\hspace{-1pt}n}}   
\newcommand{\Gk}{\varGamma_{\hspace{-1pt}k}}    
\newcommand{\gal}{\gamma}                       
\newcommand{\gall}{\delta}                      

\def\Mbartext{\kern3.6pt\overline{\phantom{J}}\kern-11.48pt M}
\def\Mbarscript{\kern2.6pt\overline{\phantom{J}}\kern-8pt M}
\def\Mbarscriptscript{\kern2.3pt\overline{\phantom{J}}\kern-6.7pt M}
\newcommand{\Mbar}{{\mathchoice{\Mbartext}{\Mbartext}{\Mbarscript}{\Mbarscriptscript}}}

\def\mbar{{\kern2.5pt\overline{\phantom{n}}\kern-8.5pt m}}

\def\Mtildetext{\kern5.8pt\widetilde{\phantom{t}}\kern-10pt M}
\def\Mtildescript{\kern1.4pt\widetilde{\phantom{N}}\kern-9pt M}
\def\Mtildescriptscript{\kern1.1pt\widetilde{\phantom{N}}\kern-7.5pt M}
\newcommand{\Mtilde}{{\mathchoice{\Mtildetext}{\Mtildetext}{\Mtildescript}{\Mtildescriptscript}}}

\newcommand{\kum}[3]{[\hspace{1.5pt}#1,#2\hspace{.5pt})_{#3}} 

\newcommand{\complot}{\mathbin{\widehat{\otimes}}} 

\newcommand{\WithThickBar}[1]{
  \mathchoice{\tholine{#1}}{\tholine{#1}}%
  {\!\scalebox{.65}{\tholine{#1}}\!}%
  {\!\scalebox{.65}{\tholine{#1}}\!}%
}

\newcommand{\WithThickBarVar}[1]{
  \mathchoice{\vartholine{#1}}{\vartholine{#1}}%
  {\!\scalebox{.65}{\vartholine{#1}}\!}%
  {\!\scalebox{.65}{\vartholine{#1}}\!}%
}

\newcommand{\Fpbar}{{\WithThickBar{\Fp}}}  
\newcommand{\Fqbar}{{\WithThickBar{\Fq}}}  
\newcommand{\Frbar}{{\WithThickBarVar{\Fr}}}  

\newcommand{\Fpres}{\Fpbar}              
\newcommand{\Fqres}{\Fqbar}              
\newcommand{\Frres}{\Frbar}              

\newcommand{\Tu}{T^\circ}                

\newcommand{\Tuad}{T^{\kern.5pt\circ}{\hspace{-1.3ex}}_{\ad}}

\newcommand{\Tux}[1]{T^{\kern.5pt\circ}{\hspace{-1.2ex}}_{#1\,}}
\newcommand{\Tup}{\Tux{\Fp}}
\newcommand{\Tuq}{\Tux{\Fq}}

\newcommand{\BFp}{\BF_{\!p}}             

\newcommand{\BDM}{\olsi{1pt}{\omega}}    

\newcommand{\Bado}{{B^{\kern.5pt\circ}{\hspace{-1ex}}_{\ad}}}

\newcommand{\Rado}{{R^{\kern1pt\circ}{\hspace{-.9ex}}_{\ad}}}

\newcommand{\Bop}{{B^{\kern.5pt\circ}{\hspace{-1ex}}_{\Fp}}}

\newcommand{\KOKsym}[1]{\CP_{\kern-1pt#1}}
\newcommand{\KOK}{\KOKsym{K}}
\newcommand{\brKOK}{\KOK} 
\newcommand{\KOKprime}{\KOKsym{K'}}

\newcommand{\KOKperf}{\KOKsym{\Kperf}}

\newcommand{\Wn}[2]{W_{\kern-1pt#1}\hspace{1pt}#2}%
\newcommand{\Wo}{\Wn{\bullet}}%
\newcommand{\Wi}[1]{\Wn{i}{#1}}%
\newcommand{\Wii}[1]{\Wn{i+1}{#1}}%
\newcommand{\Wj}[1]{\Wn{j}{#1}}%
\newcommand{\grWi}{\gr^W_{i}\!}%

\begin{document}

\title{\strut
\vskip-80pt
Local Kummer theory for Drinfeld modules
}
\author{
\begin{minipage}{.3\hsize}
Maxim Mornev\\[12pt]
\small Institute of Mathematics \\
EPFL\\
1015 Lausanne\\
Switzerland \\
maxim.mornev@epfl.ch\\[9pt]
\end{minipage}
\qquad
\begin{minipage}{.3\hsize}
Richard Pink\\[12pt]
\small Department of Mathematics \\
ETH Z\"urich\\
8092 Z\"urich\\
Switzerland \\
pink@math.ethz.ch\\[9pt]
\end{minipage}
}
\date{\today}

\maketitle

\Bigskip

\begin{abstract}\noindent
Let $\phi$ be a Drinfeld $A$-module of finite residual characteristic $\Fpres$ over a local field~$K$. We study the action of the inertia group of $K$ on a modified adelic Tate module $\Tuad(\phi)$ which  differs from the usual adelic Tate module only at the $\Fpres$-primary component. After replacing $K$ by a finite extension we can assume that $\phi$ is the analytic quotient of a Drinfeld module $\psi$ of good reduction by a lattice $M\subset K$. The image of inertia acting on $\Tuad(\phi)$ is then naturally a subgroup of $\Hom_A(M,\Tuad(\psi))$.

This subgroup is described by a canonical local Kummer pairing that we study extensively in this article. In particular we give an effective formula for the image of inertia up to finite index, and obtain a necessary and sufficient condition for this image to be open. We also determine the image of the ramification filtration.
\end{abstract}

{\renewcommand{\thefootnote}{}
}

\newpage
\tableofcontents

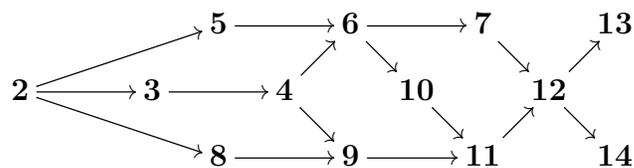
\Bigskip
$$\xymatrix@!R=3pt@!C=1pt{
&&& \bf5 \ar[rr] && \bf6 \ar[rr] \ar[dr] && \bf7 \ar[dr] && \bf13 &\\
\bf2 \ar[rr] \ar[drrr] \ar[urrr] && \bf3 \ar[rr] && \bf4 \ar[ur] \ar[dr] &&\bf10\ar[dr]&& \bf12 \ar[ur] \ar[dr] &&\\
&&& \bf8 \ar[rr] && \bf9 \ar[rr] && \bf11 \ar[ur] && \bf14 &\\}$$
\vskip1em\centerline{\bf Leitfaden}


\newpage
\section{Introduction}
\label{Intro}

Galois representations associated to Drinfeld modules play the same role for the arithmetic of function fields as Galois representations associated to elliptic curves and abelian varieties play for the arithmetic of number fields. Among the many common features of both situations is the fact that the knowledge of local Galois groups is a vital ingredient in the study of global Galois groups. Another lies in the different ways in which local Galois groups appear. Restricting ourselves to Tate modules at primes different from the residual characteristic, the action is described on the one hand by the image of Frobenius on the good reduction part, and on the other hand by the unipotent action of inertia in the case of bad reduction. In the number field case the latter comes from the tame inertia  group, which is isomorphic to $\prod_{\ell\ne p}\BZ_\ell(1)$, 
and whose action is closely linked to the reduction behavior of the abelian variety.
The aim of the present paper is to study the corresponding unipotent action of inertia in the Drinfeld module case.

%

\medskip
So let $\phi$ be a Drinfeld $A$-module over a local field~$K$ of characteristic~$p$, 
such that for all elements $a \in A$ the constant coefficient of $\phi_a$ belongs to the valuation ring~$\CO_K$. 
Being interested only in the image of Galois up to finite index, we may replace $K$ by a finite extension and $\phi$ by an isomorphic Drinfeld module. 
We may therefore assume that $\phi$ has stable reduction, i.e., that $\phi$ has coefficients in the ring $\CO_K$ and that its reduction $\bar\phi$ modulo the maximal ideal $\Fm_K$ is a Drinfeld $A$-module of positive rank over the finite residue field $k:=\CO_K/\Fm_K$. By Tate uniformization \cite[\S7]{DrinfeldEll} there then exist a unique Drinfeld $A$-module $\psi$ of good reduction over~$\CO_K$ and a unique finitely generated projective $A$-submodule $M\subset \Kalg$ for the action of $A$ through~$\psi$, such that $M\cap\CO_\Kalg=\{0\}$ and that $\phi$ is the analytic quotient of $\psi$ by~$M$. After extending $K$ again we may and do assume that $M\subset K$. 

\medskip
Recall that the adelic Tate module $T_\ad(\phi)$ is constructed from the $A$-module of all $\phi$-torsion points in~$\Kalg$. We define the \emph{modified adelic Tate module} $\Tuad(\phi)$ by using all $\phi$-torsion points in $\Kalg/\Fm_\Kalg$ instead. Like the usual Tate module this is the product of the modified $\Fp$-adic Tate modules $\Tup(\phi)$ for all primes $\Fp$ of~$A$.
For $\Fp$ different from the characteristic ideal $\Fpres$ of $\bar\phi$ this factor coincides with the usual Tate module $T_\Fp(\phi)$, while for $\Fp=\Fpres$ it is a quotient thereof.
The aim of this article is to determine the image of the Galois group $\GK = \Gal(\Ksep/K)$
up to commensurability in the action on $\Tuad(\phi)$.

\medskip
As a first step the description of $\phi$ as the analytic quotient of $\psi$ by~$M$ implies a natural isomorphism between $\Tuad(\phi)$ and the modified extended adelic Tate module $\Tuad(\psi,M)$ associated to $\psi$ and~$M$. Since this isomorphism is $\GK$-equivariant, it suffices to determine the image on the latter.

\medskip
As a second step the modified Tate module lies in a natural short exact sequence
\UseTheoremCounterForNextEquation
\begin{equation}\label{TadMIntro}
\xymatrix{ 0 \ar[r] & \Tuad(\psi) \ar[r] & \Tuad(\psi,M) \ar[r] & M\otimes_AA_\ad \ar[r] & 0\\}
\end{equation}
with $A_\ad := \prod_\Fp A_\Fp$, which is also $\GK$-equivariant. Since $\psi$ has good reduction $\bar\phi$ modulo~$\Fm_K$, there is a natural isomorphism $\Tuad(\psi) \cong \Tuad(\bar\phi)$. In particular the action on it is unramified and hence completely determined by the image of the Frobenius of~$k$. On the other hand the fact that $M\subset K$ implies that $\GK$ acts trivially on $M\otimes_AA_\ad$. To determine the image of $\GK$ it thus remains to determine the image of the inertia group under the induced homomorphism
\UseTheoremCounterForNextEquation
\begin{equation}\label{InertiaTadMIntro}
I_K := \Gal(K^\sep/K^\nr) \longto \Hom_A(M,\Tuad(\psi)).
\end{equation}
Since
its target is an $\BFp$-vector space, this homomorphism factors through the maximal quotient $J_K$ of $I_K$ that is abelian of exponent~$p$, inducing a homomorphism
\UseTheoremCounterForNextEquation
\begin{equation}\label{InertiaTadMJIntro}
\rho_M\colon J_K\ \longto\ \Hom_A\bigl(M,\:\Tuad(\psi)\bigr).
\end{equation}
This corresponds to a bilinear map $M\times J_K \to \Tuad(\psi)$ that is $A$-linear in the first variable.

\medskip
This pairing can be computed uniformly as follows. 
For later use we pass to the perfection $\Kperf$ of~$K$ at this point. So let $(\KOKperf,\psi)$ denote the $A$-module with the underlying additive group $\KOKperf:=\Kperf/\CO_\Kperf$ and the action of $A$ through~$\psi$. Then we construct a 
\emph{local Kummer pairing} 
\UseTheoremCounterForNextEquation
\begin{equation}\label{KummerProp2Intro}
\xymatrix{
\kum{\ }{\ }{\psi}\colon\ (\KOKperf,\,\psi) \times J_K \ar[r] \ & \ \Tuad(\psi) \rlap{,}}
\end{equation}
which is $A$-linear in the first variable, and show that it induces the homomorphism (\ref{InertiaTadMIntro}) via the embedding $M\into\KOKperf$, $m\mapsto [m]$.

The study of this pairing is at the heart of our paper.
The main question is to determine the image of the adjoint map
\UseTheoremCounterForNextEquation
\begin{equation}\label{KumPsiIntro}
J_K\ \longto\ \Hom_A\bigl((\KOKperf,\psi), \:\Tuad(\psi)\bigr), \quad
[\gal]\ \longmapsto\ \kum{\ }{\gal}{\psi}.
\end{equation}
Here we encounter a somewhat hidden but serious complication. The problem is that the image of (\ref{KumPsiIntro}) obeys symmetries that stem from endomorphisms of the reduction~$\bar\phi$ and not just from endomorphisms of~$\psi$. Since $\bar\phi$ has special characteristic, further symmetries can also arise from endomorphisms of $\bar\phi|A'$ for admissible coefficient rings $A'\subset A$, as in~\cite{PinkDrinSpec2}. We deal with these obstacles as follows. 

\medskip
First we study the case of the \emph{basic Drinfeld module} $\BDM$. This Drinfeld module has the coefficient ring $B:=\BFp[s]$ and is defined by the formula $\BDM_s := \tau^d$ where $d:=[k/\BFp]$ and $\tau$ represents the $p$-Frobenius endomorphism. It has the endomorphism ring $R:= k[\tau]$ with center $\BDM(B)$ and is in some sense universal among Drinfeld modules defined over~$k$.

The Drinfeld module $\BDM$ is supersingular with characteristic ideal $\Fqres:=(s)$, so $\Tu_\Fqres(\BDM) = 0$ and $s$ acts by an automorphism on $\Tuad(\BDM)$. The action of the endomorphism ring $R$ on $\Tuad(\BDM)$ thus extends uniquely to a left action of the ring $R^\circ := R[(\BDM_s)^{-1}] = k[\tau^{\pm 1}]$.
On the other hand the natural left action of $R^\circ$ on $\Kperf$ induces an action on~$\KOKperf$.
By functoriality the local Kummer pairing of $\BDM$ is $R^\circ$-linear in the first variable.
One of our main result states that this pairing is adjoint to an isomorphism
\UseTheoremCounterForNextEquation
\begin{equation}\label{KumBDMIntro}
J_K\ \longisoto\ \Hom_{R^\circ}\bigl(\KOKperf, \:\Tuad(\BDM)\bigr), \quad
[\gal]\ \longmapsto\ \kum{\ }{\gal}{\BDM}
\end{equation}
(see Theorem~\ref{PerfectAndFiltPerf}). This can be interpreted as saying that the local Kummer pairing of $\BDM$ is perfect on the side of~$J_K$. We also establish a similar perfectness property on the side of $\KOKperf$ (see Theorem~\ref{AltPerfectAndFiltPerf}).

\medskip
To achieve these results we use
the ramification filtration on the inertia group~$I_K$. Recall that this is a decreasing filtration by closed normal subgroups $I_K^\nu$ in the upper numbering. As the group $J_K$ is abelian, the Hasse-Arf theorem implies that the induced filtration of $J_K$ jumps only at integers. 

This filtration interacts with the Kummer pairing of $\BDM$ as follows. 
We construct an increasing filtration of $\KOKperf$ by finitely generated free left $R^\circ$-sub\-modules $\Wi{\KOKperf}$, whose subquotients $\grWi{\KOKperf} := \Wii{\KOKperf}/\Wi{\KOKperf}$ are free $R^\circ$-modules of rank $1$ if $p\nmid i>0$ and zero otherwise. We prove that for every $i\ge0$, the isomorphism \eqref{KumBDMIntro} identifies the ramification subgroup $J_K^i$ with the subgroup of homomorphisms that vanish on $\Wi{\KOKperf}$ and that it induces an isomorphism
\UseTheoremCounterForNextEquation
\begin{equation}\label{KumBDMIntroGr}
\gr^i J_K\ :=\ J_K^i / J_K^{i+1}\ \longisoto\ \Hom_{R^\circ}\bigl(\grWi{\KOKperf},\:\Tuad(\BDM)\bigr)
\end{equation}
(see Theorems \ref{PerfectAndFiltPerf} and~\ref{GrJKIsom}). In particular, for any $p\nmid i>0$ this induces an (uncanonical) isomorphism $\gr^i J_K \isoto \Tuad(\BDM)$.

\medskip
Next we look at the action of Frobenius. Note that any $\BFp$-vector space together with an action of $\Gal(\kalg/k)$ can be viewed as a module over the ring $B^\circ := B[s^{-1}] = \BFp[s^{\pm 1}]$ with $s$ acting through the Frobenius of~$k$. If the vector space is profinite and the action is continuous, this action extends uniquely to a continuous action of the profinite completion $\Bado \cong \prod_{\Fq\not=\Fqres} B_\Fq$.
In particular this turns $J_K$ and its ramification subgroups \smash{$J_K^i$} into topological $\Bado$-modules. Also, the induced $\Bado$-module structure on $\Tuad(\BDM)$ agrees with the one from its construction as the Tate module of~$\BDM$. 

Since the local Kummer pairing is Galois equivariant in the first variable, it follows that the isomorphisms \eqref{KumBDMIntro} and \eqref{KumBDMIntroGr} are $\Bado$-linear, and this has an interesting consequence. As $\Tuad(\BDM)$ is a free $\Bado$-module of rank~$d$, the isomorphism \eqref{KumBDMIntroGr} implies that each graded piece
$\gr^i J_K$ is a free $\Bado$-module of rank
\UseTheoremCounterForNextEquation
\begin{equation}\label{JKModStructIntro}
\rank\, \gr^iJ_K = 
\biggl\{\begin{array}{cl}
0 & \hbox{if $p\kern1pt|\kern1pti$,} \\[3pt]
d & \hbox{if $p\nmid i$.}
\end{array}
\end{equation}
Consequently, each quotient $J_K/J_K^i$ is a finitely generated free $\Bado$-module.

\medskip
Transferring the results about $\BDM$ to the arbitrary good reduction Drinfeld module $\psi$ involves two major steps. The passage from $\BDM$ to $\bar\phi$ is the easier one: There is a canonical $\Bado$-linear isomorphism
\UseTheoremCounterForNextEquation
\begin{equation}\label{TadToBasicIntro}
\Tuad(\bar\phi)\ \cong \ \Tuad(\BDM),
\end{equation}
and the composite of the local Kummer pairing of $\bar\phi$ with this isomorphism is exactly the local Kummer pairing of $\BDM$ (see Section~\ref{CompBarPhiBarOmega}). The homomorphism (\ref{KumPsiIntro}) for $\bar\phi$ is therefore injective, but for the reasons mentioned above it is only rarely surjective.

\medskip
The passage from $\bar\phi$ to $\psi$ is much more subtle. For this we first show that there exists a unique power series $x = \sum_{i\le0}x_i\tau^i\in \CO_K[[\tau^{-1}]]$ that is congruent to $1$ modulo $\Fm_K$ and satisfies $\psi_a x = x \bar\phi_a$ for all $a \in A$. This can be viewed as a kind of canonical analytic isomorphism from $\bar\phi$ to~$\psi$.
Thus for any $\xi\in\Kperf$ we have $v(x_i\xi^{p^i})\ge0$ for all $i\ll0$, and so the sum
\UseTheoremCounterForNextEquation
\begin{equation}\label{XFormIntro}
\chi([\xi])\ :=\ \sum\nolimits_{i\le0} [x_i \xi^{p^i}]\ \in\ \KOKperf
\end{equation}
has only finitely many nonzero terms and is therefore well-defined. We show that this induces a natural isomorphism of $A$-modules
\UseTheoremCounterForNextEquation
\begin{equation}\label{XIsomIntro}
\chi\colon (\KOKperf,\,\bar\phi)\ \longisoto\ (\KOKperf,\,\psi)
\end{equation}
(see Proposition~\ref{XIsom}). As this isomorphism involves negative powers of Frobenius in an essential way, it does not preserve the submodule $K/\CO_K\subset \KOKperf$ in general. This is the reason for working with $\Kperf$ in place of~$K$.

\medskip
On the other hand the fact that $\psi$ reduces to $\bar\phi$ induces a natural isomorphism 
\UseTheoremCounterForNextEquation
\begin{equation}\label{TadPsiBarPhiIntro}
\Tuad(\psi)\ \cong\ \Tuad(\bar\phi).
\end{equation}
We prove that under this isomorphism and the isomorphism~$\chi$, the local Kummer pairing for $\bar\phi$ corresponds exactly to the local Kummer pairing for~$\psi$ (see Proposition~\ref{KummerRed}). Thus the homomorphism (\ref{KumPsiIntro}) for $\psi$ is injective, but again only rarely surjective.
It may be somewhat surprising that in this way, the local Kummer pairing for an arbitrary good reduction Drinfeld module $\psi$ is induced by the local Kummer pairing for~$\BDM$, which may therefore be viewed as universal in some sense.

\medskip
In Section~\ref{TLKPITGC} we transfer some further results about the local Kummer pairing of $\BDM$ to that of~$\psi$. In particular, we prove that for every $\xi\in K$ of normalized valuation $-i < 0$ with $p\nmid i$, the local Kummer pairing of $\psi$ induces an isomorphism
\UseTheoremCounterForNextEquation
\begin{equation}\label{KumPsiIntroEv}
\gr^i J_K\ \longisoto \Tuad(\psi), \quad [\gal]\ \longmapsto\ \kum{\xi}{\gal}{\psi}.
\end{equation}
(see Theorem \ref{GrJKIsomPsi}).
This can be seen as the analog of \eqref{KumBDMIntroGr} for $\psi$. 

\medskip
Let us now return to a finitely generated $A$-submodule $M\subset(K,\psi)$ with $M\cap\CO_K=\{0\}$. We identify $M$ with its image in $\KOKperf$ and consider the left $R$-submodule $\Mbar\subset\KOKperf$ that is generated by~$\chi^{-1}(M)$. Combining various properties of the local Kummer pairings we obtain a natural commutative diagram of $\Bado$-linear maps
\UseTheoremCounterForNextEquation
\begin{equation}\label{MainImageDiagIntro}
\vcenter{\xymatrix@R+10pt{J_K \ar[r]^-{\rho_M} \ar[d]_{\rho_\Mbar} & 
\Hom_A\bigl(M,\:\Tuad(\psi)\bigr) \ar[d]^\wr \\
\Hom_R\bigl(\Mbar,\:\Tuad(\BDM)\bigr) \ar@{^{ (}->}[r] & 
\Hom_A\bigl(M,\:\Tuad(\BDM)\bigr) \\}}
\end{equation}
(see Proposition~\ref{MainImageDiag}). This has several consequences. 

On the one hand the isomorphy of (\ref{KumBDMIntro}) implies that the left vertical arrow in (\ref{MainImageDiagIntro}) has open image. Since $\Mbar$ is a free $R$-module of finite rank, this shows that $\rho_M(J_K)$ is a free $\Bado$-module of an explicit finite rank (see Theorems \ref{Mammage} and~\ref{MammagePsi}). In particular, we obtain a necessary and sufficient criterion for $\rho_M(J_K)$ to have finite index in terms of $M$ and $\Mbar$ alone (see Corollary~\ref{MammagePsiOpen}).
But in the general case the image $\rho_M(J_K)$ is a $\Bado$-submodule that may lie askew with respect to the $A_\ad$-module structure of $\Hom_A(M,\Tuad(\psi))$. 
Also, there does not seem to be an easy way 
of computing it up to finite index directly in terms of $\psi$ and~$M$, without the passage to~$\Mbar$.

On the other hand the fact that $\Mbar$ is finitely generated implies that $\Mbar\subset\Wi{\KOKperf}$ for some $i\ge0$. For this $i$ we then have $\rho_\Mbar(J_K^i) = \rho_M(J_K^i)=0$ (see Theorems \ref{MammageRam} and~\ref{MammageRamPsi}), so the smallest such $i$ is related to a kind of conductor of~$\rho_M$ (see \eqref{ConductorDef} and below).

Also, the image of the ramification filtration under $\rho_M$ has the following property:
Every $\Bado$-module subquotient $\rho_M(J_K^i)/\rho_M(J_K^{i+1})$ is either finite 
or free of rank $d$ (see Theorems \ref{Rammage} and \ref{RammagePsi}).
Moreover, the latter case occurs if and only if $i$ is a break for the filtration of $\Mbar$
induced by the $\Wo$-filtration of $\KOKperf$.

\medskip
While determining the image of $\rho_M$ in the above fashion requires the knowledge of~$\Mbar$, in Section \ref{OnlyM} we also deduce some of its properties without computing $\Mbar$. In particular we give a sufficient condition for $\rho_M(J_K)$ to be open (see Theorem \ref{iOpenness}) and a sufficient condition for an integer $i$ to appear in the set of breaks (see Theorem \ref{iMJump}).

\medskip
In the last section 
we discuss some consequences for the image of the homomorphism
\UseTheoremCounterForNextEquation
\begin{equation}\label{InertiaTadMJFpIntro}
\rho_{M,\Fp}\colon J_K\ \longto\ \Hom_A\bigl(M,\:\Tup(\psi)\bigr)
\end{equation}
induced by $\rho_M$ for a single prime $\Fp$ of~$A$, which describes the action of inertia on the $\Fp$-adic Tate module $\Tu_\Fp(\phi)$. Of course the image of $\rho_{M,\Fp}$ is open if the image of $\rho_M$ is open. Perhaps surprisingly the converse is not generally true. This can happen for instance when the image of the combined map $(\rho_{M,\Fp},\rho_{M,\Fp'})$ for $\Fp\not=\Fp'$ lies obliquely in the group $\Hom_A(M,\Tup(\psi))\oplus \Hom_A(M,\Tu_{\Fp'}(\psi))$ and surjects to both factors (see Example~\ref{OpennessDependsOnFpEx}).

In general one might expect the size of $\rho_{M,\Fp}(J_K)$ to obey a kind of independence of $\Fp$ rule, but this seems to be only approximately true.
We establish some lower and upper bounds (see Theorem~\ref{NPQBound}), but beyond these the situation seems messy. 

For $\Fp\not=\Fpres$ we can, however, prove a few other results. On the one hand we show 
that the smallest integer $i$ with $\rho_{M,\Fp}(J_K^i)=0$ is independent of~$\Fp$ (see Theorem~\ref{PConductor}). 
Thus the conductor of $M$ or of $\phi$ is determined by the $\Fp$-adic Galois representation alone.
Finally, we show that the action of inertia on $T_\Fp(\phi)$ determines the rank of~$\psi$, in analogy to the criterion of N\'eron--Ogg--Shafarevich for the reduction of an abelian variety (see Theorem~\ref{NOS}).

\medskip
There are other analogies between the results described above and known facts from number theory. For instance, the pairing \eqref{KummerProp2Intro} is the Drinfeld module analog of the local Kummer pairing for a semiabelian variety, with $J_K$ playing the role of the tame inertia group $I_K^\tame$.
Moreover, our local Kummer theory for Drinfeld modules replicates the tame local Kummer theory for the group scheme~$\BG_m$ if the latter is reformulated as follows.

On the one hand, our modified Tate module $\Tu_\ad(\psi)$ corresponds to the modified Tate module $\Tuad(\BG_m) \cong \prod_{\ell\ne p} T_\ell(\BG_m)$ constructed from all torsion points in $\BG_m(\Kalg)$ modulo those that vanish modulo $\Fm_\Kalg$.
On the other hand, the module $\KOKperf = \Kperf/\CO_\Kperf$ can be viewed as an analog of the factor group
\UseTheoremCounterForNextEquation
\begin{equation}\label{KumGmIntro1}
\CV_{\Kperf}\ :=\ \BG_m(\Kperf)/\BG_m(\CO_\Kperf).
\end{equation}
Note that the normalized valuation of $K$ induces an isomorphism $\CV_\Kperf\isoto\BZ[p^{-1}]$. Thus the tame local Kummer theory for $\BG_m$ amounts to saying that the Kummer pairing is adjoint to an isomorphism
\UseTheoremCounterForNextEquation
\begin{equation}\label{KumGmIntro2}
I_K^\tame\ \longisoto\ \Hom_{\BZ[p^{-1}]}\bigl(\CV_\Kperf,\:\Tuad(\BG_m)\bigr).
\end{equation}
This can be seen as an analog of \eqref{KumBDMIntro}. Also, evaluating the right hand side at a generator of $\CV_{\Kperf}$ we obtain an isomorphism
\UseTheoremCounterForNextEquation
\begin{equation}\label{KumGmIntro3}
I_K^\tame\ \longisoto\ \Tuad(\BG_m),
\end{equation}
which is comparable to \eqref{KumPsiIntroEv}.

\medskip
Finally, we want to address the relation with the global Kummer theory for Drinfeld modules. While the setup is similar, there is one vital difference between the two cases: In the global case the Tate conjecture, which is a theorem, asserts that, after some reductions, the image of Galois on $T_\ad(\psi)$ is open in the commutator of $\End_K(\psi)$. It is this fact, combined with a theorem of Poonen \cite{Poonen}, which allows one to attain essentially complete results using elementary Kummer theory alone (see \cite{PinkKummer2016}). 

By contrast, in the local case the image of Galois on $\Tu_\ad(\psi)$ can be much too small for this approach. 
Only after the passage from $\psi$ to $\BDM$ we can use the abundant symmetries of $\BDM$ to sufficiently constrain the image of~$\rho_M$. 



\bigskip
\noindent
\textbf{Acknowledgements.}
The first author is supported by an Ambizione fellowship of the Swiss National Science Foundation
(project PZ00P2\_202119).


\section{The setup}
\label{Setup}

Throughout this article we fix a finite field $\BFp$ of prime order~$p$.
Consider a local field  $K$ containing~$\BFp$. Let $\CO_K$ denote its valuation ring with the maximal ideal $\Fm_K$ and the finite residue field $k:=\CO_K/\Fm_K$. 

Fix an algebraic closure $\Kalg$ of~$K$. Let $K^\nr\subset K^\sep$ denote the maximal intermediate field of $\Kalg/K$ that is unramified, respectively separable over~$K$. Let $\Kperf\subset\Knrperf$ denote the perfection of~$K$, respectively of~$\Knr$, that is, the maximal intermediate field of $\Kalg/K$ that is totally inseparable over~$K$, respectively over~$K^\nr$. 

For any intermediate field $K'$ we let $\CO_{K'}$ denote its valuation ring with the maximal ideal~$\Fm_{K'}$. Then the residue field $\kalg := \CO_\Knr/\Fm_\Knr$ is an algebraic closure of~$k$, and we have $k \cong \CO_\Kperf/\Fm_\Kperf$ and $k^\alg \cong \CO_\Ksep/\Fm_\Ksep \cong \CO_\Kalg/\Fm_\Kalg$.
Since $K$ is a local field of positive characteristic, there is a unique homomorphism $k\into\CO_K$ whose composite with the projection $\CO_K\onto k$ is the identity on~$k$. We thereby view $K$ as an extension of~$k$. In the same way we view $\Knr$ and each of the above overfields thereof as an extension of~$k^\alg$.

Recall that the restriction of automorphisms induces natural isomorphisms
\UseTheoremCounterForNextEquation
\begin{eqnarray}\label{GammaKDef}
\GK\ :=\ \Gal(\Kalg/\Kperf) & \longisoto & \Gal(\Ksep/K),\\[3pt]
\UseTheoremCounterForNextEquation
\label{IKDef}
I_K\ :=\ \Gal(\Kalg/\Knrperf) & \longisoto & \Gal(\Ksep/\Knr),
\end{eqnarray}
and that we have a natural short exact sequence 
\UseTheoremCounterForNextEquation
\begin{equation}\label{GammaExactSeq}
\xymatrix{ 1 \ar[r] & I_K \ar[r] & \GK \ar[r] & \Gk := \Gal(\kalg/k) \ar[r] & 1\rlap{.}\\}
\end{equation}
Denote by $J_K$ the maximal quotient of $I_K$ that is abelian of exponent~$p$:
\UseTheoremCounterForNextEquation
\begin{equation}\label{JKDef}
J_K\ :=\ I_K^\ab / (I_K^\ab)^{\times p}.
\end{equation}
Let $\{I_K^\nu \mid \nu \in \BQ_{\geqslant 0}\}$ be the upper numbering ramification filtration of~$I_K$. This is a decreasing separated exhaustive filtration by closed normal subgroups of~$\GK$. By the Hasse--Arf theorem the induced filtration on~$J_K$ jumps only at integers; more precisely, for every rational number $\nu \geqslant 0$ we have $J_K^\nu = J_K^i$ with $i = {\lceil \nu \rceil}$. 

\medskip
Next recall that $\Gk$ is the free profinite group that is topologically generated by the Frobenius automorphism $x\mapsto x^{|k|}$. Consider a profinite $\BFp$-vector space $V$ together with a continuous action of~$\Gk$. Then we can view $V$ as a module over the ring $B^\circ := \BFp[s^{\pm1}]$ with $s$ acting through this Frobenius. Since $V$ is profinite, the action of $B^\circ$ extends uniquely to a continuous action of the profinite completion $\Bado$ of~$B^\circ$. Giving a profinite $\BFp$-vector space $V$ together with a continuous action of~$\Gk$ is therefore equivalent to giving a profinite $\Bado$-module.

In particular, consider the action of $\GK$ by left conjugation on~$I_K$. By construction this induces a continuous action on~$J_K$, and since $J_K$ is abelian, this action factors through an action of~$\Gk$. As $J_K$ is also a profinite $\BFp$-vector space, we can therefore view it as a $\Bado$-module. Moreover, as the subgroups $I_K^\nu$ are normal in~$\GK$, it follows that the subgroups ${J_K^\nu \subset J_K}$ are invariant under~$\Gk$ and are therefore closed $\Bado$-submodules. 
We will determine their structure in Section~\ref{Basic}, the main result being Theorem \ref{JKModStruct}.


\medskip
Finally, for any commutative $\BFp$-algebra $S$ we let $S[\tau]$ denote the set of finite linear combinations $\sum_{i\ge0} x_i\tau^i$ with $x_i\in S$. This becomes an $\BFp$-algebra with the usual addition and the multiplication subject to the commutation rule
\UseTheoremCounterForNextEquation
\begin{equation}\label{TauI}
\tau \, x = x^p \, \tau
\end{equation}
for all $x\in S$. Any $S$-algebra $T$ becomes a left $S[\tau]$-module by setting $f(\xi) := \sum_if_i\xi^{p^i}$ for all $f=\sum f_i\tau^i\in S[\tau]$ and all $\xi\in T$.

\section{The Tate module}
\label{TatModSec}

Let $\phi$ be a Drinfeld $A$-module over~$K$.
Any $K$-algebra $S$ becomes an $A$-module with the given addition and the
action $$A\times S\longto S,\ (a,\xi)\mapsto \phi_a(\xi).$$ To keep track of
$\phi$ we denote this $A$-module by $(S,\phi)$. We denote any $A$-submodule
thereof with the underlying set $H$ by $(H,\phi)$, and likewise for factor
modules. In particular, we have the $A$-module $(\Kalg,\phi)$.

\medskip
Let $F$ be the field of quotients of~$A$. The \emph{adelic Tate module of $\phi$} is defined as
\UseTheoremCounterForNextEquation
\begin{equation}\label{TateTotDef}
T_\ad(\phi)\ :=\ \Hom_A\bigl(F/A,\,(\Kalg,\phi)\bigr).
\end{equation}
It is naturally a finitely generated module over the ring of finite adeles
\UseTheoremCounterForNextEquation
\begin{equation}\label{AdelesDef}
A_\ad\ :=\ \End_A(F/A).
\end{equation}
This ring is naturally the product of the completions $A_\Fp$ at all primes 
$\Fp$ of~$A$, and correspondingly $T_\ad(\phi)$ is the product of the \emph{$\Fp$-adic Tate modules} 
\UseTheoremCounterForNextEquation
\begin{equation}\label{TatePDef}
T_\Fp(\phi)\ :=\ \Hom_A\bigl(F_\Fp/A_\Fp,\,(\Kalg,\phi)\bigr).
\end{equation}

From now on we assume that the characteristic homomorphism $d\phi\colon {A\to K}$ factors through~$\CO_K$.
After replacing $K$ by a finite extension we may and do assume that
$\phi$ has stable reduction. This means that $\phi$ has coefficients in $\CO_K$
and that its reduction modulo $\Fm_K$ is a Drinfeld $A$-module $\bar\phi$ of
positive rank over~$k$. Via the embedding $k\into\CO_K$ we view $\bar\phi$
again as a Drinfeld $A$-module of good reduction over~$\CO_K$.

By Tate uniformization \cite[\S7]{DrinfeldEll} there exist a
unique Drinfeld $A$-module $\psi$ with coefficients in $\CO_K$ and good
reduction, and a unique finitely generated projective $A$-submodule $M$ of $(\Kalg,\psi)$
satisfying $M\cap\CO_\Kalg=\{0\}$, such that $\phi$ is the 
\emph{analytic quotient of $\psi$ by~$M.$} 
More precisely there is a unique everywhere converging power series
$e\in\CO_K[[\tau]]$ with constant coefficient $1$ and satisfying $e\equiv1$ modulo
$\Fm_K$, such that $e\circ\psi_a = \phi_a\circ e$ for all $a\in A$ and
$\ker(e|\Kalg)=M$. After enlarging $K$ again we assume that $M$ is contained
in~$K$. 
(Then we can keep $K$ fixed for the rest of this paper.)
The conditions on $e$ imply that $\psi$ has the same reduction $\bar\phi$ and the same characteristic homomorphism as~$\phi$. 

\medskip
Since $\phi$ has coefficients in $\CO_K$, the $A$-module $(\Kalg,\phi)$ possesses the $A$-submodules $(\CO_\Kalg,\phi)$ and $(\Fm_\Kalg,\phi)$ and the corresponding factor modules.

\begin{Lem}\label{TateModLem}
All these modules are divisible.
\end{Lem}

\begin{Proof}
Take any $a\in A\setminus\{0\}$. By the assumption on $\phi$ we have $\phi_a=\sum_{i=0}^dx_i\tau^i$ for some $d\ge0$ with coefficients $x_i\in\CO_K$ and some $x_i\in\CO_K^\times$. For any $\xi\in\Kalg$ the equation $\phi_a(\eta)=\xi$ then means that $\eta$ is a root of the polynomial $f(X) := \sum_{i=0}^dx_iX^{p^i} - \xi$. As this polynomial is non-constant, it always possesses a root in~$\Kalg$, showing that $(\Kalg,\phi)$ is divisible. 

Next observe that $f$ has constant coefficient $-\xi$ and at least one higher coefficient in~$\CO_K^\times$. In the case $\xi\in\CO_\Kalg$ this implies that the first slope of the Newton polygon of $f$ is non-positive; hence $f$ has a root in $\Kalg$ with non-negative valuation. This is the desired $\eta\in\CO_\Kalg$ with $\phi_a(\eta)=\xi$, proving that $(\CO_\Kalg,\phi)$ is divisible. 
In the case $\xi\in\Fm_\Kalg$ it implies that the first slope of the Newton polygon of $f$ is negative; hence $f$ has a root with positive valuation. This is the desired $\eta\in\Fm_\Kalg$ with $\phi_a(\eta)=\xi$, proving that $(\Fm_\Kalg,\phi)$ is divisible. 

Finally, as all these submodules are divisible, so are the factor modules in question.
\end{Proof}

\medskip

As a variant of the adelic Tate module of $\phi$ we now define
\UseTheoremCounterForNextEquation
\begin{equation}\label{TateTotPrimeDef}
\Tuad(\phi)\ :=\ \Hom_A\bigl(F/A,\,(\Kalg/\Fm_\Kalg,\phi)\bigr).
\end{equation}
We use the same definition with $\psi$ or $\bar\phi$ in place of~$\phi$. Then:

\begin{Prop}\label{TateMod1Prop}
There are natural isomorphisms 
$$\Tuad(\psi)\ \cong\ \Tuad(\bar\phi)\ \cong\ T_\ad(\bar\phi).$$
\end{Prop}

\begin{Proof}
Since $\psi$ has good reduction, for every $a\in A\setminus\{0\}$ we have $\psi_a\in\CO_K[\tau]$ with highest coefficient a unit. 
For every $\xi\in\Kalg\setminus\CO_\Kalg$ we therefore have $\psi_a(\xi)\in\Kalg\setminus\CO_\Kalg$; 
hence the $A$-module $(\Kalg/\CO_\Kalg,\psi)$ is torsion free. 
It follows that 
$$\Tuad(\psi)\ \cong\ \Hom_A\bigl(F/A,\,(\CO_\Kalg/\Fm_\Kalg,\phi)\bigr).$$
But the action of $A$ through $\psi$ on $\CO_\Kalg/\Fm_\Kalg$ coincides with the action through~$\bar\phi$. Therefore this Tate module is naturally isomorphic to $\Tuad(\bar\phi)$. Moreover, since $\bar\phi$ is defined over the residue field~$k$, it possesses no non-zero torsion points in $\Fm_\Kalg$, while $\CO_\Kalg/\Fm_\Kalg  = \kalg$ consists only of torsion points. Thus we have a natural isomorphism $T_\ad(\bar\phi)\ \cong\ \Tuad(\bar\phi)$.
\end{Proof}

\medskip

Like the usual Tate module, the modified Tate module is naturally the product of its $\Fp$-adic versions
\UseTheoremCounterForNextEquation
\begin{equation}\label{AltTatePDef}
\Tup(\phi)\ :=\ \Hom_A\bigl(F_\Fp/A_\Fp,\,(\Kalg/\Fm_\Kalg,\phi)\bigr).
\end{equation}
Let $\Fpres\subset A$ denote the characteristic ideal of~$\bar\phi$, also known as the \emph{residual characteristic} of~$\phi$. Since $\bar\phi$ is defined over a finite field, this is a maximal ideal. 

\begin{Prop}\label{TateMod0Prop}
There are natural epimorphisms
$$\begin{array}{cl}
T_\ad(\phi)\ \longonto\ \Tuad(\phi),\\[3pt]
T_\Fp(\phi)\ \longonto\ \Tup(\phi),
\end{array}$$
which for all $\Fp\not=\Fpres$ give isomorphisms 
$$T_\Fp(\phi)\ \longisoarrow\ \Tup(\phi).$$
\end{Prop}

\begin{Proof}
Consider the short exact sequence of $A$-modules
$$\xymatrix{ 0 \ar[r] & (\Fm_\Kalg,\phi) \ar[r] & (\Kalg,\phi) \ar[r] & (\Kalg/\Fm_\Kalg,\phi) \ar[r] & 0\rlap{.}\\}$$
Since $(\Fm_\Kalg,\phi)$ is divisible by Lemma \ref{TateModLem}, applying the functor $\Hom_A(F/A,\underline{\ \ })$ yields a short exact sequence
$$\xymatrix{ 0 \ar[r] & \Hom_A\bigl(F/A,(\Fm_\Kalg,\phi)\bigr) \ar[r] & T_\ad(\phi) \ar[r] & \Tuad(\phi) \ar[r] & 0\rlap{.}\\}$$
In particular we obtain the stated epimorphisms.

Now take any $a\in A\setminus \Fpres$. Then $\phi_a$ lies in $\CO_K[\tau]$ and its lowest coefficient lies in~$\CO_K^\times$. Thus $\phi_a$ does not possess any non-zero root in $\Fm_\Kalg$. The torsion $A$-submodule of $(\Fm_\Kalg,\phi)$ therefore consists of $\Fpres$-power torsion only, and the last statement follows.
\end{Proof}

\medskip

Next, as a variant of the \emph{extended adelic Tate module of $(\psi,M)$} from \cite[\S2]{PinkKummer2016} we define
\UseTheoremCounterForNextEquation
\begin{equation}\label{TateTotMPrimeDef}
\Tuad(\psi,M)\ :=\ \Hom_A\bigl(F/A,\,(\Kalg/(\Fm_\Kalg+M),\psi)\bigr).
\end{equation}
Like its brothers it is naturally the product of its $\Fp$-adic versions
\UseTheoremCounterForNextEquation
\begin{equation}\label{TatePMPrimeDef}
\Tup(\psi,M)\ :=\ \Hom_A\bigl(F_\Fp/A_\Fp,\,(\Kalg/(\Fm_\Kalg+M),\psi)\bigr).
\end{equation}

\begin{Prop}\label{TateMod2Prop}
There is a natural isomorphism
$$\Tuad(\psi,M)\ \cong\ \Tuad(\phi).$$
\end{Prop}

\begin{Proof}
By Tate uniformization the power series $e$ induces a short exact sequence of $A$-modules
$$\xymatrix{ 0 \ar[r] & M \ar[r] & (\Kalg,\psi) \ar[r]^-e & (\Kalg,\phi) \ar[r] & 0\\}$$
and therefore an isomorphism $(\Kalg,\psi)/M \longisoarrow (\Kalg,\phi)$. Since $e$ possesses an 
inverse in $\CO_K[[\tau]]$, it also induces an isomorphism $(\Fm_\Kalg,\psi)/M \longisoarrow (\Fm_\Kalg,\phi)$. As $M\cap\CO_K=\{0\}$, it induces an isomorphism $(\Kalg/(\Fm_\Kalg+M),\psi) \longisoarrow (\Kalg/\Fm_\Kalg,\phi)$ as well. The desired isomorphism now results on applying the functor $\Hom_A(F/A,\underline{\ \ })$.
\end{Proof}

\begin{Prop}\label{TateMod3Prop}
There is a natural short exact sequence 
$$\xymatrix{ 0 \ar[r] & \Tuad(\psi) \ar[r] & \Tuad(\psi,M) \ar[r] & M\otimes_AA_\ad \ar[r] & 0.\\}$$
\end{Prop}

\begin{Proof}
For any $A$-module $L$ and any $A$-submodule $L'$ we set
$$\begin{array}{rcl}
L_\tor &:=& \bigl\{ x\in L \bigm| \exists a\in A\setminus\{0\}\colon ax=0 \bigr\}, \\[3pt]
\Div_L(L') &:=& \bigl\{ x\in L \bigm| \exists a\in A\setminus\{0\}\colon ax\in L' \bigr\}.
\end{array}$$
Since $(\Kalg,\psi)$ is divisible, we obtain a natural short exact sequence
$$\xymatrix{ 0 \ar[r] & (\Kalg,\psi)_\tor \ar[r] & 
\Div_{(\Kalg,\psi)}(M) \ar[r] & M\otimes_A F \ar[r] & 0\rlap{,}\\}$$
where the map on the right hand side is described by $\xi\mapsto \psi_a(\xi)\otimes\frac{1}{a}$ for any $a\in A\setminus\{0\}$ satisfying $\psi_a(\xi)\in M$. Dividing by the submodule $(\Fm_\Kalg,\psi)$, which is divisible by Lemma \ref{TateModLem}, yields a short exact sequence
$$\xymatrix{ 0 \ar[r] & (\Kalg/\Fm_\Kalg,\psi)_\tor \ar[r] & 
\Div_{(\Kalg,\psi)}(M)/(\Fm_\Kalg,\psi)_\tor \ar[r] & M\otimes_A F \ar[r] & 0\rlap{,}\\}$$
Dividing further by the $A$-module~$M$, which is projective and therefore torsion free and satisfies $M\cap\Fm_\Kalg=\{0\}$, we obtain a short exact sequence 
$$\xymatrix{ 0 \ar[r] & (\Kalg/\Fm_\Kalg,\psi)_\tor \ar[r] & 
(\Kalg/(\Fm_\Kalg+M),\psi)_\tor \ar[r] & M\otimes_A F/A \ar[r] & 0\rlap{.}\\}$$
Here all $A$-modules are divisible; hence the proposition follows on applying the functor $\Hom_A(F/A,\underline{\ \ })$.
\end{Proof}

\begin{Prop}\label{TateRank}
Let $r$ denote the rank of $\phi$ and $h$ the height of~$\bar\phi$. Then for each 
prime $\Fp$ of $A$ the module $\Tup(\phi)$ is finitely generated free over $A_\Fp$ of rank
\begin{equation*}
\rank\, \Tup(\phi) = \biggl\{\hspace{-.5ex}\begin{array}{cl}
r - h &\hbox{if $\Fp = \Fpres$,}\\[3pt]
r &\hbox{otherwise.}
\end{array}
\end{equation*}
In particular $\Tuad(\phi)$ is a finitely generated projective module over the ring~$A_\ad$.
\end{Prop}

\begin{Proof} 
Let $\bar r$ denote the rank of~$\bar\phi$. Then by the general theory of Drinfeld modules the $A_\Fp$-module $T_\Fp(\bar\phi)$ is free of rank $\bar r$ if $\Fp\not=\Fpres$, respectively free of rank $\bar r-h$ if $\Fp=\Fpres$.
By Proposition \ref{TateMod1Prop} the same follows for $\Tup(\psi)$.
Now, since $\psi$ has good reduction, its rank is again~$\bar r$, and by the theory of Tate uniformization the rank of $\phi$ is $r=\bar r+\rank_A(M)$. Thus the statements about $\Tup(\phi)$ follow directly from Propositions \ref{TateMod2Prop} and~\ref{TateMod3Prop}. In particular $\Tuad(\phi)$ is a direct summand of a free $A_\ad$-module of finite rank; hence it is finitely generated projective.
\end{Proof}

\medskip
Finally, each of the above Tate modules can be written as an
inverse limit of finite groups of the form $\Hom_A\bigl(\frac{1}{a}A/A,...\bigr)$ for $a\in A\setminus\{0\}$ and therefore carries a natural profinite topology. Moreover, it inherits a natural continuous $A_\ad$-linear action of $\GK=\Aut(\Kalg/K)$ from its action on~$\Kalg$. All the above homomorphisms and isomorphisms are $\GK$-equivariant. 

Since $\bar\phi$ is defined over~$k$, the action of $\GK$ on $\Tuad(\bar\phi)$ factors through~$\Gk$, and by Proposition \ref{TateMod0Prop} the same follows for $\Tuad(\psi)$. As explained in Section \ref{Setup}, this turns these Tate modules into $\Bado$-modules. 

As $M$ is contained in~$K$,
the action of $\GK$ is trivial on $M\otimes_A A_\ad$. Thus the inertia group $I_K$ acts trivially on the outer two modules of the exact sequence from Proposition \ref{TateMod3Prop}. Its action on the module in the middle therefore differs from the identity by a homomorphism
\begin{equation*}
I_K \longto \Hom_A\!\big(M,\:\Tuad(\psi)\big).
\end{equation*}
Here the right hand side is abelian of exponent~$p$, so the homomorphism factors through a unique homomorphism
\UseTheoremCounterForNextEquation
\begin{equation}\label{InertiaTadMJ}
\rho_M\colon J_K \longto \Hom_A\!\big(M,\:\Tuad(\psi)\big).
\end{equation}
By construction this homomorphism is equivariant under $\Gk$ and therefore $\Bado$-linear. It
is described by the \emph{local Kummer pairing} for~$\psi$, which is constructed in the next section.


\section{The local Kummer pairing}
\label{KumPair}

We keep everything from the preceding section.
For later use we will construct the local Kummer pairing with the perfection $\Kperf$ in place of~$K$.

\begin{Lem}\label{KummerLem}
Fix any $\xi\in\Kperf$ and any $\gal\in I_K$. For any $a\in A\setminus\{0\}$ choose an element $\xi_a\in\Kalg$ with $\psi_a(\xi_a) = \xi$. 
\begin{enumerate}\StatementLabels%
\item\label{KummerLemDepA} The residue class $[\gal(\xi_a)-\xi_a]\in\Kalg/\Fm_\Kalg$ depends only on $\gal$, $\xi$ and~$a$.
\item\label{KummerLemDepB} For any $a\in A\setminus\{0\}$ and $b\in A$ the residue class $[\psi_b(\gal(\xi_a)-\xi_a)]\in\Kalg/\Fm_\Kalg$ depends only on $\gamma$, $\xi$ and the residue class $\tfrac{b}{a}+A\in F/A$.
\item\label{KummerLemWelldef} There is a well-defined $A$-linear map
$$\xymatrix@R=3pt{
\llap{$\kum{\xi}{\gal}{\psi}\colon\ $} F/A \ar[r] \ & \ (\Kalg/\Fm_\Kalg,\psi) \rlap{,}\  \\
 \ \tfrac{b}{a}+A \ar@{|->}[r]\ & \ \bigl[\psi_b(\gal(\xi_a)-\xi_a)\bigr] \rlap{.}\ 
}$$
\end{enumerate}
\end{Lem}

\begin{Proof}
First we observe that any other choice for $\xi_a$ has the form $\xi_a+\eta$ for an element $\eta\in\ker(\psi_a)$. Thus $\eta$ lies in $\CO_\Kalg$ and its residue class modulo $\Fm_\Kalg$ is fixed by $\gal\in I_K$. Therefore $\gal(\eta)-\eta \in \Fm_\Kalg$, proving \ref{KummerLemDepA}.

Next consider any $a,a'\in A\setminus\{0\}$ and $b\in A$. Then the computation $\psi_a(\psi_{a'}(\xi_{aa'})) = \psi_{aa'}(\xi_{aa'}) = \xi$ together with \ref{KummerLemDepA} shows that we may without loss of generality assume that $\psi_{a'}(\xi_{aa'}) = \xi_a$. We then have
$$\psi_b(\gal(\xi_a)-\xi_a)\ =\ 
\psi_b(\gal(\psi_{a'}(\xi_{aa'}))-\psi_{a'}(\xi_{aa'}))\ =\ 
\psi_{ba'}(\gal(\xi_{aa'})-\xi_{aa'})$$
This proves that $\psi_b(\gal(\xi_a)-\xi_a)$ modulo $\Fm_\Kalg$ depends only on the quotient $\tfrac{b}{a}\in F$. As the expression is $A$-linear in~$b$, it induces a well-defined $A$-linear map $F \to (\Kalg/\Fm_\Kalg,\,\psi)$. Finally, in the case $\tfrac{b}{a}\in A$ we can compute the expression with $a=1$, so that $\xi_a=\xi\in\Kperf$ and hence $\psi_b(\gal(\xi_a)-\xi_a)= \psi_b(\xi-\xi) = 0$.
Thus the map factors through $F/A$, proving \ref{KummerLemDepB} and \ref{KummerLemWelldef}.
\end{Proof}

\begin{Prop}\label{KummerProp1}
There is a well-defined pairing
$$\xymatrix@R=3pt{
\ (\Kperf,\psi)\times I_K \ar[r] \ & \ \Tuad(\psi) \rlap{,}\  \\
\kern30pt (\xi,\gal) \kern-30pt \ar@{|->}[r]\ \ \ \ \ \ \ \ & \ \kum{\xi}{\gal}{\psi}}$$
which is $A$-linear on $(\Kperf,\psi)$ and a continuous homomorphism on~$I_K$.
\end{Prop}

\begin{Proof}
By the definition (\ref{TateTotPrimeDef}) of the restricted Tate module, the map $\kum{\xi}{\gal}{\psi}$ in Lemma \ref{KummerLem} \ref{KummerLemWelldef} is an element of $\Tuad(\psi)$, so the pairing is well-defined. Next consider $c\in A$ and set $\xi':=\psi_c(\xi)$. For any $a\in A\setminus\{0\}$ the element $\xi'_a:=\psi_c(\xi_a)$ then satisfies $\psi_a(\xi'_a)=\xi'$. By Lemma \ref{KummerLem} we thus have
$$\kum{\xi'}{\gal}{\psi} (\tfrac{b}{a}+A)
\ =\ \bigl[\psi_b\bigl(\gal(\xi'_a)-\xi'_a)\bigr)\bigr]
\ =\ \bigl[\psi_c\bigl(\psi_b(\gal(\xi_a)-\xi_a)\bigr)\bigr]
\ =\ \bigl[\psi_c\bigl(\kum{\xi}{\gal}{\psi} (\tfrac{b}{a}+A)\bigr)\bigr],$$
proving that the pairing is $A$-linear in~$\xi$.
Also, for any $\gal,\gall\in I_K$ and $a\in A\setminus\{0\}$ we have $\gall(\xi_a)-\xi_a \in\ker(\psi_a)\subset \CO_\Kalg$ and therefore
$$\gal\gall(\xi_a)-\xi_a
\ =\ \gal(\xi_a)-\xi_a + \gal(\gall(\xi_a)-\xi_a)
\ \equiv\ (\gal(\xi_a)-\xi_a) + (\gall(\xi_a)-\xi_a)$$
modulo $\Fm_\Kalg$.
Thus the pairing is a homomorphism in~$\gal$. Finally, for any fixed $\xi\in\Kperf$ and $a\in A\setminus\{0\}$ the map 
$$I_K\longto\Kalg/\Fm_\Kalg,\ \ \gal\longmapsto \kum{\xi}{\gal}{\psi} (\tfrac{1}{a}+A)\ =\  \bigl[\gal(\xi_a)-\xi_a\bigr]$$
is continuous. Thus the pairing is continuous in $\gal$ with respect to the profinite topology on $\Tuad(\psi)$.
\end{Proof}

\begin{Prop}\label{KummerGal}%
\hspace{-2pt}%
The pairing in Proposition \ref{KummerProp1} is $\GK$-equivariant;
that is, for all ${\xi\in\Kperf}$ and $\gal\in I_K$ and $\sigma\in\GK$
we have
$$\kum{\xi}{{}^\sigma\gamma}{\psi}
\ =\ \sigma\bigl(\kum{\xi}{\gamma}{\psi}\bigr).$$
\end{Prop}

\begin{Proof}
Consider any element $\frac{b}{a}\in F$.
Since $\xi\in\Kperf$, the element $\xi_a$ from Lemma \ref{KummerLem} satisfies $\psi_a(\sigma(\xi_a)) = \sigma(\psi_a(\xi_a)) = \sigma(\xi)=\xi$. By Lemma \ref{KummerLem} \ref{KummerLemDepA} the pairing for $\xi$ can therefore be computed equally with $\sigma(\xi_a)$ in place of~$\xi_a$. Thus we have
$$\begin{array}{rl}
\kum{\xi}{{}^\sigma\gamma}{\psi} (\tfrac{b}{a}+A)
&=\ \bigl[\psi_b((\sigma\gamma\sigma^{-1}(\sigma(\xi_a))-\sigma(\xi_a))\bigr] \\[3pt]
&=\ \sigma\bigl(\bigl[\psi_b(\gamma(\xi_a)-\xi_a)\bigr] \\[3pt]
&=\ \sigma\bigl(\kum{\xi}{\gamma}{\psi} (\tfrac{b}{a}+A)\bigr),
\end{array}$$
as desired.
\end{Proof}

\begin{Prop}\label{KummerProp2}
The pairing in Proposition \ref{KummerProp1} factors through a natural pairing
$$\xymatrix@R=3pt{
\ (\Kperf/\CO_\Kperf,\psi) \times J_K \ar[r] \ & \ \Tuad(\psi) \rlap{,}\  \\
\kern70pt ([\xi],[\gal]) \kern-30pt \ar@{|->}[r]\ \ \ \ \ \ \ \ & \ \kum{\xi}{\gal}{\psi},}$$
which is $A$-linear on $(\Kperf/\CO_\Kperf,\psi)$, and continuous and $\Bado$-linear on~$J_K$.
\end{Prop}

\begin{Proof}
For any $a\in A\setminus\{0\}$ the highest coefficient of $\psi_a\in\CO_K[\tau]$ is a unit in~$\CO_K$. Thus for any $\xi\in\CO_\Kperf$ the equation $\psi_a(\xi_a)=\xi$ is integral over $\CO_\Kalg$. Therefore $\xi_a$ lies in~$\CO_\Kalg$, and so for any $\gal\in I_K$ we have $\gal(\xi_a)-\xi_a\in\Fm_\Kalg$. Our pairing therefore vanishes on $\CO_\Kperf\times I_K$. Being $A$-linear in~$\xi$, it therefore factors through a map on $(\Kperf/\CO_\Kperf,\psi) \times I_K$  that is $A$-linear in the first variable.

On the other hand, since the pairing is a continuous homomorphism on $I_K$ and the target $\Tuad(\psi)$ is an abelian group of exponent~$p$, the pairing factors through $J_K$ in the second variable and is continuous on~$J_K$. Being $\GK$-equivariant by Proposition \ref{KummerGal}, the pairing is $\Bado$-linear on $J_K$ by the construction of the $\Bado$-module structures on $J_K$ and $\Tuad(\psi)$.
\end{Proof}

\medskip
We call each of the pairings in Propositions \ref{KummerProp1} and \ref{KummerProp2} the \emph{local Kummer pairing} associated to the Drinfeld module~$\psi$, and for simplicity we denote both of them by $\kum{\ }{\ }{\psi}$.


\begin{Prop}\label{KummerEndo}
The local Kummer pairing is functorial in the Drinfeld module $\psi$. In particular it is equivariant under the action of $\End_K(\psi)$ on~$\Kperf$.
\end{Prop}

\begin{Proof}
Let $\psi'$ be another Drinfeld module of good reduction over $\CO_K$.
Then each morphism $h\colon \psi \to \psi'$ has coefficients in $\CO_K$ and therefore
induces a morphism of modified Tate modules $\Tuad(\psi)\to\Tuad(\psi')$. The compatibility of the pairing with $h$ follows in the same way as $A$-linearity in Proposition~\ref{KummerProp1}.
\end{Proof}


\begin{Prop}\label{KummerM}
The homomorphism $\rho_M$ from (\ref{InertiaTadMJ}) is equal to
$$J_K \longto \Hom_A(M,\:\Tuad(\psi)),\quad
[\gal] \longmapsto \kum{\ }{\gal}{\psi}.$$
\end{Prop}

\begin{Proof}
By the proof of Proposition \ref{TateMod3Prop} the map 
$$\Div_{(\Kalg,\psi)}(M) \longto M\otimes_A F$$ 
is given by $\xi\mapsto \psi_a(\xi)\otimes\frac{1}{a}$ for any $a\in A\setminus\{0\}$ satisfying $\psi_a(\xi)\in M$. Therefore the resulting map 
$$(\Kalg/(\Fm_\Kalg+M),\psi)_\tor \longto M\otimes_A F/A$$ 
is described by $[\xi]\mapsto m\otimes\frac{1}{a}$ for any $a\in A\setminus\{0\}$ and $m\in M$ satisfying $\psi_a(\xi)\equiv m \bmod\Fm_\Kalg$. The induced map
$$\xymatrix@R-0pt{\ \Tuad(\psi,M)\ \ar[r] \ar@{=}[d]
&\ M\otimes_AA_\ad\ =\ M\otimes_A(\End_A(F/A)\ \kern60pt \ar@{=}[d]^\wr \\
\ \Hom_A\bigl(F/A,(\Kalg/(\Fm_\Kalg+M),\psi)\bigr)\ \ar[r]^-{f\,\longmapsto\,g} &\ \Hom_A(F/A,M\otimes F/A)\ \\}\kern-60pt$$
is therefore characterized by the equation $g(\frac{b}{a}+A) = m\otimes[\frac{b}{a}]$ for any $m\in M$ such that $\psi_a(f(\frac{1}{a}+A)) = [m]$ in $\Kalg/(\Fm_\Kalg+M)$. Write $f(\frac{1}{a}+A) = [\xi_a]$ with $\xi_a\in\Kalg$ and set $\xi := \psi_a(\xi_a)$, so that $[\xi]=[m]$. Then by Lemma \ref{KummerLem} and Proposition \ref{KummerProp2} we have
$$\begin{array}{rcl}
(\gal(f)-f)(\tfrac{b}{a}+A)
&=& \bigl[\psi_b(\gal(f(\tfrac{1}{a}+A))-f(\tfrac{1}{a}+A))\bigr] \\[3pt]
&=& \bigl[\psi_b(\gal(\xi_a)-\xi_a)\bigr]  \\[3pt]
&\stackrel{\smash{\ref{KummerLem}}}{=}& \kum{\xi}{\gal}{\psi}(\tfrac{b}{a}+A) \\[3pt]
&\stackrel{\smash{\ref{KummerProp2}}}{=}& \kum{m}{\gal}{\psi}(\tfrac{b}{a}+A).
\end{array}$$
This implies the desired formula.
\end{Proof}

\medskip
At last we note that the construction and all results of this section apply equally to any Drinfeld module of good reduction in place of~$\psi$, in particular to~$\bar\phi$.

\section{The dual ramification filtration}
\label{Filtr}

Throughout the following we abbreviate $R := k[\tau]$. Then the local field $K$ is a left $R$-module with $k$ acting by scalar multiplication and $\tau$ acting via the $p$-th power map, and the valuation ring $\CO_K$ is a left $R$-submodule of~$K$. We are interested in the $R$-module structure of the quotient
\UseTheoremCounterForNextEquation
\begin{equation}\label{KOKDef}
\KOK\ :=\ K/\CO_K.
\end{equation}
The notation $\KOK$ is motivated by the fact that, when viewing an element of $K$ as a Laurent series in a uniformizer, its residue class in $\KOK$ is precisely its principal part.

Likewise, consider the localized ring $R^\circ := R[\tau^{-1}] = k[\tau^{\pm1}]$.
Then $\tau$ induces automorphisms of $\Kperf$ and $\CO_\Kperf$, turning both into left $R^\circ$-modules. We are interested in the $R^\circ$-module structure of the quotient
\UseTheoremCounterForNextEquation
\begin{equation}\label{KOKPerfDef}
\KOKperf\ :=\ \Kperf/\CO_\Kperf.
\end{equation}
We will identify $\KOK$ with its image under the natural embedding $\KOK \into \KOKperf$.

\medskip
Throughout the following we let $v$ be the normalized discrete valuation on $K$ with valuation ring~$\CO_K$, and denote its unique extension to a valuation on~$\Kalg$ again by~$v$.  
For every integer $j>0$ that is not divisible by~$p$ we fix an element $\xi_j\in K$ with ${v(\xi_j)=-j}$. 

\begin{Prop}\label{KOKKOKperfBasis}
\begin{enumerate}\StatementLabels
\item\label{KKbB1}
The left $R$-module $\KOK$ is free with the residue classes $[\xi_j]$ as a basis.
\item\label{KKbB2} 
The left $R^\circ$-module $\KOKperf$ is free with the residue classes $[\xi_j]$ as a basis.
\item\label{KKbB3}
The embedding $\KOK \into \KOKperf$ induces a natural isomorphism $R^\circ\otimes_R \brKOK \isoarrow \KOKperf$.
\end{enumerate}
\end{Prop}

\begin{Proof}
For any integers $p\nmid j>0$ and~$\nu$, the element $\tau^\nu(\xi_j) = \xi_j^{p^\nu}$ has normalized valuation $-jp^\nu$ by construction. Varying $j$ and $\nu$ the residue classes of these elements therefore form a basis of $\KOKperf$ over~$k$. Since $R^\circ=k[\tau^{\pm1}]$, this proves \ref{KKbB2}. The same argument with the additional condition $\nu\ge0$ proves \ref{KKbB1}. 
Finally, \ref{KKbB1} and \ref{KKbB2} together imply \ref{KKbB3}. 
\end{Proof}


\begin{Def}\label{WiKOKperfDef}
For any integer $i\ge0$ we let $\Wi\KOKperf$ denote the free left $R^\circ$-submodule of $\KOKperf$ with the basis $\{[\xi_j] : p\nmid j<i\}$.
This defines a separated exhaustive increasing filtration of $\KOKperf$ 
with the graded pieces
$$\grWi\KOKperf\ :=\ \Wii{\KOKperf}/\Wi{\KOKperf}.$$
\end{Def}

Here the indexing of $\grWi\KOKperf$ is shifted by $1$, which is better suited to our purposes.
The construction
immediately implies:

\begin{Prop}\label{GrKOKperf}
\begin{enumerate}\StatementLabels%
\item\label{KOKBasisII1} 
For every $p\kern1pt|\kern1pti\ge0$ we have $\Wi\KOKperf = \Wii\KOKperf$.
\item\label{KOKBasisII2} 
For every $p\nmid i\ge0$ the left $R^\circ$-module $\grWi\KOKperf$ is free of rank $1$ with the basis $[\xi_i]$.
\end{enumerate}
\end{Prop}

\medskip
Also, the filtration is independent of the choice of the $\xi_j$ by the following fact:

\begin{Prop}\label{WiKOKperfGens}
For every integer $i\ge0$ the $R^\circ$-submodule $\Wi\KOKperf$ is generated by the residue classes of all elements of $K$ (sic!) of normalized valuation~${>-i}$. 
\end{Prop}

\begin{Proof}
The generators $[\xi_j]$ for all $p\nmid j<i$ in the definition of $\Wi\KOKperf$ satisfy the stated conditions $\xi_j\in K$ and $v(\xi_j)>-i$. Conversely consider an arbitrary element $\xi\in K$ with $v(\xi)>-i$. Then its residue class in $\KOK$ is a $k$-linear combination of the residue classes $[\xi_j^{p^\nu}]$ for all $p\nmid j>0$ and $\nu\ge0$ such that $jp^\nu<i$. For these we always have $j<i$ and therefore $[\smash{\xi_j^{p^\nu}}] = [\tau^\nu(\xi_j)] \in \Wi\KOKperf$. Together this implies that $[\xi]\in \Wi\KOKperf$, as desired.
\end{Proof}

\begin{Caut}\label{WiPerfCaut}
\rm The $R^\circ$-submodule $\Wi{\KOKperf}$ does not contain the residue classes of all elements $\xi\in\Kperf$ that satisfy $v(\xi) >-i$. Indeed, for any fixed $i>0$ and $p\nmid j\ge i$ we have $v(\xi_j^{p^\nu}) = -jp^\nu > -i$ for all $\nu\ll0$, although $[\xi_j^{p^\nu}] \not\in \Wi{\KOKperf}$.
\end{Caut}


\medskip
Finally we note that, in principle, any finite computation in $\Kperf$ could be done within~$K$ after replacing $K$ by a purely inseparable finite extension. But in this paper we find it simpler
to keep $K$ fixed and work within the big module $\KOKperf$. Nonetheless, we have:

\begin{Prop}\label{PureInsepExtStable}
Replacing $K$ by a finite extension contained in $\Kperf$ changes neither the Galois group~$\GK$, nor the inertia group $J_K$ with its ramification filtration~$J_K^\bullet$, nor the filtration $\Wo\KOKperf$.
\end{Prop}

\begin{Proof}
Consider a finite extension $K'\subset \Kperf$ of $K$ of degree $p^n$. Then $\tau^n$ induces an isomorphism $K'\isoto K$ which identifies the normalized valuation on $K'$ with the normalized valuation on~$K$. It also induces an isomorphism on all finite separable extensions thereof, so the induced isomorphism on Galois groups does not affect the ramification filtration. In particular it does not change $J_K$ or its subgroups~$J_K^i$.

For the same reason the elements $\tau^{-n}(\xi_j)\in K'$ have normalized valuation $-j$ with respect to~$K'$. Since $\tau^n$ is a unit in~$R^\circ$, and $\Kperf$ is a perfection of $K'$ as well, the $R$-submodule $\Wi{\KOKperf}$ therefore does not change on passage from $K$ to~$K'$.
\end{Proof}

\section{The local Kummer pairing in the basic case}
\label{Basic}

In this section we work out the local Kummer pairing for a Drinfeld module that is in some sense universal among Drinfeld modules over the finite field~$k$. 

\begin{Cons}\label{BDMCons}
\rm Consider the polynomial ring $B := \BFp[s]$.
Write $|k|=p^d$ and let $\BDM\colon B\to k[\tau]$ be the ring homomorphism with $\BDM(s)=\tau^d$. Then $\BDM$ is a supersingular Drinfeld $B$-module of rank~$d$ and characteristic ideal $\Fqres = (s)$ over~$k$.
Thus for every prime $\Fq$ of $B$ the Tate module $\Tuq(\BDM)$ is a free $B_\Fq$-module of rank~$d$ if $\Fq\not=\Fqres$, respectively $0$ if $\Fq=\Fqres$. Therefore $\Tuad(\BDM)$ is a free module of rank $d$ 
over the ring $\prod_{\Fq\not=\Fqres}B_\Fq$. 
This ring coincides with the profinite completion $\Bado$ of the ring $B^\circ = \BFp[s^{\pm1}]$ from Section \ref{Setup}.

By construction the image $\BDM(B) = \BFp[\tau^d]$ is precisely the center of the ring $R = k[\tau]$. In particular we have $\End_k(\BDM)=R$, which acts on $\Tuad(\BDM)$ on the left by functoriality. Recall from Goss \cite[\S4.12]{GossBook} that the total ring of quotients $D := \Quot(R)$ is a central division algebra of dimension $d^2$ over the field $E := \Quot(B)$, that the algebra $D$ is ramified exactly at the places $\Fqres$ and~$\infty$, and that $R$ is a maximal $B$-order in~$D$. For every prime $\Fq\not=\Fqres$ of $B$ the completion $R_\Fq := R\otimes_B B_\Fq$ is therefore isomorphic to the matrix ring $\Mat_{d\times d}(B_\Fq)$. Thus $\Rado := R\otimes_B\Bado$ is isomorphic to the matrix ring $\Mat_{d\times d}(\Bado)$, and the $\Rado$-module $\Tuad(\BDM)$ is isomorphic to the standard module $(\Bado)^{\oplus d}$.
\end{Cons}

\begin{Prop}\label{TateBado}
The $\Bado$-module structure on $\Tuad(\BDM)$ from its definition as Tate module coincides with the $\Bado$-module structure induced by the $\Gk$-action as in Section \ref{Setup}.
\end{Prop}

\begin{Proof}
Since $|k|=p^d$, the Frobenius automorphism $\sigma\in\Gk$ acts on $\kalg$ by $\sigma(x)=x^{p^d}$. As this action coincides with that of $\tau^d = \BDM(s)$, it follows that the Galois action of $\sigma$ on $\Tuad(\BDM)$ coincides with that of~$s$. By continuity the two actions of the profinite completion $\Bado$ of $B^\circ := \BFp[s^{\pm1}]$ are therefore equal.
\end{Proof}

\medskip
Recall from the previous section that for every integer $p\nmid i > 0$ we fixed an element $\xi_i \in K$ of normalized valuation $v(\xi_i)=-i$. 
Next, recall from Section \ref{Setup} that the group $J_K$ carries a decreasing ramification filtration by closed normal subgroups $J^i_K$ that jumps only at integers. The first result of this section describes the subquotients $\gr^iJ_K := J^i_K/J^{i+1}_K$ by means of the local Kummer pairing of~$\BDM$:

\begin{Thm}\label{GrJKIsom}
Consider any integer $i\ge0$.
\begin{enumerate}\StatementLabels%
\item\label{GrJKIsomVanish} If $i$ is divisible by~$p$, then $\gr^iJ_K=0$.
\item\label{GrJKIsomIsom} If $i$ is not divisible by~$p$, the homomorphism $\kum{\xi_i}{\ }{\BDM}$ vanishes on $J_K^{i+1}$ and induces a $\Bado$-linear isomorphism
$$\gr^iJ_K \longisoarrow \Tuad(\BDM), \quad
[\gal] \longmapsto \kum{\xi_i}{\gal}{\BDM}.$$
\end{enumerate}
\end{Thm}

\begin{Proof}
For the duration of the proof we use the following notation for any integer $n\ge1$:
\par\smallskip $\bullet$\ \ $K_n$ is the subextension of $\Knr/K$ of degree~$n$ over~$K$.
\par\smallskip $\bullet$\ \ $k_n$ is the residue field of~$K_n$, which has order~$p^{nd}$.
\par\smallskip $\bullet$\ \ $\GKn$ is the Galois group of $\Ksep/K_n$.
\par\smallskip $\bullet$\ \ $I_n$ is the inertia subgroup of the maximal abelian quotient of $\GKn$.
\par\smallskip $\bullet$\ \ $J_n$ is the maximal quotient of $I_n$ of exponent $p$.
\par\smallskip $\bullet$\ \ $I_n^i \subset I_n$ and $J_n^i \subset J_n$ denote the images of the ramification subgroup $I_K^i$.
\par\medskip\noindent
Then we have $J_K = \invlim_n\!J_n$, and $J^i_K = \invlim_n\!J^i_n$ for every $i\ge0$. Setting $\gr^iJ_n := J^i_n/J^{i+1}_n$, it follows that
\UseTheoremCounterForNextEquation
\begin{equation}\label{InvLimGrJKn}
\gr^iJ_K\ =\ \invlim_n \gr^iJ_n.
\end{equation}
Next recall that the reciprocity isomorphism of local class field theory identifies $I_n$ with the group of units $U_n := \CO_{K_n}^\times$, and the ramification subgroup $I^i_n$ with the subgroup $U_n^i$ of elements that are congruent to $1$ modulo $\Fm_{K_n}^i$.
By looking at residue classes we obtain isomorphisms
$$\gr^iU_n\ :=\ U^i_n/U^{i+1}_n
\ \cong\ \biggl\{\begin{array}{cl}
k_n^\times & \hbox{if $i=0$,} \\[3pt]
\Fm_{K_n}^i / \Fm_{K_n}^{i+1}
& \hbox{if $i>0$.}
\end{array}$$
In particular we have
$$|\gr^iU_n|\ =\ 
\biggl\{\begin{array}{cl}
p^{nd}-1 & \hbox{if $i=0$,} \\[3pt]
p^{nd} & \hbox{if $i>0$.}
\end{array}$$
Moreover the $p$-th power map induces an isomorphism $\gr^iU_n\isoarrow \gr^{pi}U_n$ for every $i>0$. For the maximal quotient that is abelian of exponent $p$ it follows that
\UseTheoremCounterForNextEquation
\begin{equation}\label{GrJKnOrder}
|\gr^iJ_n|\ =\ 
\biggl\{\begin{array}{cl}
1 & \hbox{if $p\kern1pt|\kern1pti$,} \\[3pt]
p^{nd} & \hbox{if $p\nmid i$.}
\end{array}
\end{equation}
Combining this with (\ref{InvLimGrJKn}) proves \ref{GrJKIsomVanish}.

\medskip\newcommand{\xin}{\xi_{i,n}}%
\medskip\newcommand{\xinprim}{\xi_{i,n'}}%
Now suppose that $i$ is not divisible by~$p$ and recall that $\xi_i\in K$ has normalized valuation~$-i$.
For every $n\ge1$ pick an element $\xin\in\Kalg$ with $\xin^{p^{nd}}-\xin=\xi_i$.

\begin{Lem}\label{ArtinSchreierLem}
The extension $K_n(\xin)/K_n$ is totally ramified of degree $p^{nd}$ and Galois. Its Galois group $G_{i,n}$ possesses an isomorphism
$$G_{i,n} \isoto k_n,\ \ \gamma\mapsto\gamma(\xin)-\xin,$$
and its upper numbering ramification filtration has a unique break at~$i$, that is, we have $(G_{i,n})^i = G_{i,n}$ and $(G_{i,n})^\nu = \{1\}$ for any $\nu > i$.
\end{Lem}

\begin{Proof}
In the case $p^{nd}=p$ this is the content of Serre \cite[Ch.\,IV Ex.\,5]{SerreLocalFields}.
The general case is established in exactly the same way.
\end{Proof}

\medskip
Consider the surjective homomorphism $\GKn\onto G_{i,n}$ that corresponds to the embedding ${K_n(\xin)\subset K^\alg}$. As the extension is abelian of exponent~$p$ and totally ramified, this homomorphism induces a surjection $J_n\onto G_{i,n}$. The information on the ramification filtration in Lemma \ref{ArtinSchreierLem} implies that this homomorphism vanishes on $J^{i+1}_n$ and induces a surjection $\gr^iJ_n\onto G_{i,n}$. Using \eqref{GrJKnOrder} and comparing orders we therefore obtain a group isomorphism 
\UseTheoremCounterForNextEquation
\begin{equation}\label{GrJKnKn}
\gr^iJ_n \isoto k_n,\ \ [\gamma]\mapsto \gamma(\xin)-\xin.
\end{equation}

To rewrite this further we observe that for every $n\ge1$ we have $\BDM_{s^n-1} = \tau^{nd}-1$ by the definition of~$\BDM$. Thus $k_n$ is the module of $(s^n-1)$-torsion points of~$\BDM$, and the chosen element $\xin$ satisfies the equation $\BDM_{s^n-1}(\xin)=\xi_i$. Moreover, as in the proof of Lemma \ref{KummerLem} \ref{KummerLemDepA} we are free to assume without loss of generality that $\BDM_{(s^{n'}-1)/(s^n-1)}(\xinprim) = \xin$ for all $n|n'$. Thus combining (\ref{GrJKnKn}) with (\ref{InvLimGrJKn}) we deduce a group isomorphism
\UseTheoremCounterForNextEquation
\begin{equation}\label{GrJKInvLimKn}
\kern12pt \gr^iJ_K \isoto \invlim_nk_n,\ \ [\gamma]\mapsto (\gamma(\xin)-\xin)_{n\ge1},
\end{equation}
where the transition maps in the inverse system are the trace maps $\BDM_{(s^{n'}-1)/(s^n-1)}\colon k_{n'}\onto k_n$ for all $n|n'$. We can identify this inverse limit as follows:

\begin{Lem}\label{GrJKIsomLem}
There is a well-defined isomorphism
\begin{equation*}
\Tuad(\BDM) \longisoarrow \invlim_nk_n, \quad \ell \longmapsto (\zeta_n)_{n\geqslant 1},
\end{equation*}
where $\zeta_n\in k_n$ is characterized by the equation $[\zeta_n] = \ell\bigl(\frac{1}{s^n-1} + B\bigr)$ in $\Kalg/\Fm_\Kalg$.
\end{Lem}

\begin{Proof}
Recall that $E$ denotes the quotient field of~$B$. We first note that, since $\BDM$ is defined over~$k$, the embedding $\kalg\into \Kalg/\Fm_\Kalg$ induces an isomorphism
$$\Hom_B\bigl(E/B,(\kalg,\BDM)\bigr)\ \longisoarrow\ 
\Hom_B\bigl(E/B,(\Kalg/\Fm_\Kalg,\BDM)\bigr)
\ \stackrel{(\ref{TateTotPrimeDef})}{=}\ 
\Tuad(\BDM),$$
as in the proof of Proposition \ref{TateMod1Prop}. Secondly let $B_{(\Fqres)}$ denote the localization and $B_\Fqres\subset E_\Fqres$ the completion at $\Fqres=(s)$, so that we have a natural decomposition of $B$-modules
$$E/B\ \cong\ B_{(\Fqres)}/B \oplus E_\Fqres/B_\Fqres.$$
As $\BDM$ is supersingular, there are no non-zero $\Fqres$-primary torsion points in $(\kalg,\BDM)$; hence we also have a natural isomorphism
$$\Hom_B\bigl(B_{(\Fqres)}/B,(\kalg,\BDM)\bigr)\ \longisoarrow\ 
\Tuad(\BDM).$$
Thirdly, every $b\in B\setminus\Fqres$ is a polynomial in $\BFp[s]$ with non-zero constant term. Thus all zeros of $b$ have finite order in $(\kalg)^\times$ and are therefore zeros of $s^m-1$ for some $m\ge1$. If all zeros have multiplicity $\le p^k$ it follows that $b$ divides $s^n-1$ for $n:= mp^k$. Thus we have
$$B_{(\Fqres)}\ =\ \BFp[s]_{(s)}\ =\ \bigcup_{n\ge1} \tfrac{1}{s^n-1}B.$$
Any $B$-linear map $\ell\colon B_{(\Fqres)}/B \longto (\kalg,\BDM)$ is therefore determined by the system of values $\bigl(\ell(\tfrac{1}{s^n-1}+B)\bigr)_{n\ge1}$. Moreover, by $B$-linearity each value $\ell(\tfrac{1}{s^n-1}+B)$ must lie in the kernel $k_n$ of $\BDM_{s^n-1}$. Conversely, a system of elements $(\zeta_n)_n\in\prod_n k_n$ comes from a $B$-linear map $B_{(\Fqres)}/B \longto (\kalg,\BDM)$ if and only if its entries satisfy $\BDM_{(s^{n'}-1)/(s^n-1)}(\zeta_{n'}) = \zeta_n$ for all $n|n'$.
Together the lemma follows.
\end{Proof}

\medskip
Finally, by the choice of the $\xin$ and the construction of the local Kummer pairing in Lemma \ref{KummerLem} \ref{KummerLemWelldef} we have the equality
\begin{equation*}
\bigl[\gal(\xin)-\xin\bigr]\ =\ \kum{\xi_i}{\gal}{\BDM}\bigl(\tfrac{1}{s^n-1}+B\bigr)
\end{equation*}
in $\Kalg/\Fm_\Kalg$ for every~$n$.
Thus the isomorphy in part \ref{GrJKIsomIsom} of Theorem \ref{GrJKIsom} follows by combining Lemma~\ref{GrJKIsomLem} with the formula \eqref{GrJKInvLimKn}. 
The $\Bado$-linearity follows from Proposition \ref{KummerProp2}.
\end{Proof}


\begin{Thm}\label{JKModStruct}
For every integer $i\geqslant 0$ the subquotient $\gr^iJ_K$ is a finitely generated free $\Bado$-module of rank
\begin{equation*}
\rank\, \gr^iJ_K = 
\biggl\{\begin{array}{cl}
0 & \hbox{if $p\kern1pt|\kern1pti$,} \\[3pt]
d & \hbox{if $p\nmid i$.}
\end{array}
\end{equation*}
\end{Thm}

\begin{Proof}
By Theorem \ref{GrJKIsom} we have $\gr^iJ_K = 0$ if $p|i$ and $\gr^iJ_K \cong \Tuad(\BDM)$ if $p \nmid i$. Since $\Tuad(\BDM)$ is a free $\Bado$-module of rank~$d$, the theorem follows.
\end{Proof}

\medskip
Passing from the subquotients $\gr^iJ_K$ to $J_K$ we obtain:

\begin{Thm}\label{PerfectAndFiltPerf}
The local Kummer pairing of $\BDM$ is adjoint to a $\Bado$-linear isomorphism
\begin{equation*}
J_K\ \longisoarrow\ \Hom_{R^\circ}\bigl(\KOKperf,\:\Tuad(\BDM)\bigr), \quad
[\gal] \longmapsto \kum{\ }{\gal}{\BDM}.
\end{equation*}
This isomorphism identifies the ramification subgroup $J_K^i$ with the subgroup of homomorphisms that vanish on the submodule $\Wi{\KOKperf}$.
\end{Thm}

\begin{Proof}
For any integer $i\ge0$, Proposition \ref{GrKOKperf} shows that $\grWi\KOKperf$ is zero if $p\kern1pt|\kern1pti$, respectively a free left $R^\circ$-module of rank $1$ generated by the residue class of $[\xi_i]$ if $p\nmid i$. In the second case evaluation at $[\xi_i]$ therefore induces an isomorphism
$$\Hom_{R^\circ}\bigl(\grWi\KOKperf, \:\Tuad(\BDM)\bigr)\ \longisoarrow\ \Tuad(\BDM).$$
Theorem \ref{GrJKIsom} thus implies that for every $i\ge0$ we have an isomorphism
$$\gr^iJ_K\ \longisoarrow\ \Hom_{R^\circ}\bigl(\grWi\KOKperf, \:\Tuad(\BDM)\bigr), \quad
[\gal] \longmapsto \kum{\ }{\gal}{\BDM}$$
Using the 5-Lemma and induction this yields an isomorphism
$$J_K/ J_K^i\ \longisoarrow\ \Hom_{R^\circ}\bigl(\Wi\KOKperf, \:\Tuad(\BDM)\bigr), \quad
[\gal] \longmapsto \kum{\ }{\gal}{\BDM}$$
for every $i\ge1$. Taking the limit the theorem follows.
The $\Bado$-linearity follows from Proposition \ref{KummerProp2}.
\end{Proof}

\medskip
To summarize, we see that the local Kummer pairing
$$\KOKperf \times J_K \longto \Tuad(\BDM),\quad
([\xi],[\gal]) \longmapsto \kum{\xi}{\gal}{\BDM}$$
is $R^\circ$-linear on $\KOKperf$ by Proposition \ref{KummerEndo} and $\Bado$-linear in $J_K$ by Proposition \ref{KummerProp2}. Theorem \ref{PerfectAndFiltPerf} can be interpreted as saying that this pairing is perfect on the side of $J_K$. On the other side we have the following form of perfectness.

Recall that $\KOKperf$ is a free $R^\circ$-module of countably infinite rank by Proposition \ref{KOKKOKperfBasis} \ref{KKbB2}. Consequently 
$(\KOKperf,\:\BDM)$ is a free $B^\circ$-module of countably infinite rank. We equip it with the discrete topology and denote the completed tensor product by the symbol~$\complot$. We also denote the subset of continuous homomorphisms by the superscript~$(\ )^\cont$. 

\begin{Thm}\label{AltPerfectAndFiltPerf}
The local Kummer pairing of $\BDM$ is adjoint to an isomorphism of topological $\Bado$-modules
$$\Bado \complot_B (\KOKperf,\,\BDM)\ \longisoarrow\ \Hom_\Bado^\cont\bigl(J_K,\:\Tuad(\BDM)\bigr),\quad 1\otimes\xi \longmapsto \kum{\xi}{\ }{\BDM}.$$
For every $i$ this isomorphism identifies the submodule $\Bado\otimes_B (\Wi{\KOKperf},\,\BDM)$ with the submodule of homomorphisms that vanish on the ramification subgroup $J_K^i$.
\end{Thm}

\begin{Proof}
Recall that $\Rado = \Bado\otimes_BR$ is isomorphic to the matrix ring $\Mat_{d\times d}(\Bado)$, and that the $\Rado$-module $\Tuad(\BDM)$ is isomorphic to the standard module $(\Bado)^{\oplus d}$. Also, we have $\Rado = \Bado\otimes_BR^\circ$, because $R^\circ = B^\circ\otimes_B R$ and $B^\circ$ is a localization of~$B$.
Thus the action of $\Rado$ induces a natural isomorphism of $\Bado$-modules
$$\Bado\otimes_B R^\circ \longisoarrow \End_\Bado(\Tuad(\BDM)).$$
For any integer $i\ge0$ that is not divisible by~$p$, Theorem \ref{GrJKIsom} \ref{GrJKIsomIsom} therefore yields an isomorphism of $\Bado$-modules
$$\Bado\otimes_B R^\circ \longisoarrow \Hom_\Bado\bigl(\gr^iJ_K,\:\Tuad(\BDM)\bigr),\quad1\otimes1\longmapsto\kum{\xi_i}{\ }{\BDM}.$$
By Proposition \ref{GrKOKperf} \ref{KOKBasisII2} it follows that the local Kummer pairing induces an isomorphism of $\Bado$-modules
\UseTheoremCounterForNextEquation
\begin{equation}\label{AltPerfectAndFiltGri}
\Bado\otimes_B\bigl(\grWi\KOKperf,\,\BDM\bigr) \longisoarrow \Hom_\Bado\bigl(\gr^iJ_K,\:\Tuad(\BDM)\bigr),\quad 1\otimes[\xi]\longmapsto\kum{\xi}{\ }{\BDM}.
\end{equation}
By Proposition \ref{GrKOKperf} \ref{KOKBasisII1} and Theorem \ref{GrJKIsom} \ref{GrJKIsomVanish} this is also an isomorphism for $p\kern1pt|\kern1pti$, where both sides are zero.
In addition, the module $\grWi\KOKperf$ is flat over~$B$, because it is free over $R^\circ$. 
Using the 5-Lemma and induction, \eqref{AltPerfectAndFiltGri} thus yields an isomorphism
\UseTheoremCounterForNextEquation
\begin{equation}\label{AltPerfectAndFilti}
\Bado\otimes_B\bigl(\Wi{\KOKperf},\,\BDM\bigr) \longisoarrow \Hom_\Bado\bigl(J_K/J_K^i,\:\Tuad(\BDM)\bigr)
\end{equation}
for every $i\ge0$. Consequently, we get an isomorphism
\UseTheoremCounterForNextEquation
\begin{equation}\label{AltPerfectAndFiltiModn}
B^\circ\kern-2pt/\Fn\otimes_B \bigl(\Wi{\KOKperf},\,\BDM\bigr) \longisoarrow \Hom_\Bado\bigl(J_K/J_K^i,\:\Tuad(\BDM)/\Fn\Tuad(\BDM)\bigr).
\end{equation}
for each nonzero ideal $\Fn\subset B^\circ$. 

Now recall that every open neighbourhood of $1$ in $I_K$ contains a ramification subgroup $I_K^i$ of some index~$i$. As the group $\Tuad(\BDM)/\Fn\Tuad(\BDM)$ is discrete, it follows that every continuous homomorphism from $J_K$ to this group factors through the quotient $J_K/J_K^i$ for some $i$. Taking the direct limit of \eqref{AltPerfectAndFiltiModn} over $i$ therefore yields an isomorphism
$$B^\circ\kern-2pt/\Fn\otimes_B (\KOKperf,\,\BDM) \longisoarrow \Hom_\Bado^\cont\bigl(J_K,\:\Tuad(\BDM)/\Fn\Tuad(\BDM)\bigr).$$
Passing to the inverse limit over $\Fn$ we thus obtain the desired isomorphism. 
%
The statement about the ramification subgroups follows from (\ref{AltPerfectAndFilti}). 
\end{Proof}

\medskip
One can show the isomorphy in Theorem \ref{AltPerfectAndFiltPerf} independently using abstract Galois cohomology and the five-term inflation-restriction sequence of the pair $I_K \lhd \GK$. By contrast we do not know how to avoid local class field theory in the proof of the isomorphy in Theorem \ref{PerfectAndFiltPerf} or in the part of Theorem \ref{AltPerfectAndFiltPerf} that concerns the ramification filtration.
\section{The image of inertia in the basic case}
\label{RamFilt}

Now let $\Mbar$ be a finitely generated left $R$-submodule of $\KOKperf = \Kperf/\CO_\Kperf$.
We shall investigate the homomorphism
\UseTheoremCounterForNextEquation
\begin{equation}\label{RhoMbarFactors}
\rho_\Mbar\colon J_K\ \longto\ \Hom_R\bigl(\Mbar,\:\Tuad(\BDM)\bigr), \quad
[\gal] \ \longmapsto\ \kum{\ }{\gal}{\BDM}.
\end{equation}
This homomorphism is $\Bado$-linear by Proposition \ref{KummerProp2}.

\medskip
Next consider the separated exhaustive increasing filtration of $\Mbar$ by the $R$-submodules
\UseTheoremCounterForNextEquation
\begin{equation}\label{WiMDef}
\Wi{\Mbar}\ :=\ \Wi{\KOKperf} \cap \Mbar.
\end{equation}
that is induced by the filtration from Section \ref{Filtr}.
For every $i\ge0$ consider the subquotient 
\UseTheoremCounterForNextEquation
\begin{equation}\label{grWiMDef}
\grWi\Mbar := \Wii{\Mbar}/\Wi{\Mbar}.
\end{equation}
We are interested in the set of breaks
\UseTheoremCounterForNextEquation
\begin{equation}\label{SMDef}
S_\Mbar\ := \bigl\{ i\in\BZ_{\ge0} \bigm| \grWi\Mbar \ne 0 \bigr\}.
\end{equation}

\begin{Lem}\label{BreaksMbar}
\begin{enumerate}\StatementLabels%
\item\label{BreaksMbarPrime} For any $i\in S_\Mbar$ we have $p\nmid i$, and $\grWi\Mbar$ is a free $R$-module of rank~$1$.
\item\label{BreaksMbarRank} The $R$-module $\Mbar$ is free of rank $|S_\Mbar|$.
\end{enumerate}
\end{Lem}

\begin{Proof}
For every $i\ge0$ the subquotient $\grWi\Mbar$ is isomorphic to a finitely generated left $R$-submodule of $\grWi\KOKperf$. By Proposition \ref{GrKOKperf} the latter is zero for $p\kern1pt|\kern1pti$, respectively a free left $R^\circ$-module of rank $1$ with basis $[\xi_i]$ for $p\nmid i$. For each $i\in S_\Mbar$ we therefore have $p\nmid i$. Since $R^\circ = \bigcup_{n>0} R\tau^{-nd}$ and $R$ is a left principal ideal domain by Goss \cite[Cor.\,1.6.3]{GossBook}, it follows that $\grWi\Mbar$ is a free left $R$-module of rank~$1$. This proves \ref{BreaksMbarPrime}. It also shows that $\Mbar$ is a successive extension of $|S_\Mbar|$ free left $R$-modules of rank~$1$, proving \ref{BreaksMbarRank}.
\end{Proof}


\begin{Thm}\label{Rammage}
For every $i\ge0$ we have:
\begin{enumerate}\StatementLabels%
\item\label{RammageNot} If $i\not\in S_\Mbar$, then $\rho_\Mbar(J_K^i)/\rho_\Mbar(J_K^{i+1})$ is a finite $\Bado$-module.
\item\label{RammageIn} If $i\in S_\Mbar$, then $\rho_\Mbar(J_K^i)/\rho_\Mbar(J_K^{i+1})$ is a free $\Bado$-module of rank~$d$.
\item\label{RammageRank} The image $\rho_\Mbar(J_K^i)$ is a free $\Bado$-module of rank $d\cdot\bigl|\{j\in S_\Mbar\mid j\ge i\}\bigr|.$
\end{enumerate}
\end{Thm}

\begin{Proof}
For all $i$ the images $\rho_\Mbar(J_K^i)$ are $\Bado$-submodules of $\Hom_R(\Mbar,\Tuad(\BDM))$, and their subquotients are $\Bado$-modules.
For every $i$ we can identify $\Hom_R(\Mbar/\Wi{\Mbar},\Tuad(\BDM))$ with the subgroup of homomorphisms in $\Hom_R(\Mbar,\Tuad(\BDM))$ that vanish on $\Wi{\Mbar}$. Thus $\rho_\Mbar$ induces a homomorphism
$$\begin{tikzcd}[column sep=4em] J_K^i \rar["\rho_\Mbar^i"] & \Hom_R\bigl(\Mbar/\Wi{\Mbar},\:\Tuad(\BDM)\bigr)\end{tikzcd}$$
by Theorem~\ref{PerfectAndFiltPerf}. 
On the other hand Lemma \ref{BreaksMbar} implies that we have a short exact sequence of free $R$-modules $0\to \grWi\Mbar \to \Mbar/\Wi{\Mbar} \to \Mbar/\Wii{\Mbar} \to 0$. In particular this sequence splits, so we obtain the following commutative diagram with exact columns:
\UseTheoremCounterForNextEquation
\begin{equation}\label{Rammage1}
\begin{tikzcd}[column sep=5em, row sep=1.5em]
0\dar & 0\dar \\
J_K^{i+1} \rar["\rho_\Mbar^{i+1}"] \dar 
& \Hom_R\bigl(\Mbar/\Wii{\Mbar},\:\Tuad(\BDM)\bigr) \dar \\
J_K^i \rar["\rho_\Mbar^i"] \dar
& \Hom_R\bigl(\Mbar/\Wi{\Mbar},\:\Tuad(\BDM)\bigr) \dar \\
\gr^iJ_K \rar["\gr^i\rho_\Mbar"] \dar & 
\Hom_R\bigl(\grWi\Mbar,\:\Tuad(\BDM)\bigr) \dar \\
0 & 0 
\end{tikzcd}
\end{equation}
In the case $i\in S_\Mbar$, recall from the proof of Lemma \ref{BreaksMbar} that $\grWi\Mbar$ is a free left $R$-submodule of rank $1$ of the free left $R^\circ$-module $\grWi\KOKperf$ with basis $[\xi_i]$. Consequently $\grWi\Mbar$ is generated by $[f_i(\xi_i)]$ for some $f_i\in R^\circ\setminus\{0\}$. Evaluation at $[f_i(\xi_i)]$ thus induces an isomorphism $\Hom_R(\grWi\Mbar,\Tuad(\BDM)) \isoarrow \Tuad(\BDM)$. Its composite with the homomorphism $\gr^i\rho_\Mbar$ from (\ref{Rammage1}) is therefore $f_i$ times the isomorphism $\gr^iJ_K \isoarrow \Tuad(\BDM)$ from Theorem \ref{GrJKIsom} \ref{GrJKIsomIsom}. Since multiplication by $f_i$ induces an injective homomorphism $\Tuad(\BDM) \to \Tuad(\BDM)$ with finite cokernel, it follows that $\gr^i\rho_\Mbar$ is injective with finite cokernel.

Now observe that in the case $i\not\in S_\Mbar$ we have ${\grWi\Mbar=0}$ and hence ${\coker(\gr^i\rho_\Mbar)=0}$. Thus $\coker(\gr^i\rho_\Mbar)$ is finite for all $i\ge0$. Moreover, for all $i\gg0$ we have $\Mbar/\Wi{\Mbar}=0$ and therefore $\smash{\coker(\rho_\Mbar^i)}=0$. Applying the snake lemma to the diagram (\ref{Rammage1}) and using descending induction on $i$ it thus follows that $\smash{\coker(\rho_\Mbar^i)}$ is finite for every $i\ge0$. 

We can now prove the theorem. For any $i\not\in S_\Mbar$, we have $\Mbar/\Wii{\Mbar} = \Mbar/\Wi{\Mbar}$; hence $\rho_\Mbar(J_K^{i+1})$ and $\rho_\Mbar(J_K^i)$ are submodules of finite index of the same ambient $\Bado$-module. This directly implies \ref{RammageNot}.
By contrast, for any $i\in S_\Mbar$ the injectivity of $\gr^i\rho_\Mbar$ implies that the induced homomorphism $\gr^iJ_K \to \rho_\Mbar(J_K^i)/\rho_\Mbar(J_K^{i+1})$ is injective as well. As it is surjective by construction, it is therefore an isomorphism. 
Since $\gr^iJ_K$ is a free $\Bado$-module of rank $d$ by Theorem \ref{JKModStruct}, this proves \ref{RammageIn}.

Finally, Lemma \ref{BreaksMbar} implies that for every $i\ge0$ the quotient $\Mbar/\Wi{\Mbar}$ is a free left $R$-module of rank $|\{j\in S_\Mbar\mid j\ge i\}|$. Thus $\Hom_R(\Mbar/\Wi{\Mbar},\Tuad(\BDM))$ is a free $\Bado$-module of rank $d\cdot|\{j\in S_\Mbar\mid j\ge i\}|$. Since $\smash{\coker(\rho_\Mbar^i)}$ is finite, the image $\rho_\Mbar(J_K^i)$ is a $\Bado$-submodule of finite index thereof. As $\Bado$ is a product of principal ideal domains, it follows that $\rho_\Mbar(J_K^i)$ is itself a free $\Bado$-module of rank $d\cdot|\{j\in S_\Mbar\mid j\ge i\}|$, proving \ref{RammageRank}.
\end{Proof}


\begin{Thm}\label{Mammage}
The image $\rho_\Mbar(J_K)$ is a $\Bado$-module of finite index in $\Hom_R(\Mbar,\Tuad(\BDM))$. In particular it is a free $\Bado$-module of rank $d\cdot\rank_R(\Mbar)$.
\end{Thm}

\begin{Proof}
In the proof of Theorem \ref{Rammage} we have seen that $\smash{\coker(\rho_\Mbar^i)}$ is finite for every $i\ge0$. For $i=0$ this means that $\rho_\Mbar\colon J_K \longto \Hom_R(\Mbar,\Tuad(\BDM))$ has finite cokernel, proving the first statement. In view of Lemma \ref{BreaksMbar} \ref{BreaksMbarRank} the second statement is the special $i=0$ case of Theorem \ref{Rammage} \ref{RammageRank}.
\end{Proof}


\begin{Thm}\label{MammageRam}
For every $i\ge0$, we have $\rho_\Mbar(J_K^i)=0$ if and only if $\Mbar\subset \Wi{\KOKperf}$.
In particular we have $\rho_\Mbar(J_K^i) = 0$ for all $i\gg0$. 
\end{Thm}

\begin{Proof}
By the definition of $S_\Mbar$ we have $\Mbar\subset \Wi{\KOKperf}$ if and only if $i>\sup S_\Mbar$.
Thus the assertion is a special case of Theorem \ref{Rammage} \ref{RammageRank}.
\end{Proof}


\begin{Rem}\label{Bach}
\upshape
The image $\rho_\Mbar(J_K)$ admits the following description.
As $\Mbar$ is finitely generated, there exists a finite extension $K'\subset\Kperf$ of $K$ such that $\Mbar \subset \KOKprime = K'/\CO_{K'}$. Since $\KOKprime$ is a free left $R$-module, the submodule $\Mbar$ is then contained in a direct summand of $\KOKprime$ that is free of finite rank.
Next consider the \emph{saturation}
$$\Mtilde\ :=\ \bigl\{[\xi]\in \KOKprime \bigm| \exists f\in R\setminus\{0\}\colon\ [f(\xi)] \in \Mbar \bigr\}$$
of the $R$-module $\Mbar$ in $\KOKprime$. It is easy to see that $\Mdir$ is a left $R$-submodule of $\KOKprime$ and contained in the same direct summand. From a non-commutative version of the elementary divisor theorem \cite[Ch.\,3, Thm.\,16]{Jacobson1948} it now follows that $\Mdir$ is itself a direct summand of $\KOKprime$ and that the quotient $\Mdir/\Mbar$ is finite. 
The isomorphism $R^\circ\otimes_R \KOKprime \isoarrow \KOKperf$ from Proposition \ref{KOKKOKperfBasis} \ref{KKbB3} thus identifies $R^\circ\otimes_R \Mtilde$ with a direct summand of $\KOKperf$, and so the homomorphism
$$J_K\ \longto\ \Hom_{R^\circ}\bigl(R^\circ\otimes_R\Mtilde,\:\Tuad(\BDM)\bigr)\ \cong\ \Hom_R\bigl(\Mtilde,\:\Tuad(\BDM)\bigr), \quad
[\gal] \longmapsto \kum{\ }{\gal}{\BDM}$$
is surjective by Theorem \ref{PerfectAndFiltPerf}. The image $\rho_\Mbar(J_K)$ therefore is precisely the submodule of finite index
$$\Hom_R\bigl(\Mdir,\:\Tuad(\BDM)\bigr)\ \subset\ \Hom_R\bigl(\Mbar,\:\Tuad(\BDM)\bigr).$$
More generally, letting $\Mdiri$ be the saturation of the  submodule 
$\Mbar/\Wi{\Mbar}\ \subset\ \KOKprime/\Wi{\KOKprime}$,
the image $\rho_\Mbar(J_K^i)$ turns out to be the submodule
$$\Hom_R\bigl(\Mdiri,\:\Tuad(\BDM)\bigr)\ \subset\ \Hom_R\bigl(\Mbar,\:\Tuad(\BDM)\bigr).$$
\end{Rem}

\bigskip
The following example shows that the subquotients $\rho_\Mbar(J_K^i)/\rho_\Mbar(J_K^{i+1})$ in Theorem \ref{Rammage} \ref{RammageNot} can be non-zero torsion modules.

\begin{Ex}\rm
Suppose that $p>2$ and pick an element $b\in B\setminus\Fqbar$. Let $\Mbar$ be the free left $R$-submodule of $\KOKperf$ of rank $2$ that is generated by $[\xi_1]$ and $[\BDM_b(\xi_2)]$. Evaluation at these generators in order induces an isomorphism 
$$\kappa\colon \Hom_R\bigl(\Mbar,\:\Tuad(\BDM)\bigr)\ \longisoarrow\ \Tuad(\BDM) \oplus \Tuad(\BDM).$$
The computation in the proof of Theorem \ref{Rammage} shows that up to a unipotent automorphism of $\Tuad(\BDM) \oplus \Tuad(\BDM)$ that preserves the second summand we have
$$\kappa(\rho_\Mbar(J^i_K))\ =\ \ 
\scriptstyle\left\{\textstyle
\begin{array}{cl}
\phantom{\Tuad(\BDM)}\llap{$0$} \oplus \rlap{$0$}\phantom{b\,\Tuad(\BDM)} 
& \hbox{if $i\ge3$,} \\[3pt]
\phantom{\Tuad(\BDM)}\llap{$0$} \oplus b\,\Tuad(\BDM) 
& \hbox{if $i=2$,} \\[3pt]
\Tuad(\BDM) \oplus b\,\Tuad(\BDM)
& \hbox{if $i\le1$.}
\end{array}\right.$$
Now set $\xi:=\xi_1+\BDM_b(\xi_2)$ and consider the free $R$-submodule $\Mbar{}' \subset\Mbar$ of rank $1$ that is generated by $[\xi]$. Evaluation at $[\xi]$ then induces an isomorphism 
$$\kappa'\colon \Hom_R\bigl(\Mbar{}',\:\Tuad(\BDM)\bigr)\ \longisoarrow\ \Tuad(\BDM).$$
Since $\kappa'\circ\rho_{\Mbar{}'}$ is the composite of $\kappa\circ\rho_\Mbar$ with the addition map $\Tuad(\BDM) \oplus \Tuad(\BDM) \to \Tuad(\BDM)$, it follows that
$$\kappa'(\rho_{\Mbar{}'}(J^i_K))\ =\ \ 
\scriptstyle\left\{\textstyle
\begin{array}{cl}
$0$
& \hbox{if $i\ge3$,} \\[3pt]
b\,\Tuad(\BDM) 
& \hbox{if $i=2$,} \\[3pt]
\Tuad(\BDM)
& \hbox{if $i\le1$.}
\end{array}\right.$$
In particular we have 
$$\rho_{\Mbar{}'}(J_K^1)/\rho_{\Mbar{}'}(J_K^2) 
\ \cong\ \Tuad(\BDM)/b\,\Tuad(\BDM).$$
As $\BDM$ is a Drinfeld $B$-module of rank~$d$, this is a free module of rank~$d$ over $B/(b)$. Varying $b$ the subquotient $\rho_{\Mbar{}'}(J_K^1)/\rho_{\Mbar{}'}(J_K^2)$ can thus become an arbitrarily large finite $\Bado$-module.
\end{Ex}

\section{Preparatory computations with power series}
\label{Prep}


This section is preparatory to the next. 

First consider any commutative $\BFp$-algebra $S$ that is \emph{perfect} in the sense that the map $S\to S$, $x\mapsto x^p$ is bijective. Let $S[[\tau^{-1}]]$ be the set of all formal power series $\sum_{i\le0} x_i\tau^i$ with $x_i\in S$, and let $S\pptauone$ be the set of all formal Laurent series $\sum_{i\in\BZ} x_i\tau^i$ with $x_i\in S$ and $x_i=0$ for all $i\gg0$. Both are again $\BFp$-algebras with the usual addition and the multiplication subject to the commutation rule (\ref{TauI}) for all $i\in\BZ$. Both $S[\tau]$ and $S[[\tau^{-1}]]$ are subrings of $S\pptauone$.


\begin{Prop}\label{RTauInvUnits}
An element $x=\sum_{i\le0} x_i\tau^i \in S[[\tau^{-1}]]$ possesses a two-sided inverse $x^{-1} \in S[[\tau^{-1}]]$ if and only if $x_0$ is a unit in~$S$.
\end{Prop}

\begin{Proof}
The map $\sum_{i\le0} x_i\tau^i \mapsto x_0$ is a homomorphism of unitary rings $S[[\tau^{-1}]] \to S$ and therefore maps units to units. This implies that ``only if'' part. For the ``if'' part assume that $x_0$ is a unit in~$S$ and consider a second element $y=\sum_{j\le0} y_j\tau^j \in S[[\tau^{-1}]]$. Then 
$$xy\ =\ \sum_{i\le0} \sum_{j\le0} x_i\tau^i y_j\tau^j
\ =\ \sum_{i\le0} \sum_{j\le0} x_iy_j^{p^i}\tau^{i+j}
\ =\ \sum_{k\le0} \biggl(\,\sum_{k\le j\le0} \!\! x_{k-j}y_j^{p^{k-j}}\biggr)\tau^k.$$
Thus we have $xy=1$ if and only if
$$x_0y_0=1\quad\hbox{and}\quad
\sum_{\ell\le j\le0} x_{\ell-j}y_j^{p^{\ell-j}} = 0\quad\hbox{for all $\ell<0$}.$$
Since $x_0$ is a unit, these equations are equivalent to
\UseTheoremCounterForNextEquation
\begin{equation}\label{RTauInvUnits1}
y_0=x_0^{-1}\quad\hbox{and}\quad
y_\ell = -x_0^{-1}\cdot\!\! \sum_{\ell<j\le0} \!\! x_{\ell-j}y_j^{p^{\ell-j}}\quad\hbox{for all $\ell<0$}.
\end{equation}
Solving these equations inductively therefore shows that there exists $y\in S[[\tau^{-1}]]$ with $xy=1$. In the same way one shows that there exists $z\in S[[\tau^{-1}]]$ with $zx=1$. The computation $y = 1\cdot y = zx\cdot y = z\cdot xy = z\cdot1 = z$ then implies that $y=z$ is a two-sided inverse of~$x$. This proves the ``if'' part.
\end{Proof}

\medskip
Now we return to our situation of the local field~$K$. 

\begin{Lem}\label{SolveLem}
\begin{enumerate}\StatementLabels%
\item\label{SolveLemMult} For any $y\in1+\Fm_K$ and any integer $n$ not divisible by~$p$ there exists a unique $x\in1+\Fm_K$ with ${x^n=y}$, and this satisfies $v(x-1) \ge v(y-1)$.
\item\label{SolveLemRoot} For any $y\in\CO_K$ and $z\in\Fm_K$ and any integer $n>1$ there exists a unique $x\in\Fm_K$ with ${x-yx^n=z}$, and this satisfies $v(x) \ge v(z)$.
\end{enumerate}
\end{Lem}

\begin{Proof}
In \ref{SolveLemMult} the residue class of $1$ is a solution of the congruence $x^n\equiv y$ modulo $(y-1)\subset \Fm_K$, and since $p \nmid n$ it is a simple solution modulo~$\Fm_K$. By Hensel's lemma there is therefore a unique solution in $1+\Fm_K$ and this solution satisfies $v(x-1) \ge v(y-1)$.

Similarly, in \ref{SolveLemRoot} the residue class of $0$ is a solution of the congruence $x-yx^n\equiv z$ modulo $(z)\subset \Fm_K$ and a simple solution modulo~$\Fm_K$. By Hensel's lemma there is therefore a unique solution in $\Fm_K$ and this solution satisfies $v(x) \ge v(z)$.
\end{Proof}

%
%
%


\begin{Prop}\label{XExists}
Consider an element $f=\sum_{i=0}^d f_i\tau^i$ with $f_i\in\CO_K$ and $f_d\in\CO_K^\times$ and $d>0$. For each~$i$ let $\bar f_i\in k$ denote the residue class of $f_i$ modulo~$\Fm_K$, viewed again as an element of $\CO_K$ via the embedding $k\into\CO_K$, and set $\bar f := \sum_{i=0}^d \bar f_i\tau^i$. Then there exists a unique element $x=\sum_{j\le0} x_j\tau^j$ with all $x_j\in\CO_K$ and $x\equiv 1$ modulo~$\Fm_K$ such that $fx = x \bar f$.
\end{Prop}

\begin{Proof}
The desired equality $fx = x \bar f$ is equivalent to
$$0\ =\ fx - x \bar f
\ =\ \sum_{j\le0} \sum_{0\le i\le d} \bigl(f_i\tau^i x_j\tau^j - x_j\tau^j\bar f_i\tau^i\bigr)
\ =\ \sum_{j\le0} \sum_{0\le i\le d} \bigl(f_ix_j^{p^i} - x_j\bar f_i^{p^j}\bigr)\tau^{i+j}$$
and therefore to
\UseTheoremCounterForNextEquation
\begin{equation}\label{XExists1}
0\ = \!\!\sum_{j\le0,\;0\le i\le d \atop i+j=d+\ell} \!\!\!\bigl(f_ix_j^{p^i} - x_j\bar f_i^{p^j}\bigr)
\ =\ \bigl(f_dx_\ell^{p^d} - x_\ell\bar f_d^{p^\ell}\bigr)
+ \!\!\!\!\sum_{\ell<j\le0,\; i\ge0 \atop i+j=d+\ell} \!\!\!\bigl(f_ix_j^{p^i} - x_j\bar f_i^{p^j}\bigr)
\end{equation}
for all $\ell\le0$. We will show that this can be solved uniquely by downward induction on~$\ell$, subject to the condition $x_\ell\equiv\delta_{\ell,0}\bmod\Fm_K$, where $\delta$ denotes the Kronecker delta function.

In the case $\ell=0$ the equation (\ref{XExists1}) reduces to $f_dx_0^{p^d} - x_0\bar f_d=0$. Here $f_d$ is a unit by assumption; hence so is~$\bar f_d$. As $x_0$ must also be a unit, the equation is equivalent to
\UseTheoremCounterForNextEquation
\begin{equation}\label{XExists2}
x_0^{p^d-1} = \bar f_d f_d^{-1}.
\end{equation}
Since $\bar f_d f_d^{-1} \equiv 1 \bmod\Fm_K$, this has a unique solution $x_0\equiv1\bmod\Fm_K$ by Lemma \ref{SolveLem} \ref{SolveLemMult}, as desired.

Now take $\ell<0$ and suppose that the desired $x_j$ have already been found for all $j>\ell$. Then we observe that since the $p$-Frobenius is bijective on~$k$, all terms $\smash{\bar f_i{}^{p^j}}$ in (\ref{XExists1}) lie in $k\subset K$, although $j$ may be negative. Since $i\ge0$, this shows that all terms except the desired $x_\ell$ already lie in~$K$. Next, in the last sum of (\ref{XExists1}) all summands with $\ell<j<0$ in (\ref{XExists1}) lie in $\Fm_K$, because $x_j\in\Fm_K$. Also, if the term for $j=0$ in (\ref{XExists1}) occurs, it is equal to
\UseTheoremCounterForNextEquation
\begin{equation}\label{XExists3}
f_ix_0^{p^i} - x_0\bar f_i\ =\ f_i(x_0^{p^i}-x_0) + (f_i-\bar f_i)x_0
\end{equation}
with $x_0^{p^i}-x_0\in\Fm_K$ and $f_i-\bar f_i \in \Fm_K$ and thus itself lies in $\Fm_K$. Since $f_d$ is a unit, the equation (\ref{XExists1}) is therefore equivalent to 
\UseTheoremCounterForNextEquation
\begin{equation}\label{XExists4}
x_\ell - f_d \bar f_d{}^{-p^\ell}\!\cdot x_\ell^{p^d}
\ =\ (\hbox{something in $\Fm_K$}).
\end{equation}
This has a unique solution $x_\ell\in\Fm_K$ by Lemma \ref{SolveLem} \ref{SolveLemRoot}, as desired.
\end{Proof}

\medskip
Next, for any element $x = \sum_{i\in\BZ} x_i\tau^i$ of $\CO_\Kperf\pptauone$ we set 
\UseTheoremCounterForNextEquation
\begin{equation}\label{VTauDef}
v(x)\ := \inf\bigl\{ v(x_i) \bigm| i\in\BZ \bigr\}.
\end{equation}

\begin{Prop}\label{XValBound}
The solution $x$ in Proposition \ref{XExists} satisfies $v(x-1) \ge v(f-\bar f) > 0$.
\end{Prop}

\begin{Proof}
Abbreviate $w := v(f-\bar f)$ according to (\ref{VTauDef}), so that $w>0$ by the construction of~$\bar f$. We must then show that $v(x_j-\delta_{j,0}) \ge w$ for all $j\le0$.

First the solution of (\ref{XExists2}) satisfies $v(x_0-1) \ge v(\bar f_d f_d^{-1}-1)\ge w$ by Lemma \ref{SolveLem} \ref{SolveLemMult}, yielding the desired bound for $j=0$. In particular, for all $i\ge0$ we have $\smash{v(x_0^{p^i}-x_0)} \ge w$.
Now take $\ell<0$ and suppose that the desired bound has already been shown for all $j>\ell$. Then all terms with $\ell<j<0$ in (\ref{XExists1}) have valuation $\ge v(x_j) \ge w$. Also, if the term for $j=0$ in (\ref{XExists1}) occurs, by (\ref{XExists3}) its valuation is $\ge\min\{v(x_0^{p^i}-x_0), v(f_i-\bar f_i) \} \ge w$ as well. Thus the right hand side of (\ref{XExists4}) has valuation $\ge w$. Its solution therefore satisfies $v(x_\ell) \ge w$ by Lemma \ref{SolveLem} \ref{SolveLemRoot}. This shows the desired bound for $j=\ell$, and hence by downward induction for all~$j$.
\end{Proof}

\begin{Prop}\label{XInvValBound}
The solution $x$ in Proposition \ref{XExists} possesses a two-sided inverse $x^{-1}$ in $\CO_\Kperf[[\tau^{-1}]]$ which satisfies $v(x^{-1}-1) \ge v(f-\bar f)$.
\end{Prop}

\begin{Proof}
Since $\Kperf$ is perfect, we can apply Proposition \ref{RTauInvUnits} to the ring~$\CO_\Kperf$. As the constant coefficient $x_0$ in Proposition \ref{XExists} is a unit in~$\CO_K$, it follows that $x$ possesses a two-sided inverse $x^{-1} \in \CO_\Kperf[[\tau^{-1}]]$. Abbreviate $w := v(f-\bar f)$ and write out $x^{-1} = \sum_{j\le0} y_j\tau^j$ with all $y_j\in\CO_\Kperf$. We must then show that $v(y_j-\delta_{j,0}) \ge w$ for all $j\le0$.

By Proposition \ref{XValBound} we already have ${v(x_i-\delta_{i,0}) \ge w}$ for all $i\le0$. Thus the first equation $y_0=x_0^{-1}$ in (\ref{RTauInvUnits1}) implies that ${v(y_0-1)\ge w}$ as well. 
Now take $\ell<0$ and suppose that the desired bound has already been shown for all $j>\ell$. Then for all $j$ with $\ell< j\le 0$ we have $\ell-j<0$ and therefore $v(x_{\ell-j})\ge w$. Thus the second equation in (\ref{RTauInvUnits1}) implies that $v(y_\ell)\ge w$ as well. By downward induction the desired bound therefore follows for all~$j$.
\end{Proof}


\begin{Prop}\label{XInvPerf}
We have $x^{-1}=\sum_{j\le0} \tau^jz_j$ with all $z_j\in\CO_K$.
\end{Prop}

\begin{Proof}
By Proposition \ref{XExists} the series $x=\sum_{j\le0} x_j\tau^j$ has coefficients in $\CO_K$. Raising the equations (\ref{RTauInvUnits1}) to the $p^{-\ell}$-th power yields
$$y_0=x_0^{-1}\quad\hbox{and}\quad
y_\ell^{p^{-\ell}} = -x_0^{-p^{-\ell}}\cdot\!\! \sum_{\ell<j\le0} \!\! x_{\ell-j}^{p^{-\ell}} y_j^{p^{-j}}\quad\hbox{for all $\ell<0$}.$$
Since $-\ell\ge0$ and $-j\ge0$, descending induction on $\ell$ thus shows that $y_\ell^{p^{-\ell}} \in \CO_K$ for all~$\ell$. Therefore 
 $x^{-1} = \sum_{j\le0} y_j\tau^j =  \sum_{j\le0} \tau^jy_j^{p^{-j}}$ has the desired property.
\end{Proof}

\medskip 
But in general $x^{-1}$ does not have coefficients in~$K$, as the following example shows.

\begin{Ex}\label{EExistsInvEx}
\rm Choose $\theta\in\Fm_K\setminus\{0\}$ and set $f := \theta+\tau$, so that $\bar f=\tau$. 
Then the equation $fx = x \bar f$ is equivalent to $x^{-1}f = \bar f x^{-1}$, which by a computation similar to that in the proof of Proposition \ref{XExists} with $x^{-1} = \sum_{j\le0} y_j\tau^j$ is equivalent to
$$y_0 - y_0^p=0\quad\hbox{and}\quad
y_\ell - y_\ell^p + y_{\ell+1}\theta^{p^{\ell+1}} = 0\quad\hbox{for all $\ell<0$}.$$
Thus $y_0=1$ and $v(y_\ell) = v(y_\ell-y_\ell^p) = v(y_{\ell+1}\theta^{p^{\ell+1}}) = v(y_{\ell+1}) + p^{\ell+1} v(\theta)$ for all $\ell<0$. By induction we deduce that 
$$v(y_\ell)\ =\ \sum_{\ell<j\le 0} p^j \cdot v(\theta)
\ =\ p^{\ell+1}\cdot\frac{p^{|\ell|}-1}{p-1}\cdot v(\theta)$$
for all $\ell\le 0$. Here $\frac{p^{|\ell|}-1}{p-1}$ is an integer that is not divisible by~$p$, but the exponent of the term $p^{\ell+1}$ decreases linearly with $\ell\le0$. Thus the degree of inseparability of the extension $K(y_\ell)/K$ grows linearly with~$\ell$, and there does not exist any finite extension of~$K$ that contains all~$y_\ell$.
\end{Ex}



\begin{Lem}\label{CommBar}
Any element $h$ of the field $\Kperf\pptauone$ that commutes with some element $\bar f \in k\pptauone\setminus k[[\tau^{-1}]]$ already lies in the subfield $k\pptauone$. 
\end{Lem}

\begin{Proof}
By assumption we have $\bar f = \sum_{i \leqslant d} \bar f_i \tau^i$ with $\bar f_i\in k$ and $\bar f_d \ne 0$ and $d>0$. Write $h = \sum_{j \leqslant e} h_j \tau^j$ with $h_j\in\Kperf$. Equating the coefficients of the series $\bar f h$ and $h \bar f$ we deduce that for every integer $\ell \leqslant e$ we have
$$\sum_{i\le d,\;j\le e\atop i+j = d+\ell} \bigl(\bar f_i\,h_j^{p^i} - h_j\,\bar f_i^{p^j}\bigr) = 0.$$
Solving for the terms with $j=\ell$ we obtain 
\UseTheoremCounterForNextEquation
\begin{equation}\label{CommBarInd}
\bar f_d\,h_\ell^{p^d} - h_\ell\,\bar f^{p^\ell}_d\ =\ - \sum_{\ell<j\le e} \bigl(\bar f_{d+\ell-j}\,h_j^{p^{d+\ell-j}} \! - h_j\,\bar f_{d+\ell-j}^{p^j}  \bigr).
\end{equation}
We will deduce from this that $h_\ell\in k$ by downward induction on~$\ell$. So take any $\ell\le e$ such that $h_j\in k$ for all $j$ satisfying $\ell<j\le e$. Then all terms in the last sum of (\ref{CommBarInd}) lie in~$k$. Since $\bar f_d \in k^\times$ and $d>0$ by assumption, the equation thus implies that $h_\ell$ is algebraic over~$k$. As $k$ is algebraically closed within $\Kperf$, it follows that $h_\ell\in k$, as desired.
\end{Proof}

\begin{Prop}\label{XExistsG}
Let $f$ and $x$ be as in Proposition \ref{XExists} and let $g\in\Kperf\pptauone$ be an element such that $fg = gf$.
Then $g$ lies in $\CO_\Kperf\pptauone$ and its reduction $\bar g$ modulo $\Fm_\Kperf$, viewed again as element of $\CO_\Kperf\pptauone$ via the embedding $k\into\CO_K$, satisfies $g x = x \bar g$.
\end{Prop}

\begin{Proof}
Using the two-sided inverse $x^{-1}$ from Proposition \ref{XInvValBound}, the equation $fx = x \bar f$ from Proposition \ref{XExists} is equivalent to $\bar f = x^{-1}fx$. Since $fg=gf$ by assumption, the element $h := x^{-1}gx$ satisfies $\bar fh = h\bar f$. By Lemma~\ref{CommBar} we therefore have $h \in k\pptauone$. The fact that $x$ and $x^{-1}$ are elements of $\CO_\Kperf\pptauone$ implies that $g = x h x^{-1}$ is also an element of $\CO_\Kperf\pptauone$. As $x$ and $x^{-1}$ reduce to $1$ in $k\pptauone$ and $h$ reduces to itself we deduce that $\bar g = h$.
Hence $g x = x h= x \bar g$, as desired.
\end{Proof}

\section{Reduction of the Kummer pairing}
\label{RedKumPair}

Now we return to our Drinfeld modules $\psi$ and $\bar\phi$ over~$\CO_K$ from Section~\ref{KumPair}, and recall that $\bar\phi$ is the reduction of $\psi$ modulo~$\Fm_K$. Set
\UseTheoremCounterForNextEquation
\begin{equation}\label{WDef}
w\ := \inf\bigl\{ v(\psi_a-\bar\phi_a) \bigm| a\in A \bigr\}.
\end{equation}
Since the valuation is discrete on~$K$ and $v(\psi_a-\bar\phi_a)>0$ for all~$a\in A$, we have $w>0$. 
This number measures the (valuative) distance from $\psi$ to $\bar\phi$.

\begin{Prop}\label{XExistsPsi}
\begin{enumerate}\StatementLabels%
\item\label{XExistsPsiExists} There exists a unique element $x=\sum_{j\le0} x_j\tau^j$ with all $x_j\in\CO_K$ and $x\equiv 1$ modulo~$\Fm_K$ such that $\psi_a x = x \bar\phi_a$ for all $a\in A$.
\item\label{XExistsPsiUnit} This element possesses a two-sided inverse $x^{-1} = \sum_{j\le0} y_j\tau^j$ in $\CO_\Kperf[[\tau^{-1}]]$.
\item\label{XExistsPsiBounds} These elements satisfy the bounds $v(x-1) \ge w$ and $v(x^{-1}-1) \ge w$.
\end{enumerate}
\end{Prop}

\begin{Proof}
Choose any non-constant element $a\in A$. Then, since $\psi$ is a Drinfeld module of positive rank, we have $f := \psi_a = \sum_{i=0}^d f_i\tau^i$ with $f_i\in\CO_K$ and $f_d\neq0$ and $d>0$. Moreover, since $\psi$ has good reduction, we have $f_d\in\CO_K^\times$. Thus $f$ satisfies the assumptions of Proposition \ref{XExists}, and there exists a unique element $x$ with the properties in \ref{XExistsPsiExists} but only for the given choice of $a\in A$. For any other $b\in A$ set $g := \psi_b$. Then since $A$ is commutative we have $fg=gf$, and so we also have $\psi_b x = x \bar\phi_b$ by Proposition \ref{XExistsG}. This proves \ref{XExistsPsiExists}.

Parts \ref{XExistsPsiUnit} and \ref{XExistsPsiBounds} follow directly from Propositions \ref{XValBound} and Proposition \ref{XInvValBound}.
\end{Proof}

\begin{Rem}\label{FormalModuleRem}
\rm The above construction has the following deeper meaning.
Let $F_\infty$ denote the completion at infinity of the coefficient field~$F$
and equip the skew field $K^\perf\pptauone$ with the $\tau^{-1}$-adic topology.
Following \cite[\S7]{MornevT} we consider the category of \emph{formal Drinfeld modules},
whose objects are the continuous homomorphisms $\phi\colon F_\infty\to \Kperf\pptauone$
and whose morphisms $\phi\to\phi'$ are the elements $x \in \Kperf\pptauone$
satisfying $\varphi'_a x = x \varphi_a$ for all $a \in F_\infty$.
Every Drinfeld module $A \to \Kperf[\tau]$ 
extends uniquely to a formal Drinfeld module 
by $\infty$-adic completion (see the proof of \cite[Thm.\,9.2.2]{MornevT}).
In our situation the power series $x$ is a canonical isomorphism of the associated formal Drinfeld modules $\bar\phi\isoto\psi$ 
that is characterized as the unique lift to $\CO_\Kperf$ of the identity morphism $\bar\phi\to\bar\phi$ over~$k$.
\end{Rem}



\begin{Prop}\label{XIsom}
There are natural $A$-module isomorphisms
$$\begin{array}{rll}
\chi\colon (\Kalg/\Fm_\Kalg,\bar\phi) &\longisoarrow& (\Kalg/\Fm_\Kalg,\psi), \\[3pt]
\chi\colon (\Kperf/\Fm_\Kperf,\bar\phi) & \longisoarrow& (\Kperf/\Fm_\Kperf,\psi), \\[3pt]
\chi\colon (\KOKperf,\bar\phi) &\longisoarrow& (\KOKperf,\psi),
\end{array}$$
which on any residue class $[\xi]$ are given by the formulas
$$\chi([\xi]) = \sum_{j\le0} \bigl[x_j \xi^{p^j}\bigr]
\qquad\hbox{and}\qquad
\chi^{-1}([\xi]) = \sum_{j\le0} \bigl[y_j \xi^{p^j}\bigr],$$
where $x=\sum_{j\le0} x_j\tau^j$ and $x^{-1} = \sum_{j\le0} y_j\tau^j$ are as in Proposition \ref{XExistsPsi}.
\end{Prop}

\begin{Proof}
Let $S\subset \CO_\Kperf\pptauone$ be the subset consisting of all series $\sum_j s_j \tau^j$ that satisfy $\liminf_{j\to-\infty} v(s_j) > 0$. One easily shows that $S$ is a subring.
Next observe that for any $\xi \in \Kalg$ we have $\liminf_{j\to-\infty} v(\xi^{p^j})\ge0$. For any $s = \smash{\sum_j}\,s_j \tau^j \in S$ it follows that $v(s_j \xi^{p^j})>0$ for all $j\ll0$.  Thus the residue class $[s_j \xi^{p^j}]$ in $\Kalg/\Fm_\Kalg$ vanishes for all $j\ll0$ and we get a well-defined sum
$$s([\xi])\ :=\ \sum\nolimits_j [s_j \xi^{p^j}] \ \in\ \Kalg/\Fm_\Kalg.$$
By manipulation of double series one shows that this formula defines a left $S$-module structure on $\Kalg/\Fm_{\Kalg}$. The same follows for the quotients $\Kperf/\Fm_\Kperf$ and $\KOKperf = \Kperf/\CO_{\Kperf}$.

Now Proposition \ref{XExistsPsi} \ref{XExistsPsiBounds} shows that $x, x^{-1} \in S$. By construction we have $\chi([\xi]) = x([\xi])$ and $\chi^{-1}([\xi]) = x^{-1}([\xi])$, so the maps $\chi$ and $\chi^{-1}$ are well-defined, and the identities $x x^{-1} = 1 = x^{-1} x$ imply that these maps are mutually inverse bijections.
Finally note that $\CO_{\Kperf}[\tau] \subset S$. Thus for each $a \in A$ the polynomials $\psi_a$ and $\bar\phi_a$ belong to $S$, so the identity $\psi_a x = x \bar\phi_a$ implies that $\chi$ is $A$-linear. 
\end{Proof}

\begin{Prop}\label{XEndos}
Any $h\in\End_K(\psi)$ has coefficients in $\CO_K$. Let $\bar h\in k[\tau]$ be its reduction modulo~$\Fm_K$. Then the isomorphisms in Proposition \ref{XIsom} satisfy $h\circ\chi = \chi\circ\bar h$.
\end{Prop}

\begin{Proof}
The first statement follows from the fact that $\psi$ has coefficients in~$\CO_K$. The same proof as for Proposition \ref{XExistsPsi} with $h$ in place of $\phi_b$ shows that $hx=x\bar h$. The last statement follows as in the proof of the $A$-linearity in Proposition \ref{XIsom}.
\end{Proof}


\begin{Prop}\label{KummerRed}
The following diagram commutes:
$$\xymatrix@C+20pt{
\ (\KOKperf,\psi)\times J_K \ar[r]^-{\kum{\ }{\ }{\psi}} \ \ar@{=}@<+32pt>[d] & \ \Tuad(\psi)\ \ar[d]^\wr_\res \\
\ (\KOKperf,\bar\phi)\times J_K \ar[r]^-{\kum{\ }{\ }{\bar\phi}} \ar@<+10pt>[u]^\wr_\chi\ & \ \Tuad(\bar\phi) \rlap{,}\  \\
}$$
where $\res$ denotes the natural isomorphism from Proposition \ref{TateMod1Prop}.
\end{Prop}

\begin{Proof}
Consider any $\xi,\xi'\in\Kperf$ such that $\chi([\xi])=[\xi']$ under the middle isomorphism in Proposition \ref{XIsom}. For any $a\in A\setminus\{0\}$ pick an element $\xi_a\in\Kalg$ with $\bar\phi_a(\xi_a) = \xi$. Choose $\xi_a''\in\Kalg$ such that $\chi([\xi_a])=[\xi''_a]$ under the first isomorphism in Proposition \ref{XIsom}. 

\begin{Lem}\label{KummerRedLem}
There exists $\xi_a'\in\Kalg$ such that $\psi_a(\xi_a') = \xi'$ and $\xi_a'\equiv\xi_a'' \bmod \Fm_\Kalg$.
\end{Lem}

\begin{Proof}
Choose any $\xi_a'\in\Kalg$ with $\psi_a(\xi_a') = \xi'$. Then the $A$-linearity in Proposition \ref{XIsom} implies that 
$$[\psi_a(\xi_a'')]\ =\ \psi_a([\xi_a''])\ =\ \psi_a(\chi([\xi_a]))\ 
\mathrel{\smash{\stackrel{!}{=}}}\ \chi([\bar\phi_a(\xi_a)])\ =\ \chi([\xi])\ =\ [\xi']\ =\ [\psi_a(\xi_a')]$$
in $\Kalg/\Fm_\Kalg$. Thus $\psi_a(\xi_a'') \equiv \psi_a(\xi_a') \bmod \Fm_\Kalg$ and therefore $\psi_a(\xi_a''-\xi_a') \in \Fm_\Kalg$. Since the highest coefficient of $\psi_a \in \CO_K[\tau]$ is a unit, this implies that $\xi_a''-\xi_a'$ lies in $\CO_\Kalg$ and is congruent modulo $\Fm_\Kalg$ to a unique $\eta\in\ker(\psi_a)$. Thus $\xi_a'' \equiv \xi_a'+\eta \bmod \Fm_\Kalg$, and so the lemma holds with $\xi_a'+\eta$ in place of~$\xi_a'$.
\end{Proof}

\medskip
For any $b\in A$ and $[\gal]\in J_K$ we can now compute
\begin{eqnarray*}
\bigl[\kum{\xi}{\gal'}{\psi}(\tfrac{b}{a}+A)\bigr]
&\stackrel{\ref{KummerLem}}{=}& \bigl[\psi_b(\gal(\xi_a')-\xi_a')\bigr] \\[3pt]
&\stackrel{\ref{KummerRedLem}}{=}& \bigl[\psi_b(\gal(\xi_a'')-\xi_a'')\bigr] \\[3pt]
&=& \psi_b\bigl(\gal([\xi_a''])-[\xi_a'']\bigr) \\[3pt]
&=& \psi_b\bigl(\gal(\chi([\xi_a]))-\chi([\xi_a])\bigr) \\[3pt]
&\stackrel{\smash{(*)}}{=}& \chi\bigl(\bigl[\bar\phi_b(\gal(\xi_a)-\xi_a)]\bigr]\bigr) \\[3pt]
&\stackrel{\ref{KummerLem}}{=}& 
\chi\bigl(\bigl[\kum{\xi}{\gal}{\bar\phi}(\tfrac{b}{a}+A)\bigr]\bigr)
\end{eqnarray*}
within $\Kalg/\Fm_\Kalg$, where at $(*)$ we have used the $A$-linearity and the Galois equivariance of the first map in Proposition \ref{XIsom}. Next we observe that, since the element $x$ from Proposition \ref{XExistsPsi} is congruent to $1$ modulo~$\Fm_K$, the map $\chi$ from Proposition \ref{XIsom} induces the identity on $\CO_\Kalg/\Fm_\Kalg$. Thus the above computation implies that 
$$\llap{$\kum{\xi}{\gal'}{\psi}(\tfrac{b}{a}+A)$}
\ \equiv\ \rlap{$\kum{\xi}{\gal}{\bar\phi}(\tfrac{b}{a}+A)
\ \bmod\ \Fm_\Kalg$} \kern35pt$$
and therefore
$$\llap{$\res(\kum{\xi}{\gal'}{\psi})(\tfrac{b}{a}+A)$}
\ =\ \rlap{$\kum{\xi}{\gal}{\bar\phi}(\tfrac{b}{a}+A).$} \kern35pt$$
Varying $a$ and $b$ this shows that $\res(\kum{\xi}{\gal'}{\psi}) = \kum{\xi}{\gal}{\bar\phi}$, and so the diagram commutes.
\end{Proof}

\section{Passage to the basic Drinfeld module}
\label{CompBarPhiBarOmega}

In this section we relate the Kummer pairing of the Drinfeld module $\bar\phi$ with that of~$\BDM$. 
Keeping the notation of the earlier sections, recall that $F$ is the quotient field of~$A$, that $D$ is the total ring of quotients of the non-commutative integral domain~$R$, and that $\bar\phi$ is an embedding $A\into R$.

\begin{Lem}\label{RotimesF=D}
There is a natural isomorphism of right $F$-vector spaces
$$\kappa\colon R\otimes_A F\longisoto D,\ f\otimes ab^{-1}\longmapsto f\bar\phi(a)\bar\phi(b)^{-1}.$$
\end{Lem}

\begin{Proof}
Clearly $\kappa$ is well-defined. Since $R$ is a torsion free right $A$-module and $\bar\phi$ is injective, the map $\kappa$ is injective as well. Let $D'$ denote its image. 

Next recall that $\tau^d$ lies in the center of~$R$ and hence in the finitely generated $A$-module $\End_k(\bar\phi)$. Thus the subring $\bar\phi(A)[\tau^d]\subset R$ is a finitely generated $A$-module. On the other hand $R$ is a finitely generated module over the subring $\BDM(B)=\BFp[\tau^d]$ of $\bar\phi(A)[\tau^d]$. Together this implies that $R$ is a finitely generated right $A$-module. Therefore $R\otimes_AF$ and hence $D'$ is a finite dimensional right $F$-vector space.

Now observe that left multiplication by any non-zero $f\in R$ induces an injective endomorphism of the right $A$-module~$R$. It therefore also induces an injective endomorphism of the right $F$-vector space~$D'$. By finite dimensionality this endomorphism must be surjective. For every $x\in D'$ there therefore exists $y\in D'$ with $fy=x$. 

But in $D$ this last equation is equivalent to $y=f^{-1}x$. Thus the subset $D'\subset D$ is invariant under left multiplication by~$f^{-1}$. As it is already a left $R$-submodule, it is therefore a left $D$-submodule of~$D$. Since it also contains the identity element, it is therefore equal to~$D$, and we are done.
\end{Proof}

\medskip
Next observe that the embedding $\bar\phi\colon A\into R$ extends uniquely to a continuous ring homomorphism $A_\ad\into R_\ad := R\otimes_BB_\ad$, turning $\Tuad(\BDM)$ into an $A_\ad$-module.

\begin{Prop}\label{TadToBasic}
There is a natural Galois-equivariant isomorphism of $A_\ad$-modules 
$$\Tuad(\bar\phi)\ \longisoarrow\ \Tuad(\BDM).$$
\end{Prop}

\begin{Proof}
Recall from the proof of Lemma \ref{RotimesF=D} that $R$ is finitely generated as a right $A$-module. Being torsion free, it is thus a projective right $A$-module. Using Lemma \ref{RotimesF=D}, the short exact sequence $0\to A\to F\to F/A\to 0$ therefore induces an isomorphism of short exact sequences of left $R$-modules
$$\xymatrix@R-10pt{
0 \ar[r] & R \ar[r] & R\otimes_AF \ar[r] & R\otimes_AF/A \ar[r] & 0 \\
0 \ar[r] & R \ar[r] \ar@{=}[u] & D \ar[r] \ar@{=}[u]^\wr_\kappa & D/R \ar[r]  \ar@{=}[u]^\wr & 0.}$$
The same holds with $(B,\BDM)$ in place of $(A,\bar\phi)$; hence together we obtain natural isomorphisms of left $R$-modules
\UseTheoremCounterForNextEquation
\begin{equation}\label{TadToBasicCons}
R\otimes_A F/A\ \cong\ D/R\ \cong\ R\otimes_{B}E/B
\end{equation}
with $E=\Quot(B)$. Letting $i$ denote the tautological action of $R$ on~$\kalg$, using adjunction twice we deduce natural isomorphisms
\UseTheoremCounterForNextEquation
\begin{equation}\label{TadToBasicIsos}
\begin{array}{r@{\hspace{1.2ex}}c@{\hspace{1.2ex}}l}%
T_\ad(\bar\phi) & \stackrel{\rm def}{=} & \Hom_A\bigl(F/A,(\kalg,\,\bar\phi)\bigr) \\[.66em]
& \stackrel{{\rm adj}}{=} & \Hom_R\bigl(R\otimes_A F/A,\,(\kalg,i)\bigr) \\[.66em]
& \cong & \Hom_R\bigl(R\otimes_{B}E/B,(\kalg,i)\bigr) \\[.66em]
& \stackrel{{\rm adj}}{=} & \Hom_{B}\bigl(E/B,\,(\kalg,\,\BDM)\bigr)
\qquad\stackrel{\rm def}{=}\ T_\ad(\BDM).
\end{array}
\end{equation}
By construction these isomorphisms are Galois-equivariant. Using Proposition \ref{TateMod1Prop} we obtain natural Galois-equivariant isomorphisms $\Tuad(\bar\phi) \cong T_\ad(\bar\phi) \cong T_\ad(\BDM) \cong \Tuad(\BDM)$, as desired.
\end{Proof}


\begin{Prop}\label{KummerSame}
The following diagram commutes:
$$\xymatrix@C+20pt{
\ (\KOKperf,\bar\phi)\times J_K\ \ar[r]^-{\kum{\ }{\ }{\bar\phi}} \ar@{=}@<+32pt>[d] \ar@{=}@<-14pt>[d] & \ \Tuad(\bar\phi)\ \ar[d]^\wr_{\ref{TadToBasic}} \\
\ (\KOKperf,\BDM)\times J_K\ \ar[r]^-{\kum{\ }{\ }{\BDM}} & \ \Tuad(\BDM) \rlap{.}\  \\
}$$
\end{Prop}

\begin{Proof}
For any $\xi \in \Kperf$ and $[\gal]\in J_K$ we must show that, under the isomorphism from Proposition \ref{TadToBasic}, the element $\kum{\xi}{\gal}{\bar\phi} \in \Tuad(\bar\phi)$ corresponds to the element ${\kum{\xi}{\gal}{\BDM} \in \Tuad(\BDM)}$.
Here $\kum{\xi}{\gal}{\bar\phi}$ is an $A$-linear map $F/A \to (\kalg,\,\bar\phi)$, 
which by the adjunction in \eqref{TadToBasicIsos} corresponds 
to the left $R$-linear map
\newcommand{\kumti}[3]{L_{#3}}%
\UseTheoremCounterForNextEquation
\begin{equation}\label{KummerSame1}
\kumti{\xi}{\gal}{\bar\phi}\colon R\otimes_A F/A \longto (\kalg,i),
\ \ f\otimes[\tfrac{1}{a}]\longmapsto f\bigl(\kum{\xi}{\gal}{\bar\phi}(\tfrac{1}{a})\bigr).
\end{equation}
Similarly the $B$-linear map $\kum{\xi}{\gal}{\BDM}\colon E/B\to(\kalg,\,\BDM)$ corresponds to the left $R$-linear map
\UseTheoremCounterForNextEquation
\begin{equation}\label{KummerSame2}
\kumti{\xi}{\gal}{\BDM}\colon R\otimes_B E/B \longto (\kalg,i),
\ \ f\otimes[\tfrac{1}{b}]\longmapsto f\bigl(\kum{\xi}{\gal}{\BDM}(\tfrac{1}{b})\bigr).
\end{equation}

Consider any element $b\in B\setminus\{0\}$. Then $\BDM_b R$ is a right ideal of finite index in~$R$; hence there exists an element $a\in A\setminus\{0\}$ with $\bar\phi_a\in \BDM_b R$. Choose $g\in R$ such that $\bar\phi_a = \BDM_b g$. Then in the ring $D$ we have $\BDM_b^{-1} = g\bar\phi_a^{-1}$. This means that $g\otimes[\frac{1}{a}]\in R\otimes_AF/A$ equates to $1\otimes[\frac{1}{b}]\in R\otimes_BE/B$ under the isomorphism \eqref{TadToBasicCons}. By construction the isomorphism \eqref{TadToBasicCons} is $R$-linear on the left, so for all $f\in R$ 
it is given by
\UseTheoremCounterForNextEquation
\begin{equation}\label{KummerSame2Bis}
R\otimes_AF/A \longto R\otimes_BE/B,\ \ 
fg\otimes\bigl[\tfrac{1}{a}\bigr] \longmapsto f\otimes\bigl[\tfrac{1}{b}\bigr]
\end{equation}
Now choose an element $\xi_a\in\kalg$ with $\psi_a(\xi_a) = \xi$. Then by Lemma \ref{KummerLem} for $\bar\phi$ in place of $\psi$ we have
\UseTheoremCounterForNextEquation
\begin{equation}\label{KummerSame4}
\kum{\xi}{\gal}{\bar\phi}(\tfrac{1}{a})\ =\ \rlap{$[\gal(\xi_a)-\xi_a].$}
\phantom{[\gal(g(\xi_a))-f(\xi_a)].}
\end{equation}
On the other hand, the equation $\BDM_b g = \bar\phi_a$ implies that $\BDM_b(g(\xi_a)) = \bar\phi_a(\xi_a) = \xi$. By applying Lemma \ref{KummerLem} to $\BDM$ we therefore have 
\UseTheoremCounterForNextEquation
\begin{equation}\label{KummerSame5}
\kum{\xi}{\gal}{\BDM}(\tfrac{1}{b})\ =\ [\gal(g(\xi_a))-g(\xi_a)].
\end{equation}
Combining everything we now deduce that
$$\xymatrix@R-10pt@C-10pt{
\kumti{\xi}{\gal}{\bar\phi}\bigl(fg\otimes[\tfrac{1}{a}]\bigr)\ 
\ar@{=}[r]^-{(\ref{KummerSame1})}
& \ fg\bigl(\kum{\xi}{\gal}{\bar\phi}(\tfrac{1}{a})\bigr)\ 
\ar@{=}[r]^-{(\ref{KummerSame4})}
& \ fg\bigl([\gal(\xi_a)-\xi_a]\bigr)\ \kern10pt \\
\kern7pt \kumti{\xi}{\gal}{\BDM}\bigl(f\otimes[\tfrac{1}{b}]\bigr)\ 
\ar@{=}[r]^-{(\ref{KummerSame2})}
& \kern10pt f\bigl(\kum{\xi}{\gal}{\BDM}(\tfrac{1}{b})\bigr)\ 
\ar@{=}[r]^-{(\ref{KummerSame5})}
& \ f\bigl([\gal(g(\xi_a))-g(\xi_a)]\bigr)\ \kern-14pt \ar@{=}[u]
}$$
The formula \eqref{KummerSame2Bis} now implies that the maps $\kumti{\xi}{\gal}{\bar\phi}$ and $\kumti{\xi}{\gal}{\BDM}$ correspond to each other, as desired.
\end{Proof}

\section{The local Kummer pairing in the general case}
\label{TLKPITGC}

Now we will compare the local Kummer pairings of the Drinfeld modules $\psi$ and $\BDM$
and study the implications.

\begin{Prop}\label{TateDoubleIso}
There are natural isomorphisms of $\Bado$-modules
$$\Tuad(\psi)\ \longisoto\ \Tuad(\bar\phi)\ \longisoto\ \Tuad(\BDM).$$
In particular $\Tuad(\psi)$ is a free $\Bado$-module of rank~$d$.
\end{Prop}

\begin{Proof}
The isomorphisms are those from Propositions \ref{TateMod1Prop} and \ref{TadToBasic}. Both are Galois-equivariant and therefore $\Bado$-linear. But since $\BDM$ is a supersingular Drinfeld $B$-module of rank~$d$, Proposition \ref{TateBado} implies that $\Tuad(\BDM)$ is a free $\Bado$-module of rank~$d$. Thus the last statement follows.
\end{Proof}

\medskip
Let $\chi$ be the natural isomorphism from Proposition \ref{XIsom}.

\begin{Prop}\label{KummerSamePsi}
The following diagram commutes:
$$\xymatrix@C+20pt{
\ (\KOKperf,\psi)\times J_K \ar[r]^-{\kum{\ }{\ }{\psi}} \ \ar@{=}@<+33pt>[d] & \ \Tuad(\psi)\ \ar[d]^\wr_{\ref{TateDoubleIso}} \\
\ (\KOKperf,\BDM)\times J_K\ \ar[r]^-{\kum{\ }{\ }{\BDM}} \ar@<+15pt>[u]^\wr_\chi & \ \Tuad(\BDM) \rlap{.}\  \\
}$$
\end{Prop}

\begin{Proof}
Combine Propositions \ref{KummerRed} and \ref{KummerSame}.
\end{Proof}

\medskip
The next auxiliary result describes how the isomorphism $\chi$ interacts with the filtration $\Wo{\KOKperf}$ in some instances.

\begin{Prop}\label{ChiXiXi}
For any element $\xi\in K$ (sic!) of normalized valuation $-i\le0$ we have 
$$\chi^{-1}([\xi]) - [\xi]\ \in\ \Wi{\KOKperf}.$$
\end{Prop}

\begin{Proof}
Recall from Proposition \ref{XInvPerf} that $x^{-1}=\sum_{j\le0} \tau^jz_j$ with all $z_j\in\CO_K$. By Proposition \ref{XIsom} we therefore have 
$$\chi^{-1}([\xi]) - [\xi]\ =\ \sum_{j\le0} [\tau^j((z_j-\delta_{j,0})\xi)]$$
with $(z_j-\delta_{j,0})\xi\in K$.
Now by Proposition \ref{XExistsPsi} \ref{XExistsPsiBounds} 
we have $v(z_j-\delta_{j,0})>0$. Thus $(z_j-\delta_{j,0})\xi$ has normalized valuation $>-i$.
Its residue class therefore lies in $\Wi\KOKperf$ by Proposition \ref{WiKOKperfGens}, and hence so does its image under~$\tau^j$.
\end{Proof}


\medskip
The following result generalizes Theorem \ref{GrJKIsom} \ref{GrJKIsomIsom} to all Drinfeld modules $\psi$ of good reduction and sharpens Theorem 3.2.4 of \cite{MornevLM}. 

\begin{Thm}\label{GrJKIsomPsi}
Consider any element $\xi\in K$ of normalized valuation $-i<0$.
\begin{enumerate}\StatementLabels%
\item\label{GrJKIsomPsiVanish} The homomorphism $\kum{\xi}{\ }{\psi}$ vanishes on the ramification subgroup $J_K^{i+1}$.
\item\label{GrJKIsomPsiIsom} If in addition $p\nmid i$, then the homomorphism $\kum{\xi}{\ }{\psi}$ induces an isomorphism
$$\gr^iJ_K \longisoarrow \Tuad(\psi).$$
\end{enumerate}
\end{Thm}

\begin{Proof}
Write $\chi^{-1}([\xi])=[\xi']$ for $\xi'\in\Kperf$. First observe that $[\xi]\in \Wii\KOKperf$ by Proposition \ref{WiKOKperfGens}. From Proposition \ref{ChiXiXi} and the fact that $\Wo$ is an increasing filtration we deduce that $[\xi']\in \Wii{\KOKperf}$. By Theorem \ref{AltPerfectAndFiltPerf} this implies that the homomorphism $\kum{\xi'}{\ }{\BDM}$ vanishes on~ $J_K^{i+1}$. By Proposition \ref{KummerSamePsi} it follows that $\kum{\xi}{\ }{\psi}$ vanishes on~$J_K^{i+1}$, proving  \ref{GrJKIsomPsiVanish}.

Next, Proposition \ref{ChiXiXi} also shows that $[\xi'-\xi]\in \Wi{\KOKperf}$. By Theorem \ref{AltPerfectAndFiltPerf} this implies that the homomorphism $\kum{\xi'-\xi}{\ }{\BDM}$ vanishes on~$J_K^i$. Thus the maps $\kum{\xi'}{\ }{\BDM}$ and $\kum{\xi}{\ }{\BDM}$ agree on~$J_K^i$.
Now assume that $p\nmid i$. Then by Theorem \ref{GrJKIsom} \ref{GrJKIsomIsom} the map $\kum{\xi}{\ }{\BDM}$ vanishes on $J_K^{i+1}$ and induces an isomorphism $\gr^iJ_K \isoarrow \Tuad(\BDM)$. The same therefore holds for the map $\kum{\xi'}{\ }{\BDM}$. By Proposition \ref{KummerSamePsi} it follows that $\kum{\xi}{\ }{\psi}$ vanishes on $J_K^{i+1}$ and induces an isomorphism $\gr^iJ_K \isoarrow \Tuad(\psi)$, proving  \ref{GrJKIsomPsiIsom}.
\end{Proof}

\begin{Rem}\label{GrJKIsomPsiRem}
\rm One can use the isomorphism of Theorem \ref{GrJKIsomPsi} \ref{GrJKIsomPsiIsom} to endow the group $\gr^iJ_K$ with a structure of an $A_\ad$-module. However, this structure is uncanonical, because it depends not only on the Drinfeld module~$\psi$, but also on the residue class $[\xi]\in\KOKperf$.
\end{Rem}

\medskip
We can also 
generalize the isomorphy part of Theorem \ref{AltPerfectAndFiltPerf}. As a preparation we need:

\begin{Lem}\label{CompletLem}
There is a natural isomorphism $A_\ad\otimes_AR \cong B_\ad\otimes_BR$.
\end{Lem}

\begin{Proof}
Let $C$ be the subring of $R$ that is generated by $\bar\phi(A)$ and
$\BDM(B)$. Since $\BDM(B)$ lies in the center of~$R$, this ring is commutative.
The fact that $R$ is a finitely generated left module
over $A$ and over~$B$ implies that $C$ is also a finitely generated module
over $A$ and over~$B$. Thus the profinite completion of $R$ as a left
$A$-module coincides with the profinite completion as a left $C$-module
and again with the profinite completion as a left $B$-module. Translated
into tensor products this is just the desired formula.
\end{Proof}

\begin{Thm}\label{AltPerfectPsi}
The local Kummer pairing of $\psi$ is adjoint to an isomorphism of
topological $A_\ad$-modules
$$A_\ad \complot_A (\KOKperf,\,\psi)\ \longisoto\
\Hom_\Bado^\cont\bigl(J_K,\:\Tuad(\psi)\bigr),\quad 1\otimes[\xi]
\longmapsto \kum{\xi}{\ }{\psi}.$$
\end{Thm}

\begin{Proof}
We will construct the following commutative diagram:
$$\begin{tikzcd}[column sep=3.5em, row sep=1.24em]
A_\ad\complot_A (\KOKperf,\,\psi) \rar["\sim"', "\kum{\ }{\ }{\psi}"] \dar[equals,"\wr"]&
\Hom_\Bado^\cont(J_K,\,\Tuad(\psi)) \dar[equals,"\wr"] \\
A_\ad\complot_A (\KOKperf,\,\bar\phi) \rar["\sim"', "\kum{\ }{\ }{\bar\phi}"] \dar[equals,"\wr"]&
\Hom_\Bado^\cont(J_K,\,\Tuad(\bar\phi)) \dar[equals,"\wr"] \\
B_\ad\complot_B (\KOKperf,\,\BDM) \rar["\sim"', "\kum{\ }{\ }{\BDM}"]&
\Hom_\Bado^\cont(J_K,\,\Tuad(\BDM)).
\end{tikzcd}$$

For this we first note that since $R^\circ$ is already a module over the
localization $B^\circ=B[s^{-1}]$, there is a natural isomorphism
$B_\ad\otimes_BR^\circ \cong \Bado\otimes_BR^\circ$. Since $\KOKperf$ is a
free $R^\circ$-module, this yields a natural isomorphism
$B_\ad\otimes_B(\KOKperf,\BDM) \cong \Bado\otimes_B(\KOKperf,\BDM)$.
Theorem \ref{AltPerfectAndFiltPerf} thus induces the lower horizontal
isomorphism.

Next, applying ${(\ )\complot_R \KOKperf}$ to the isomorphism from Lemma
\ref{CompletLem} yields the vertical isomorphism at the lower left. The
vertical isomorphism at the lower right and the commutativity of the lower
half result by adjunction and continuity from Proposition \ref{KummerSame}.

Finally, the upper half of the commutative diagram results by adjunction
from Proposition \ref{KummerRed}. The upper edge now yields the theorem.
\end{Proof}

\begin{Rem}\label{AltPerfectPsiRem}
\rm By construction, the isomorphism in Theorem \ref{AltPerfectPsi}
identifies the submodule $A_\ad\complot_A (\chi(\Wi{\KOKperf}),\psi)$ with
the submodule of homomorphisms that vanish on the ramification
subgroup~$J_K^i$. To get a full analog of Theorem \ref{AltPerfectAndFiltPerf} one should also supply a
description of the modules $\chi(\Wi{\KOKperf})$ in terms of~$\psi$. At
the moment we do not know how to do this. A partial result in this
direction is Proposition \ref{ChiXiXi}.
\end{Rem}

\section{The image of inertia in the general case}
\label{MainRes}

Now let $M\subset (\Kperf,\psi)$ be a finitely generated $A$-submodule satisfying $M\cap\CO_\Kperf=\{0\}$. We want to determine the image of the homomorphism
\UseTheoremCounterForNextEquation
\begin{equation}\label{RhoMFactors}
\rho_M\colon J_K\ \longto\ \Hom_A\bigl(M,\:\Tuad(\psi)\bigr),\quad\gal\ \longmapsto\ \kum{\gal}{\ }{\psi}.
\end{equation}
By Proposition \ref{KummerM} this includes the homomorphism \eqref{InertiaTadMJ} as a special case. The homomorphism $\rho_M$ is $\Bado$-linear by Proposition \ref{KummerProp2}.

\medskip
Since $M\cap\CO_\Kperf=\{0\}$, we may identify $M$ with its image in 
$\KOKperf=\Kperf/\CO_\Kperf$.
Here $(\KOKperf,\psi)$ is a torsion-free $A$-module by the same argument as in the proof of Proposition \ref{TateMod1Prop}.
Thus the $A$-module $M$ is finitely generated torsion free and therefore projective.
Via the natural isomorphism $\chi$ of Proposition \ref{XIsom} it
corresponds to the submodule $\chi^{-1}(M)$ of $(\KOKperf,\bar\phi)$. Let $\Mbar$ be the left $R$-submodule of $\KOKperf$ that is generated by $\chi^{-1}(M)$. 

\begin{Prop}\label{MMbarRank}
The left $R$-module $\Mbar$ is free, and we have  $\rank_R(\Mbar) \le \rank_A(M)$. 
\end{Prop}

\begin{Proof}
By construction $\Mbar$ is finitely generated, so it is free by Lemma \ref{BreaksMbar} \ref{BreaksMbarRank}. The inequality of ranks follows from the surjection $R\otimes_AM\onto\Mbar$, $f \otimes m \mapsto f(\chi^{-1}(m))$. 
\end{Proof}



\begin{Prop}\label{MainImageDiag}
There is a natural commutative diagram of $\Bado$-linear maps
$$\xymatrix@R+10pt{J_K \ar[r]^-{\rho_M} \ar[d]_{\rho_\Mbar}^{(\ref{RhoMbarFactors})} & 
\Hom_A\bigl(M,\:\Tuad(\psi)\bigr) \ar[d]^\wr_{\ref{TateDoubleIso}} \\
\Hom_R\bigl(\Mbar,\:\Tuad(\BDM)\bigr) \ar@{^{ (}->}[r] & 
\Hom_A\bigl(M,\:\Tuad(\BDM)\bigr)\rlap{,} \\}$$
where the lower horizontal map is induced by the embedding $\chi^{-1}\colon M\into\Mbar$. Furthermore, the image $\rho_\Mbar(J_K)$ has finite index in $\Hom_R(\Mbar,\,\Tuad(\BDM))$.
\end{Prop}

\begin{Proof}
Since $\rho_M$ and $\rho_\Mbar$ are given by the respective Kummer pairings, the specified commutative diagram is obtained by adjunction from the diagram of Proposition~\ref{KummerSamePsi}.
The last statement is a repetition of Theorem~\ref{Mammage}.
\end{Proof}


\begin{Thm}\label{MammagePsi}
The image $\rho_M(J_K)$ is a free $\Bado$-module of rank $d\cdot\rank_R(\Mbar)$. Moreover $\rho_M(J_K)$ has finite index in a $\Bado$-module direct summand of $\Hom_A(M,\,\Tuad(\psi))$.
\end{Thm}

\begin{Proof}
The left $R$-module $\Mbar$ is free by Proposition \ref{MMbarRank}; hence the canonical surjection $R\otimes_A M \onto \Mbar$ possesses a right inverse. The image of the induced homomorphism $\Hom_R(\Mbar,\Tuad(\BDM)) \longto \Hom_R(R\otimes_AM,\Tuad(\BDM))$ is therefore a direct summand. But by adjunction this homomorphism is just the lower homomorphism in Proposition \ref{MainImageDiag}. Both statements thus follow from Theorem \ref{Mammage}.
\end{Proof}


\begin{Thm}\label{MammageRamPsi}
For all $i\gg0$ we have $\rho_M(J_K^i) = 0$.
\end{Thm}

\begin{Proof}
Combine Proposition \ref{MainImageDiag} with Theorem \ref{MammageRam}.
\end{Proof}

\medskip
Let $S_\Mbar$ be the set of breaks of the module $\Mbar$ as in \eqref{SMDef}.

\begin{Thm}\label{RammagePsi}
For every $i\ge0$ we have:
\begin{enumerate}\StatementLabels%
\item\label{RammagePsiNot} If $i\not\in S_\Mbar$, then $\rho_M(J_K^i)/\rho_M(J_K^{i+1})$ is a finite $\Bado$-module.
\item\label{RammagePsiIn} If $i\in S_\Mbar$, then $\rho_M(J_K^i)/\rho_M(J_K^{i+1})$ is a free $\Bado$-module of rank~$d$.
\item\label{RammagePsiRank} The image $\rho_M(J_K^i)$ is a free $\Bado$-module of rank $d\cdot\bigl|\{j\in S_\Mbar\mid j\ge i\}\bigr|.$
\end{enumerate}
\end{Thm}

\begin{Proof}
Combine Proposition \ref{MainImageDiag} with Theorem \ref{Rammage}.
\end{Proof}


\begin{Cor}\label{MammagePsiOpen}
\begin{enumerate}\StatementLabels%
\item\label{MammagePsiOpenCrit} The image $\rho_M(J_K)$ is open if and only if $\rank_A(M)=\rank_R(\Mbar)$.
\item\label{MammagePsiOpenOne} The image $\rho_M(J_K)$ is open if $\rank_A(M) = 1$.
\end{enumerate}
\end{Cor}

\begin{Proof}
Proposition \ref{MainImageDiag} implies that
$\rho_M(J_K)$ is open in $\Hom_A(M,\Tuad(\psi))$ if and only if the embedding $\Hom_R(\Mbar,\Tuad(\BDM)) \into \Hom_A(M,\Tuad(\BDM))$ has finite index. By adjunction this map corresponds to the embedding $\Hom_R(\Mbar,\Tuad(\BDM)) \into \Hom_R(R\otimes_AM,\Tuad(\BDM))$ induced by  the canonical surjection $R\otimes_AM\onto\Mbar$. Thus $\rho_M(J_K)$ is open if and only if the two modules have the same rank. This shows \ref{MammagePsiOpenCrit}, and \ref{MammagePsiOpenOne} is a direct consequence thereof.
\end{Proof}

\begin{Rem}\label{MornevLMRem1}
\rm Theorem \ref{MammageRamPsi} implies the central result of the article \cite{MornevLM}:
For every Drinfeld module $\phi$ over $K$ of finite residual characteristic $\Fpres$ 
there is a rational number $\nu\ge0$ such that the inertia subgroup $I_K^\nu$
acts trivially on the Tate modules $T_\Fp(\phi)$ with $\Fp\ne\Fpres$.
Indeed, the argument of \cite[Theorem 4.2.1]{MornevLM} reduces one to the 
case that $\phi$ has stable reduction and its period lattice $M$ is contained in $\Kperf$.
In this situation
Theorem \ref{MammageRamPsi} applies.
\end{Rem}

\begin{Rem}\label{MornevLMRem2}
\rm Corollary \ref{MammagePsiOpen} \ref{MammagePsiOpenOne} implies that for each $\xi\in\Kperf\setminus\CO_\Kperf$ the homomorphism $\kum{\xi}{\ }{\psi}\colon J_K \to \Tuad(\psi)$ has open image. This extends Theorem 3.3.2 of \cite{MornevLM} to the prime~$\Fpres$.

Corollary \ref{MammagePsiOpen} also implies Theorem 3 of \cite{MornevLM}, which applies to any Drinfeld module $\phi$ of rank $r > 1$ and finite residual characteristic $\Fpres$ over $K$ whose period lattice $M$ has $A$-rank~$1$. It shows that for every prime $\Fp\ne\Fpres$ there is a closed subgroupscheme $U\cong \mathbb{G}_a^{\times(r-1)}$ of 
$\underline{\rm Aut}(T_\Fp(\phi))$ 
such that the image of $I_K$ is commensurable to $U(A_\Fp)$.
\end{Rem}

\begin{Rem}
\rm There does not seem to be a way of computing the image $\rho_M(J_K)$ up to finite index directly in terms of $\psi$ and~$M$. The best that we can do is to take generators $f_i = \sum_j f_{ij}\otimes m_{ij}$ of the kernel of the map $R\otimes_AM \onto \Mbar$, $f\otimes m\mapsto f(\chi^{-1}(m))$ and to let each $f_{ij}\in R$ act on $\Tuad(\psi)$ through its action on $\Tuad(\BDM)$ via the isomorphisms (\ref{TateDoubleIso}). Then $\rho_M(J_K)$ is a free $\Bado$-submodule of finite index in the $\Bado$-module
$$\textstyle \bigl\{\ell\in\Hom_A\bigl(M,\:\Tuad(\psi)\bigr) \bigm| 
\forall i\colon \sum_j f_{ij}(\ell(m_{ij}))=0 \bigr\}.$$
\end{Rem}

\medskip
The following fact will be useful:

\begin{Prop}\label{MammagePsiMMPrime}
For any $A$-submodule $M'\subset M$ there is a natural surjection
$$\rho_M(J_K)\ \longonto\ \rho_{M'}(J_K).$$
This is an isomorphism if and only if every element of $\chi^{-1}(M)$ is $R$-linearly dependent on $\chi^{-1}(M')$.
\end{Prop}

\begin{Proof}
Let $\Mbar{}'$ be the left $R$-submodule of $\Mbar$ that is generated by $\chi^{-1}(M')$. Then the inclusions $M'\subset M$ and $\Mbar{}'\subset\Mbar$ and the commutative diagram from Proposition \ref{MainImageDiag} induce a commutative diagram
$$\vcenter{\vskip-3pt\xymatrix@C=20pt@R=2pt{
&& \text{\large\strut}\Hom_R\bigl(\Mbar,\:\Tuad(\BDM)\bigr) \ar[dd] \ar@{^{ (}->}[r] & \Hom_A\bigl(M,\:\Tuad(\psi)\bigr) \ar[dd] \\
J_K\  \ar@<5pt>[urr]^-{\rho_\Mbar} \ar@<-5pt>@{->}[drr]^-{\rho_{\Mbar{}'}} &&& \\
&& \text{\large\strut}\Hom_R\bigl(\Mbar{}',\:\Tuad(\BDM)\bigr) \ar@{^{ (}->}[r] & \Hom_A\bigl(M',\:\Tuad(\psi)\bigr)\rlap{.} \\}
\vskip-3pt}$$
This directly yields a natural surjection $\rho_M(J_K) \onto \rho_{M'}(J_K)$. Moreover, the two arrows at the left have open image by Theorem \ref{Mammage}. Thus the surjection is an isomorphism if and only if the vertical arrow in the middle is injective. As $\Mbar$ and $\Mbar{}'$ are free left $R$-modules of finite rank, this is so if and only if their ranks are equal, which is equivalent to the stated condition.
\end{Proof}


\begin{Caut}\label{NotAadSubmodCaut}
\rm In general the images $\rho_M(J_K)$ and $\rho_M(J_K^i)$ are not $A_\ad$-submodules, not even up to finite index, as the following example shows:
\end{Caut}

\begin{Ex}\label{NotAadSubmodEx}
\rm Let $\psi$ be the Drinfeld $A$-module of rank $1$ with $A=\BFp[t]$ and $\psi_t=\tau$. Then $\bar\phi=\psi$; hence $\chi$ is the identity on $\KOKperf$. Assume that $|k|=p^2$, so that $B=\BFp[s]$ with $\BDM_s=\tau^2$. Identifying $s=t^2$ then makes $B$ a subring of $A$ such that $\BDM=\psi|B$. 
As in Section \ref{Filtr} let $\xi_1\in K$ be an element of normalized valuation~$-1$. Choose an element $\alpha\in k\setminus\BFp$ and let $M$ be the free $A$-submodule of rank $2$ of $(\KOKperf,\psi)$ that is generated by $\xi_1$ and $\alpha\xi_1$. This is then also the free $R$-submodule of rank $1$ of $\KOKperf$ that is generated by~$\xi_1$ and hence equal to~$\Mbar$. Evaluation at $\xi_1$, respectively at $\xi_1$ and $\alpha\xi_1$, thus induces a commutative diagram
$$\def\scriptmapsto{{\begin{turn}{-90}$\scriptstyle\mapsto$\end{turn}}}
\vcenter{\vskip-3pt\xymatrix@C=23pt@R=22pt{
&& \text{\large\strut}\Hom_R\bigl(\Mbar,\:\Tuad(\BDM)\bigr) 
\ar@{^{ (}->}[d] \ar[rrr]^-\sim_-{h\kern1pt\mapsto h(\xi_1)} 
&&& \Tuad(\BDM)\text{\Large\mathstrut}
\ar@{^{ (}->}[d]_{{\raisebox{-3pt}{$\scriptstyle x$}\atop\scriptmapsto\kern1pt}\atop\scriptstyle(x,\alpha x)} \\
J_K\  \ar@{->>}[urr]!<6pt,17pt>^-{\rho_\Mbar} \ar[drr]!<6pt,-17pt>^-{\rho_M} 
&& \text{\large\strut}\Hom_A\bigl(M,\:\Tuad(\BDM)\bigr) 
\ar[rrr]^-\sim_-{h\kern1pt\mapsto (h(\xi_1),h(\alpha\xi_1))} \ar@{=}[d]^\wr
&&& \Tuad(\BDM) \rlap{${}^{\oplus2}$} \ar@{=}[d]^\wr \\
&& \text{\large\strut}\Hom_A\bigl(M,\:\Tuad(\psi)\bigr) 
\ar[rrr]^-\sim_-{h\kern1pt\mapsto (h(\xi_1),h(\alpha\xi_1))} 
&&& \Tuad(\psi)\rlap{${}^{\oplus2}$.} \\}
\vskip-3pt}$$
Here the upper left arrow $\rho_\Mbar$ is surjective by Theorem \ref{GrJKIsom}. 
The image of $J_K$ in $\Tuad(\BDM)^{\oplus2}$ is therefore the graph of multiplication by~$\alpha$. Since $\alpha$ does not commute with $\psi(A)=\BFp[\tau]$ within $R=k[\tau]$, this graph is not a $\psi(A)$-submodule. It follows that $\rho_M(J_K)$ is not an $A_\ad$-submodule of $\Hom_A(M,\Tuad(\psi))$.

Note also that in this case $\Tuad(\BDM)$ is a free $\Bado$-module of rank~$2$. Thus $\Hom_A(M,\Tuad(\psi))$ is a free $\Bado$-module of rank $4$ and $\rho_M(J_K)$ is a free $\Bado$-submodule of rank~$2$.
\end{Ex}

\section{Some simple criteria for the image of inertia}
\label{OnlyM}

The main results of Section \ref{MainRes}
concern the image of inertia under~$\rho_M$, but require the knowledge of the left $R$-module~$\Mbar$. In the present section we explore to some extent what one can say about $\rho_M$ without computing~$\Mbar$.

\medskip
First we work out a lower bound on $\rank_R(\Mbar)$ in terms of the $A$-module $M$ and deduce a sufficient condition for $\rho_M(J_K)$ to be open.
Recall that the valuation $v$ is normalized on~$K$.
For any $\xi\in{\Kperf\setminus\CO_\Kperf}$ the value $v(\xi)$ is a negative rational number with at most a power of $p$ in the denominator. We write it in the form $v(\xi)=-p^\nu j$ with unique integers $p\nmid j>0$ and~$\nu$. These integers depend only on the residue class $m = [\xi] \in \KOKperf$ and we set $j(m) := j$. 
(Caution: The value $j(m)$ is only vaguely related to the smallest $i$ such that $m\in \Wi\KOKperf$.)

\begin{Lem}\label{iLemma1}
For all non-zero $m\in\KOKperf$ we have $j(\chi^{-1}(m)) = j(m)$.
\end{Lem}

\begin{Proof}
Write $m=[\xi]$ for $\xi\in\Kperf$ with $v(\xi)<0$, and recall from Propositions  \ref{XExistsPsi} and \ref{XIsom} that $\chi^{-1}([\xi]) = \smash{\sum_{\nu\le0} [y_\nu \xi^{p^\nu}]}$ for elements $y_\nu\in\CO_\Kperf$ with $v(y_0)=0$. For all $\nu<0$ we thus have $\smash{v(y_\nu \xi^{p^\nu})} \ge p^\nu\cdot v(\xi) > v(\xi)$. This implies that $\chi^{-1}([\xi]) = [\xi']$ with $v(\xi')=v(\xi)$. Taking prime-to-$p$ parts shows that $j([\xi']) = j([\xi])$, as desired.
\end{Proof}

\begin{Lem}\label{iLemma2}
For all non-zero $m\in\KOKperf$ and $f\in R\setminus\{0\}$ we have ${j(f(m)) = j(m)}$.
\end{Lem}


\begin{Proof}
Write $m=[\xi]$ for $\xi\in\Kperf$ with $v(\xi)<0$, and write $f=\smash{\sum_{\nu=0}^n x_\nu\tau^\nu}$ with $x_\nu\in k$ and $x_n\not=0$. Then we have $f(\xi) = \smash{\sum_{\nu=0}^n x_\nu\xi^{p^\nu}}$ with $\smash{v(x_n\xi^{p^n})} = p^n\cdot v(\xi) < p^\nu\cdot v(\xi) \le \smash{v(x_\nu\xi^{p^\nu})}$ for all $\nu<n$. This implies that $v(f(\xi))=p^n\cdot v(\xi)$. Taking prime-to-$p$ parts on both sides shows that $j(f(m)) = j(m)$, as desired.
\end{Proof}

\begin{Lem}\label{iLemma3}
Any non-zero elements $m_1,\ldots,m_n\in\KOKperf$ with pairwise distinct values $j(m_1),\ldots,j(m_n)$ are $R$-linearly independent.
\end{Lem}


\begin{Proof}
If not, there exist $f_1,\dotsc,f_n\in R$, not all zero, such that $\sum_{\nu=1}^n f_\nu(m_\nu)=0$. After removing irrelevant terms we may assume that $f_\nu\not=0$ for all~$\nu$. By Lemma \ref{iLemma2} we then have $j(f_\nu(m_\nu)) = j(m_\nu)$. Write $f_\nu(m_\nu) = [\eta_\nu]$ with $\eta_\nu\in\Kperf$ and $v(\eta_\nu)<0$. Then as the numbers $j([\eta_\nu])$ are pairwise distinct, their definition implies that the valuations $v(\eta_\nu)$ are pairwise distinct. Thus $v\bigl(\sum_{\nu=1}^n\eta_\nu\bigr)$ is the smallest of these values and therefore $<0$. But this implies that $\sum_{\nu=1}^n f_\nu(m_\nu)\not=0$ and yields a contradiction.
\end{Proof}

\medskip
Now we return to the setting of Section \ref{MainRes}. Consider a finitely generated $A$-submodule $M\subset (\KOKperf,\,\psi)$ and set $j(M) := \{ j(m) \mid m \in M\setminus\{0\}\}$. As in Section \ref{MainRes} let $\Mbar$ be the left $R$-submodule of $\KOKperf$ that is generated by $\chi^{-1}(M)$. 
(Caution: The set $j(M)$ is not directly related to the set $S_\Mbar$ from (\ref{SMDef}).)

\begin{Prop}\label{iRankBound}
We have $|j(M)| \le \rank_R(\Mbar) \le \rank_A(M)$.
\end{Prop}

\begin{Proof}
Consider any non-zero elements $m_1,\ldots,m_n\in M$ with pairwise distinct values $j(m_1),\ldots,j(m_n)$. Then $j(\chi^{-1}(m_1)),\ldots,j(\chi^{-1}(m_n))$ are pairwise distinct by Lemma \ref{iLemma1}. Therefore $\chi^{-1}(m_1),\ldots,\chi^{-1}(m_n)$ are $R$-linearly independent by Lemma \ref{iLemma3}, and hence ${n\le\rank_R(\Mbar)}$. Varying $n$ we deduce the first inequality. The second is a mere repetition of Proposition \ref{MMbarRank}.
\end{Proof}

\begin{Thm}\label{iOpenness}
If $|j(M)| = \rank_A(M)$, then $\rho_M(J_K)$ is open in $\Hom_A(M,\,\Tuad(\psi))$.
\end{Thm}

\begin{Proof}
By Proposition \ref{iRankBound} the assumption implies that $\rank_R(\Mbar) = \rank_A(M)$, so the conclusion follows from Corollary \ref{MammagePsiOpen} \ref{MammagePsiOpenCrit}.
\end{Proof}

\begin{Rem}\label{iRankBoundEx}
\rm The first inequality in Proposition \ref{iRankBound} can be strict. Suppose for instance that $A=\BFp[t]$ with $p>2$ and $\psi_t=\bar\phi_t=\tau^2$, and let $M$ be the $A$-submodule generated by the elements $[\pi^{-1}]$ and $[\pi^{-2}+\pi^{-p}]$ for a uniformizer $\pi\in\CO_K$. Then $j(M)=\{1\}$, while $\Mbar$ is free of rank $2$ over $R$ with basis $[\pi^{-1}]$ and $[\pi^{-2}]$.
\end{Rem}

\medskip
The last result of this section detects certain jumps in the image of the ramification filtration:

\begin{Thm}\label{iMJump}
If there exists $[\xi]\in M$ such that $\xi$ lies in $K$ and has normalized valuation $-i$ with $p\nmid i>0$, then the subquotient $\rho_M(J_K^i)/\rho_M(J_K^{i+1})$ is a free $\Bado$-module of rank~$d$.
\end{Thm}

\begin{Proof}
By assumption $[\xi]$ is an element of $\Wii\KOKperf$ but not of $\Wi\KOKperf$. On the other hand the assumptions imply that ${\chi^{-1}([\xi]) - [\xi]}$ lies in $\Wi{\KOKperf}$ by Proposition \ref{ChiXiXi}. Thus $\chi^{-1}([\xi])$ is an element of $\Wii{\KOKperf}$ but not of $\Wi{\KOKperf}$.
From (\ref{WiMDef}) and (\ref{SMDef}) we deduce that $i\in S_\Mbar$. The conclusion is now a direct consequence of Theorem \ref{RammagePsi} \ref{RammagePsiIn}.
\end{Proof}

\section{Consequences for the \texorpdfstring{$\Fp$}{p}-adic Tate module}
\label{PAdic}

In this last section we examine the image of inertia at a single prime $\Fp$ of~$A$. For this consider the $\Bado$-linear homomorphism
\UseTheoremCounterForNextEquation
\begin{equation}\label{InertiaTadMJFp}
\rho_{M,\Fp}\colon J_K \longto \Hom_A\bigl(M,\:\Tup(\psi)\bigr)
\end{equation}
obtained by combining the homomorphism $\rho_M$ from \eqref{RhoMFactors} with the projection $\Tuad(\psi)\onto \Tup(\psi)$. 
We will discuss the image $\rho_{M,\Fp}(J_K)$, then look at the image of the ramification filtration, and end with the dual properties of the local Kummer pairing.

\medskip
First, from Corollary \ref{MammagePsiOpen} we directly deduce:

\begin{Prop}\label{MammagePsiOpenP}
The image $\rho_{M,\Fp}(J_K)$ is open if $\rank_A(M)=\rank_R(\Mbar)$. In particular it is open if $\rank_A(M) = 1$.
\end{Prop}


To say a little more recall that, since $\psi$ has good reduction, the action of $\GK$ on $\Tuad(\psi)$ factors through~$\Gk$. Let $\Bop$ denote the closure of the subring of $\End(\Tup(\psi))$ that is generated by the Frobenius of~$k$. By Proposition \ref{TateDoubleIso} this is just the image of the composite homomorphism
$$\Bado\ \longinto\ \End(\Tuad(\BDM))\ \cong\ \End(\Tuad(\psi))\ \longonto\ \End(\Tup(\psi)).$$
One should keep in mind that in general $\Bop$ is not an $A_\Fp$-subalgebra, except in the special case $\bar\phi(A)\subset \BDM(B)$. To give an explicit description of this ring, consider the ring homomorphisms
$$\vcenter{\vskip-5pt\hbox{$\xymatrix@C-10pt@R=-7pt{
A\ \ar@{^{ (}->}[drr]^-{\bar\phi} & \\
&& \ \bar\phi(A)[\tau^d] \subset \End_k(\bar\phi) \\ 
B\ \ar@{^{ (}->}[urr]^-{\BDM} & \\}$}\vskip-0pt}.$$
Let $C$ denote the integral closure of $\bar\phi(A)[\tau^d]$ in the (commutative!) subfield $F(\tau^d) \subset \End_k(\bar\phi)\otimes_AF$. Let $Q^\circ_\Fp$ be the finite set of primes $\Fq\not=\Fqres$ of $B$ for which there exists a prime of $C$ that lies over $\Fp$ and over~$\Fq$.

\begin{Prop}\label{BCircFpDecomp} 
The decomposition $\Bado\cong \prod_{\Fq\not=\Fqres}B_\Fq$ induces an isomorphism
$$\Bop\ \cong {\prod_{\smash{\Fq\in Q^\circ_\Fp}} \! B_\Fq}.$$
\end{Prop}

\begin{Proof}
By Proposition \ref{TateMod1Prop} we have natural isomorphisms $\Tup(\psi) \cong \Tup(\bar\phi) \cong T_\Fp(\bar\phi)$. Thus it suffices to prove the proposition with $\bar\phi$ in place of~$\psi$, and with $T_\Fp$ in place of~$\Tup$.

Next, by \cite[Prop.\,4.3]{DevicPink} there exist a Drinfeld $C$-module $\tilde\phi$ and an isogeny $\bar\phi \to \tilde\phi|A$, both defined over~$k$.
This isogeny induces an embedding of finite index $T_\Fp(\bar\phi) \into T_\Fp(\tilde\phi|A)$ which is both $A_\Fp$- and $B$-linear. Replacing $\bar\phi$ by $\tilde\phi|A$ therefore does not change $\Bop$ up to isomorphism, and so it suffices to prove the proposition under the assumption $\bar\phi=\tilde\phi|A$.

In that case $T_\Fp(\bar\phi)$ is naturally isomorphic to the product of the Tate modules $T_\Fr(\tilde\phi)$ for all primes $\Fr$ of $C$ that lie over~$\Fp$. We write this decomposition in the form
\UseTheoremCounterForNextEquation
\begin{equation}\label{TateDecompACB}
T_\Fp(\bar\phi) = \prod_{\smash{\Fq}} T_{\Fp,\Fq}
\qquad\hbox{with}\qquad
T_{\Fp,\Fq} \cong  \! \prod_{\smash{\Fr|\Fp,\,\Fr|\Fq}} \! T_\Fr(\tilde\phi),
\end{equation}
where $\Fq$ runs over all primes of~$B$. Here $T_\Fr(\tilde\phi)$ is a free $C_\Fr$-module of finite rank and non-zero for any prime $\Fr$ of $C$ that is different from the characteristic ideal $\Frres$ of~$\tilde\phi$. Likewise the equality $\BDM=\tilde\phi|B$ implies that for every prime $\Fq$ of $B$ the Tate module $T_\Fq(\BDM)$ decomposes as the product of the Tate modules $T_\Fr(\tilde\phi)$ for all primes $\Fr$ of $C$ that lie over~$\Fq$. 

Applying this to $\Fq=\Fqres$, the fact that  $T_\Fqres(\BDM)=0$ implies that $\Frres$ is the only prime of $C$ over~$\Fqres$ and that $T_{\Frres}(\tilde\phi)=0$. In particular this shows that $T_{\Fp,\Fqres}=0$. By the definition of $Q^\circ_\Fp$ we therefore have $T_{\Fp,\Fq}=0$ for all $\Fq\not\in Q^\circ_\Fp$. 

By contrast, for any prime $\Fr$ of $C$ over $\Fq\not=\Fqres$ the Tate module $T_\Fr(\tilde\phi)$ is a non-zero free $C_\Fr$-module of finite rank and thus a non-zero free $B_\Fq$-module of finite rank. Varying $\Fr$ thus implies that $T_{\Fp,\Fq}$ is a non-zero free $B_\Fq$-module of finite rank for every $\Fq\in Q^\circ_\Fp$.
In view of the factorization $\Bado \cong \prod_{\Fq\not=\Fqres}B_\Fq$ the decomposition (\ref{TateDecompACB}) now implies the proposition.
\end{Proof}

\begin{Prop}\label{BCircFpT}
Under the decomposition in Proposition \ref{BCircFpDecomp}, there is an isomorphism of $\Bop$-modules
$$\Tup(\psi)\ \cong \prod_{\smash{\Fq\in Q^\circ_\Fp}} B_\Fq^{d_{\Fp,\Fq}}
\qquad\hbox{for}\qquad
d_{\Fp,\Fq}\ =\ \smash{\frac{d}{[C/B]}}\cdot\! \sum_{\smash{\Fr|\Fp,\,\Fr|\Fq}} [C_\Fr/B_\Fq],$$
where $[\ /\ ]$ denotes the extension degree and $\Fr$ runs through all primes of $C$ that lie simultaneously over $\Fp$ and over $\Fq$.
\end{Prop}

\begin{Proof}
First assume that  $\bar\phi=\tilde\phi|A$ as in the proof of Proposition \ref{BCircFpDecomp}. Then the fact that $\BDM$ is a Drinfeld $B$-module of rank $d$ implies that $\tilde\phi$ is a Drinfeld $C$-module of rank $d/e$ for $e:=[C/B]$. For each prime $\Fr\not=\Frres$ of $C$ the Tate module $T_\Fr(\tilde\phi)$ is therefore a free $C_\Fr$-module of rank $d/e$. For $\Fr|\Fq$ it is then a free $B_\Fq$-module of rank $d/e\cdot[C_\Fr/B_\Fq]$. With the decomposition (\ref{TateDecompACB}) this shows that $T_{\Fp,\Fq}$ is a free $B_\Fq$-module of rank $d_{\Fp,\Fq}$ and the statement follows.

In the general case we have an embedding of finite index $\Tup(\psi) \cong T_\Fp(\bar\phi) \into T_\Fp(\tilde\phi|A)$, as in the proof of Proposition \ref{BCircFpDecomp}. Since every $B_\Fq$-submodule of finite index of $\smash{B_\Fq^{d_{\Fp,\Fq}}}$ is again isomorphic to $\smash{B_\Fq^{d_{\Fp,\Fq}}}$, the statement follows in general.
\end{Proof}

\begin{Prop}\label{BCircFpMT}
Any $\Bop$-submodule of $\Hom_A(M,\,\Tup(\psi))$ is isomorphic to
$$\prod_{\smash{\Fq\in Q^\circ_\Fp}} B_\Fq^{n_{\Fp,\Fq}}$$
for unique integers $n_{\Fp,\Fq}\ge0$.
\end{Prop}

\begin{Proof}
Since $M$ is a finitely generated projective $A$-module, Proposition \ref{BCircFpT} implies that $\Hom_A(M,\Tup(\psi))$ is a finitely generated projective $\Bop$-module. As each $B_\Fq$ is a principal ideal domain, the claim follows.
\end{Proof}

\medskip
In general the images $\rho_{M,\Fp}(J_K)$ and $\rho_{M,\Fp}(J_K^i)$ are $\Bado$-sumbodules 
of $\Hom_A(M,\Tup(\psi))$ and therefore $\Bop$-submodules. Determining their structure thus amounts to determining the associated numbers $n_{\Fp,\Fq}$. They do not, however, obey a nice general formula; in particular they do not satisfy a simple kind of independence of~$\Fp$. The following results describe some phenomena that may occur.

\begin{Prop}\label{NPQOpen}
If $\rank_A(M)=\rank_R(\Mbar)$, then for any $\Fp$ and $\Fq\in Q^\circ_\Fp$ the number $n_{\Fp,\Fq}$ for $\rho_{M,\Fp}(J_K)$ satisfies
$$n_{\Fp,\Fq}\ =\ d_{\Fp,\Fq}\cdot\rank_R(\Mbar).$$
\end{Prop}

\begin{Proof}
Under this assumption the image $\rho_{M,\Fp}(J_K)$ is open in $\Hom_A(M,\Tup(\psi))$ by Proposition \ref{MammagePsiOpenP}, and the latter is isomorphic to a direct sum of $\rank_A(M)$ copies of $\Tup(\psi)$. The statement thus directly follows from Proposition \ref{BCircFpT}.
\end{Proof}

\begin{Prop}\label{BQFree}
If $\chi^{-1}(M)$ has finite index in~$\Mbar$, then $\rho_{M,\Fp}(J_K)$ is a free $\Bop$-module of rank 
$$d\cdot\rank_R(\Mbar).$$
\end{Prop}

\begin{Proof}
Set $n:=\rank_R(\Mbar)$. Then the assumption implies that there exist $R$-linearly independent elements $\mbar_1,\ldots,\mbar_n\in\Mbar$ such that each $\Mbar_i := R\mbar_i$ is contained in $\chi^{-1}(M)$. Setting $M_i := \chi(\Mbar_i)\subset M$, we then have a commutative diagram
 $$\def\scriptmapsto{{\begin{turn}{-90}$\scriptstyle\mapsto$\end{turn}}}
 \vcenter{\vskip-3pt\xymatrix@C=23pt@R=17pt{
 &&\text{\large\strut}\Hom_R\bigl(\Mbar,\Tuad(\BDM)\bigr)\ 
 \ar@{^{ (}->}[d] \ar@{^{ (}->}[rr]
 &&\ \text{\large\strut}\smash{\bigoplus\limits_{i=1}^n \Hom_R\bigl(\Mbar_i,\Tuad(\BDM)\bigr)}\ 
 \ar@{^{ (}->}[d] \\
 J_K\  \ar@{->>}[urr]!<6pt,17pt>^-{\rho_\Mbar} \ar[drr]!<6pt,-17pt>_-{\rho_{M,\Fp}} 
 &&\text{\large\strut}\Hom_A\bigl(M,\Tuad(\psi)\bigr)\ 
 \ar@{->>}[d] \ar@{^{ (}->}[rr]
 &&\ \text{\large\strut}\smash{\bigoplus\limits_{i=1}^n \Hom_A\bigl(M_i,\Tuad(\psi)\bigr)}\ 
 \ar@{->>}[d] \\
 &&\text{\large\strut}\Hom_A\bigl(M,\Tup(\psi)\bigr)\ 
 \ar@{^{ (}->}[rr]
 &&\ \text{\large\strut}\smash{\bigoplus\limits_{i=1}^n \Hom_A\bigl(M_i,\Tup(\psi)\bigr)}\ \\}}$$
 where all horizontal arrows are inclusions of finite index. By Theorem \ref{Mammage} the image of $J_K$ in $\Hom_R(\Mbar,\Tuad(\BDM))$ has finite index. Up to finite index it is therefore the direct product of its images in $\Hom_R(\Mbar_i,\Tuad(\BDM))$ for all~$i$. As the vertical arrows along the right edge respect the direct sum decomposition, the desired statement for $M$ follows from the corresponding statement for all~$M_i$. 

For the rest of the proof we may thus assume that $\Mbar=\chi^{-1}(M)$ and that it is a free $R$-module of rank~$1$.
Choose a generator $\mbar$ and let $\vartheta$ denote the composite $\Bado$-linear map
$$\xymatrix@R=0pt@C-5pt{\Tuad(\BDM) 
\ar[r]^-\sim & \Hom_R\bigl(\Mbar,\Tuad(\BDM)\bigr)
 \ar@{^{ (}->}[r]
 & \Hom_A\bigl(M,\Tuad(\psi)\bigr)
\ar@{->>}[r] &
\Hom_A\bigl(M,\Tup(\psi)\bigr)\\
\ t\ \ar@{|->}[r] & \ (x\mbar\mapsto xt) & & \\}$$
with varying $x\in R$. Then 
Proposition \ref{MainImageDiag} implies that $\rho_{M,\Fp}(J_K)$ is a submodule of finite index in $\mathop{\rm im}(\vartheta)$. Moreover, for any $\Fq\in Q^\circ_\Fp$ the $B_\Fq$-primary part of $\Tuad(\BDM)$ is the direct factor $\Tuq(\BDM)$. The $B_\Fq$-primary part of $\rho_{M,\Fp}(J_K)$ is therefore a submodule of finite index in $\vartheta(\Tuq(\BDM))$. So we must prove that $\rank_{B_\Fq}(\vartheta(\Tuq(\BDM)))=d$. Since ${\rank_{B_\Fq} (\Tuq(\BDM)) =d}$, this amounts to showing that $\vartheta|\Tuq(\BDM)$ is injective.

To see this let $\pi_\Fp \colon \Tuad(\BDM) \isoto \Tuad(\psi) \onto \Tup(\psi)$ denote the isomorphism from Proposition \ref{TateDoubleIso} composed with the projection to~$\Fp$. 
Since $\Mbar$ is the left $R$-module generated by~$\mbar$, the map $\vartheta$ is given by the formula $t\mapsto (\chi(x\mbar) \mapsto \pi_\Fp(xt))$ with $x\in R$ arbitrary.
For any element $t\in\ker(\vartheta)$ 
and any $y\in R$, we thus have 
$$\vartheta(yt)(\chi(x\mbar))\ =\ \pi_\Fp(xyt)\ =\ \vartheta(t)(\chi(xy\mbar))\ =\ 0$$
for all $x\in R$ 
and hence $yt\in\ker(\vartheta)$. Therefore $\ker(\vartheta)$ is an $R$-submodule. 

Since $\vartheta$ is $\Bado$-linear, it follows that $\ker(\vartheta|\Tuq(\BDM))$ is an $R\otimes_BB_\Fq$-submodule of $\Tuq(\BDM)$. 
But since $\Fq\not=\Fqres$, we know that $R_\Fq := R\otimes_B B_\Fq$ is isomorphic to the matrix ring $\Mat_{d\times d}(B_\Fq)$ and $\Tuq(\BDM)$ is isomorphic to the standard module $B_\Fq^{\oplus d}$. This implies that any non-zero $R_\Fq$-submodule of $\Tuq(\BDM)$ has finite index. 

Thus if $\vartheta|\Tuq(\BDM)$ is not injective, its image is finite, 
and then the $B_\Fq$-primary part of $\rho_{M,\Fp}(J_K)$ is finite. Proposition \ref{MammagePsiMMPrime} then implies that the $B_\Fq$-primary part of $\rho_{M',\Fp}(J_K)$ is finite for the $A$-submodule $M' := A\chi(\mbar) \subset M$. But this contradicts Proposition \ref{NPQOpen} for $M'$ im place of~$M$. 
Therefore $\vartheta|\Tuq(\BDM)$ is injective, and we are done.
\end{Proof}

\begin{Thm}\label{NPQBound}
For any $\Fp$ and any $\Fq\in Q^\circ_\Fp$ the number $n_{\Fp,\Fq}$ for $\rho_{M,\Fp}(J_K)$ satisfies
$$d_{\Fp,\Fq}\cdot\rank_R(\Mbar)\ \le\ n_{\Fp,\Fq}\ \le\ d\cdot\rank_R(\Mbar).$$
\end{Thm}

\begin{Proof}
Choose $R$-linearly independent elements $\mbar_i\in\chi^{-1}(M)$ for $1\le i\le n:=\rank_R(\Mbar)$ and consider the $A$-submodule $M' := \bigoplus_{i=1}^n A\chi(\mbar_i) \subset M$. Then as in Proposition \ref{MammagePsiMMPrime} there is a natural surjection $\rho_{M,\Fp}(J_K) \onto \rho_{M',\Fp}(J_K)$. Since $M'$ satisfies the assumption of Proposition \ref{NPQOpen}, the inequality on the left follows. On the other hand with $M'' := \chi(\Mbar)$ we also have a natural surjection $\rho_{M'',\Fp}(J_K) \onto \rho_{M,\Fp}(J_K)$. Since $M''$ satisfies the assumption of Proposition \ref{BQFree}, the inequality on the right follows.
\end{Proof}


\medskip
The following examples illustrate some ways in which the numbers $n_{\Fp,\Fq}$ can be positioned in the interval allowed by Theorem \ref{NPQBound}.

\begin{Ex}\label{RankDependsOnFpEx}
\rm Suppose that $\psi=\bar\phi$ and has rank~$1$. Then 
$\bar\phi(A)$ contains $\BDM(B)$, so we can identify $B$ with a subring of~$A$ and get $\BDM=\bar\phi|B$. As we independently know that $\rank(\BDM)=d$, it follows that $d=[A/B]$. On the other hand we now get $C=A$. Thus for every $\Fp\not=\Fpres$, the set $Q^\circ_\Fp$ consists of the unique prime $\Fq$ of $B$ below~$\Fp$ and we have $\Bop=B_\Fq$. By Proposition \ref{BCircFpT} the $B_\Fq$-primary part of $\Tup(\psi)$ is therefore free of rank $d_{\Fp,\Fq} = [A_\Fp/B_\Fq]$ over~$B_\Fq$. 

Suppose in addition that $\rank_A(M)=1$, so that $n_{\Fp,\Fq} = d_{\Fp,\Fq} = [A_\Fp/B_\Fq]$ by Proposition \ref{NPQOpen}. As this local degree can vary with~$\Fp$, we conclude that the number $n_{\Fp,\Fq}$ and thus the rank of $\rho_{M,\Fp}(J_K)$ as a $\Bop$-module is in general not independent of~$\Fp$.
\end{Ex}


\begin{Ex}\label{OpennessDependsOnFpEx}
\rm In the situation of Example \ref{RankDependsOnFpEx} let us instead assume that $M=\Mbar$ and $\rank_R(\Mbar)=1$, so that $n_{\Fp,\Fq} = d$ for all $\Fp|\Fq\not=\Fqres$ by Proposition \ref{BQFree}. Then the equations $[A/B]=d$ and $\rank_B(R)=d^2$ imply that $R$ and hence $M$ is a projective left $A$-module of rank~$d$. Thus $\Hom_A(M,\Tup(\psi))$ is a free $B_\Fq$-module of rank~$d\cdot d_{\Fp,\Fq} = d\cdot[A_\Fp/B_\Fq]$. Since $\rho_{M,\Fp}(J_K)$ is a free $B_\Fq$-submodule of rank $d$ thereof, we find that the ranks coincide if and only if $[A_\Fp/B_\Fq]=1$. If $A$ is separable over~$B$ and $[A/B]>1$, this happens for infinitely many $\Fp$ and fails for infinitely many others. We conclude that $\rho_{M,\Fp}(J_K)$ is open in $\Hom_A(M,\Tup(\psi))$ for infinitely many~$\Fp$, but not for infinitely many others. (Compare also Example \ref{NotAadSubmodEx}.)
\end{Ex}


\begin{Ex}\label{NumbersDependOnFqEx}
\rm Begin with a Drinfeld $C$-module $\tilde\phi$ of rank $1$ over~$k$, for which the ring extension $B\subset\tilde\phi(C)$ is separable of degree $e>1$. Choose distinct primes $\Fq,\Fq'\not=\Fqres$ of $B$ and primes $\Fr|\Fq$ and $\Fr'|\Fq'$ of $C$ whose local degrees $[C_\Fr/B_\Fq]$ and $[C_{\Fr'}/B_{\Fq'}]$ are different. Choose an element $t\in C$ whose prime divisors are precisely $\Fr$ and~$\Fr'$, which is possible because the class number of $C$ is finite. Then $\Fr$ and $\Fr'$ are precisely the primes of $C$ over the prime $\Fp := (t)$ of the subring $A := \BFp[t]$. For the Drinfeld $A$-module $\psi:=\bar\phi := \tilde\phi|A$ we then have $Q^\circ_\Fp=\{\Fq,\Fq'\}$, and the relevant numbers from Proposition \ref{BCircFpMT} are $d_{\Fp,\Fq}=d/e\cdot[C_\Fr/B_\Fq]$ and $d_{\Fp,\Fq'}=d/e\cdot[C_{\Fr'}/B_{\Fq'}]$. If now $\rank_A(M)=1$, Proposition \ref{NPQOpen} implies that the $B_\Fq$- and $B_{\Fq'}$-primary parts of $\rho_{M,\Fp}(J_K)$ have different ranks. In this case $\rho_{M,\Fr}(J_K)$ is not a free module over  $\Bop=B_\Fq\times B_{\Fq'}$.
\end{Ex}

\medskip
Now we look at the image of the ramification filtration. First, Proposition \ref{BCircFpMT} directly implies:

\begin{Prop}\label{ZeroOrInfinite}
For any $\Fp$ and $i\ge0$ the image $\rho_{M,\Fp}(J_K^i)$ is either zero or infinite.
\end{Prop}

Next we can view the number
\UseTheoremCounterForNextEquation
\begin{equation}\label{ConductorDef}
\Ff_M\ :=\ \max(S_\Mbar\cup\{0\})
\end{equation}
as the \emph{conductor of $M$ over~$K$}. Since $S_\Mbar$ consists of positive integers, this definition ensures that $\Ff_M=0$ if and only if $M=0$. Otherwise by Theorem \ref{RammagePsi} \ref{RammagePsiRank} we have $\rho_M(J_K^i)=0$ if and only if $i>\Ff_M$, so the conductor measures the depth of the ramification filtration induced on $\rho_M(J_K)$. It turns out that the conductor is already determined by the local representation $\rho_{M,\Fp}$ for every $\Fp$ except possibly for $\Fp=\Fpres$:

\begin{Thm}\label{PConductor}
For any $\Fp$ with $\Tu_\Fp(\psi)\not=0$ and any $i>0$, we have ${\rho_{M,\Fp}(J_K^i)=0}$ if and only if $i>\Ff_M$.
\end{Thm}


\begin{Proof}
By Theorem \ref{RammagePsi} \ref{RammagePsiRank} we have $\rho_M(J_K^i)=0$ and therefore $\rho_{M,\Fp}(J_K^i)=0$ whenever $i>\Ff_M$. The converse holds trivially if $M=0$, so assume that $j:=\Ff_M>0$. It then suffices to show that $\rho_{M,\Fp}(J_K^j) \ne 0$.

For this observe that by the definition (\ref{SMDef}) of $S_\Mbar$ we have $\Mbar\not\subset \Wj{\KOKperf}$. Since $\Wj{\KOKperf}$ is an $R$-submodule of $\KOKperf$ and $\Mbar$ is the $R$-submodule generated by $\chi^{-1}(M)$, there exists an element $m\in M$ with $\chi^{-1}(m)\not\in \Wj{\KOKperf}$. Let $M'$ be the $A$-submodule of $M$ that is generated by~$m$. Then the associated $R$-submodule $\Mbar{}' := R \chi^{-1}(m)$ has the unique break~$j$. Thus Theorem \ref{RammagePsi} \ref{RammagePsiNot} implies that $\rho_{M',\Fp}(J_K^j)$ has finite index in $\rho_{M',\Fp}(J_K)$. At the same time $\rho_{M',\Fp}(J_K)$ has finite index in $\Hom_A(M',\Tup(\psi)) \cong \Tup(\psi)$ by Proposition \ref{MammagePsiOpenP}. As we assumed that $\Tup(\psi)\ne 0$, it follows that $\rho_{M',\Fp}(J^j_K)$ is infinite. Finally, as in Proposition \ref{MammagePsiMMPrime} we have a natural surjection $\rho_M(J_K^j) \onto \rho_{M'}(J_K^j)$ and hence also a natural surjection $\rho_{M,\Fp}(J_K^j) \onto \rho_{M',\Fp}(J_K^j)$. Thus $\rho_{M,\Fp}(J_K^j)$ is infinite, as desired.
\end{Proof}

\begin{Rem}\label{PConductorRem}
\rm
Since the Drinfeld module $\phi$ is the analytic quotient of $\psi$ by~$M$, one can view $\Ff_M$ equally as the {local conductor of $\phi$ over~$K$}. Theorem \ref{PConductor} implies that this agrees with the definition in \cite[\S4.1]{MornevLM}.
\end{Rem}

It is natural to ask whether the images $\rho_{M,\Fp}(J_K^i)$
also determine the set of breaks~$S_\Mbar$. The following example shows that this is in general not~the~case.

\begin{Ex}\label{RammageFpCounterEx}
\rm As in Example \ref{NotAadSubmodEx} let $\psi$ be the Drinfeld $A$-module of rank $1$ with $A=\BFp[t]$ and $\psi_t=\tau$. Then $\bar\phi=\psi$; hence $\chi$ is the identity on $\KOKperf$. Assume that $|k|=p^2$, so that $B=\BFp[s]$ with $\BDM_s=\tau^2$. Identifying $s=t^2$ then makes $B$ a subring of $A$ such that $\BDM=\psi|B$. Choose
an element $\alpha\in k^\times$ with $\alpha^p+\alpha=0$. Then conjugation by $\alpha$ within $R$ induces the non-trivial Galois automorphism of $A$ over~$B$. 

Take integers $j>\ell>0$ that are not divisible by~$p$, and as before let $\xi_j$ and $\xi_\ell$ be elements of $K$ of respective normalized valuations $-j$ and~$-\ell$. Let $M$ be the free $A$-submodule of rank $2$ of $(\KOKperf,\psi)$ that is generated by $\xi_j$ and $\alpha\xi_j+\xi_\ell$. Then $\Mbar$ is the free $R$-submodule of rank $2$ of $\KOKperf$ that is generated by $\xi_j$ and~$\xi_\ell$. By construction we obtain the following commutative diagram:
$$\def\scriptmapsto{{\begin{turn}{-90}$\scriptstyle\mapsto$\end{turn}}}
\vcenter{\vskip-3pt\xymatrix@C=30pt@R=22pt{
&& \text{\large\strut}\Hom_R\bigl(\Mbar,\:\Tuad(\BDM)\bigr) 
\ar[d]^-\wr \ar[rrr]^-\sim_-{h\kern1pt\mapsto (h(\xi_j),\,h(\xi_\ell))\phantom{+\alpha\xi_j}} 
&&&\ \Tuad(\BDM) \rlap{${}^{\oplus2}$} \text{\Large\mathstrut}
\ar[d]^-\wr_{{\raisebox{-3pt}{$\scriptstyle(x,\,y)$}\atop\scriptmapsto\kern1pt}\atop\scriptstyle(x,\,\alpha x+y)} \\
J_K\  \ar[urr]!<6pt,17pt>^-{\rho_\Mbar} \ar[drr]!<6pt,-17pt>^-{\rho_M} 
&& \text{\large\strut}\Hom_A\bigl(M,\:\Tuad(\BDM)\bigr) 
\ar[rrr]^-\sim_-{h\kern1pt\mapsto (h(\xi_j),\,h(\alpha\xi_j+\xi_\ell))} \ar@{=}[d]^\wr
&&&\ \Tuad(\BDM) \rlap{${}^{\oplus2}$} \ar@{=}[d]^\wr \\
&& \text{\large\strut}\Hom_A\bigl(M,\:\Tuad(\psi)\bigr) 
\ar[rrr]^-\sim_-{h\kern1pt\mapsto (h(\xi_j),\,h(\alpha\xi_j+\xi_\ell))} 
&&&\ \Tuad(\psi) \rlap{${}^{\oplus2}$.} \\}
\vskip-3pt}$$
Consider an arbitrary integer $i\ge0$. Then by Theorem \ref{GrJKIsom} \ref{GrJKIsomIsom} the image of $J_K^i$ under the two arrows along the upper edge is 
$$\scriptstyle\left\{\displaystyle\begin{array}{ll}
\Tuad(\BDM)^{\oplus2} & \hbox{if \ $i\le\ell$,}\\[3pt]
\Tuad(\BDM)\oplus\{0\} & \hbox{if \ $\ell<i\le j$,}\\[3pt]
\{0\} & \hbox{if \ $i>j$.}
\end{array}\right.$$
For any prime $\Fp$ of $A$ it follows that $\rho_{M,\Fp}|J_K^i$ is surjective for $i\le\ell$ and zero for $i>j$. 

So suppose now that $\ell<i\le j$. Then the image of $J_K^i$ in the module at the middle of the right edge is the graph of the map $\Tuad(\BDM)\to \Tuad(\BDM)$, $x\mapsto \alpha x$. For any prime $\Fp$ of $A$ that is inert over the corresponding prime $\Fq$ of~$B$, it follows that $\rho_{M,\Fp}(J_K^i)$ corresponds to the graph of multiplication by $\alpha$ on $\Tu_\Fq(\BDM) \cong \Tu_\Fp(\psi)$. Thus $\rho_{M,\Fp}|J_K^i$ is neither zero nor surjective, and so the set $S_\Mbar=\{j,\ell\}$ is determined by $\rho_{M,\Fp}$.
By contrast, for any prime $\Fp$ of $A$ that is split over the corresponding prime $\Fq$ of~$B$, we have $\Fq A=\Fp\Fp'$ with $\Fp\not=\Fp'$ and $\Tu_\Fq(\BDM) \cong \Tu_\Fp(\psi) \oplus \Tu_{\Fp'}(\psi)$. 
Here the choice of $\alpha$ implies that multiplication by $\alpha$ interchanges the two summands. The graph of multiplication by $\alpha$ within $\Tu_\Fq(\BDM)^{\oplus2}$ therefore surjects to $\Tu_\Fp(\psi)^{\oplus2}$.
In this case it follows that $\rho_{M,\Fp}|J_K^i$ is surjective for all $i\le j$, and so $\ell$ cannot be determined from the images $\rho_{M,\Fp}(J_K^i)$.
\end{Ex}

\medskip
Finally, we exploit the dual properties of the local Kummer pairing. 
As a preparation we need:

\begin{Lem}\label{NOSLem}
For any prime $\Fp\not=\Fpres$ of $A$ the embedding $R\into R^\circ$ induces an isomorphism
$$A_\Fp\otimes_AR \isoto A_\Fp\otimes_AR^\circ.$$
\end{Lem}

\begin{Proof}
Since $A_\Fp$ is a flat $A$-module, the short exact sequence $0\to R\to R^\circ\to R^\circ\!/R\to 0$ induces a short exact sequence $0\to A_\Fp\otimes_AR\to A_\Fp\otimes_A R^\circ\to A_\Fp\otimes_AR^\circ\!/R\to 0$. It therefore suffices to show that the third term vanishes. 
For this recall that $\Fpres$ is the kernel of the homomorphism $A\to R/R\tau \cong k$ induced by~$\bar\phi$. The assumption $\Fp\not=\Fpres$ thus implies that $A_\Fp\otimes_AR/R\tau=0$. Since $R\tau^{-n}\!/R$ is isomorphic to an extension of $n$ copies of $R/R\tau$ for any $n\ge0$, it follows that $A_\Fp\otimes_AR\tau^{-n}\!/R=0$ as well. Letting $n\to\infty$ this shows that 
$A_\Fp\otimes_AR^\circ\!/R=0$, as desired.
\end{Proof}

\begin{Thm}\label{SON}
For any prime $\Fp\not=\Fpres$ of $A$, the local Kummer pairing of $\psi$ induces an injective homomorphism
$$A_\Fp\otimes_A(\KOKperf,\,\psi)\ \longinto\ \Hom\bigl(J_K,\:T_\Fp(\psi)\bigr).$$
\end{Thm}

\begin{Proof}
Taking $\Fp$-primary parts on both sides in Theorem \ref{AltPerfectPsi}, the local Kummer pairing of $\psi$ induces an isomorphism of topological $A_\Fp$-modules
\UseTheoremCounterForNextEquation
\begin{equation}\label{AltPerfectPsiKOKFp}
A_\Fp \complot_A (\KOKperf,\,\psi)\ \longisoto\ \Hom_\Bado^\cont\bigl(J_K,\:T_\Fp(\psi)\bigr)
\end{equation}
where $\complot_A$ denotes the completed tensor product. Thus it suffices to show that the natural map
\UseTheoremCounterForNextEquation
\begin{equation}\label{NOSSON}
A_\Fp \otimes_A (\KOKperf,\,\psi) \longto A_\Fp \complot_A (\KOKperf,\,\psi)
\end{equation}
is injective. 

Here the left hand side is isomorphic to $A_\Fp \otimes_A (\KOKperf,\,\bar\phi)$ by Proposition \ref{XIsom}. Since $\KOKperf$ is a free $R^\circ$-module of countably infinite rank by Proposition \ref{KOKKOKperfBasis} \ref{KKbB2}, this tensor product is isomorphic to a countably infinite direct sum of copies of $A_\Fp\otimes_AR^\circ$. By Lemma \ref{NOSLem} it is further isomorphic to a countably infinite direct sum of copies of $A_\Fp\otimes_AR$. Here $R$ is a finitely generated torsion free $A$-module via~$\bar\phi$ and therefore locally free of finite rank. Thus $A_\Fp\otimes_AR$ is a free $A_\Fp$-module of finite rank. Together this shows that $A_\Fp \otimes_A (\KOKperf,\psi)$ is a free $A_\Fp$-module of countably infinite rank. It therefore embeds into its $\Fp$-adic completion, proving the desired injectivity of~\eqref{NOSSON}.
\end{Proof}

\medskip
Now observe that by Propositions  \ref{TateMod0Prop} and \ref{TateMod2Prop} the $\Fp$-primary part of the short exact sequence from Proposition \ref{TateMod3Prop} reads
\UseTheoremCounterForNextEquation
\begin{equation}\label{TateMod3FpProp}
\xymatrix{ 0 \ar[r] & T_\Fp(\psi) \ar[r] &
T_\Fp(\phi) \ar[r] & A_\Fp\otimes_AM \ar[r] & 0.\\}
\end{equation}
The following result determines $\rank(\psi)$ from the Galois representation on $T_\Fp(\phi)$ alone, in analogy to the criterion of N\'eron--Ogg--Shafarevich for the reduction of an abelian variety.

\begin{Thm}\label{NOS}
For every prime $\Fp\not=\Fpres$ of $A$, the image of 
$T_\Fp(\psi)$ is precisely the submodule of $I_K$-invariants in $T_\Fp(\phi)$.
\end{Thm}

\begin{Proof}
Since $A_\Fp$ is a flat $A$-module, the inclusion $M\subset(\KOKperf,\psi)$ induces an injection $A_\Fp\otimes_AM \into A_\Fp \otimes_A (\KOKperf,\psi)$. By Theorem \ref{SON} the local Kummer pairing of $\psi$ therefore induces an injective homomorphism
$$A_\Fp\otimes_AM\ \longinto\ \Hom\bigl(J_K,\:T_\Fp(\psi)\bigr).$$
The theorem is now a direct consequence of Proposition \ref{KummerM}.
\end{Proof}

\begin{Rem}\label{NOSRem}
\rm In general Theorem \ref{NOS} does not hold for $\Fp=\Fpres$. For instance, it fails trivially if $M\not=0$ and $T_\Fpres(\psi)=0$. In that case it also fails for $\psi|A'$ in place of $\psi$ for any admissible coefficient ring $A'\subset A$. Thus it can fail even when the reduction $\bar\phi$ does not have $\Fpres$-rank zero.
\end{Rem}

\begin{Rem}\label{Gardeyn}\rm
A dual version of Theorem \ref{NOS} can be deduced from the theory of Gardeyn \cite{GardeynT},
stating that the quotient module of $I_K$-coinvariants of $T_\Fp(\phi)$ is
precisely $A_\Fp\otimes_A M$.
%
\end{Rem}

\end{document}